\documentclass[3p]{elsarticle}
\usepackage{amssymb,amsmath,amsthm}
\usepackage{color}
\usepackage{url}

\newcommand{\Rset}{\mathbb{R}}

\newcommand{\Nset}{\mathbb{N}}

\usepackage{bm}

\newcommand{\FF}{\mathbf{F}}

\newcommand{\ii}{\mathbf{i}}

\newcommand{\kk}{\mathbf{k}}

\newcommand{\nn}{\mathbf{n}}
\newcommand{\pp}{\mathbf{p}}
\newcommand{\PP}{\mathbf{P}}

\newcommand{\xx}{\mathbf{x}}

\newcommand{\aalpha}{\bm{\alpha}}
\newcommand{\bbeta}{\bm{\beta}}
\newcommand{\xxi}{\bm{\xi}}

\newcommand{\mcH}{\mathcal{H}_{mix}}

\newcommand{\norm}[2]{\|#1\|_{#2}}

\newtheorem{prob}{Problem}

\newtheorem{example}{Example}

\def\ConvSize{0.39}
\def\OptSetSize{0.29}
\def\GammaSize{0.39}
\def\SpaceSize{3pt}

\newcommand{\La}[1]{#1}

\begin{document}

\begin{frontmatter}
\title{A sparse-grid isogeometric solver}

\author[kaust]{Joakim Beck}
\ead{joakim.beck@kaust.edu.sa}

\author[unipv,imati]{Giancarlo Sangalli}
\ead{giancarlo.sangalli@unipv.it}

\author[unipv,imati]{Lorenzo Tamellini\corref{cor1}}
\ead{tamellini@imati.cnr.it}

\cortext[cor1]{Corresponding author}
\address[kaust]{CEMSE, King Abdullah University of Science and Technology (KAUST), Thuwal 23955-6900, Saudi Arabia}
\address[unipv]{Dipartimento di Matematica ``F. Casorati'', Universit\`a di Pavia, Via Ferrata 5, 27100, Pavia, Italy}
\address[imati]{Consiglio Nazionale delle Ricerche - Istituto di Matematica Applicata e Tecnologie Informatiche
``E. Magenes'' (CNR-IMATI), Via Ferrata 1, 27100, Pavia, Italy}
\begin{abstract}
Isogeometric Analysis (IGA) typically adopts tensor-product splines and NURBS
as a basis for the approximation of the solution of PDEs.
In this work, we investigate to which extent IGA solvers can benefit from
the so-called sparse-grids construction in its combination technique form, which was first introduced in the early 90's 
in the context of the approximation of high-dimensional PDEs. 

The tests that we report show that, in accordance to the literature, a sparse-grid construction can 
indeed be useful if the solution of the PDE at hand is sufficiently smooth.
Sparse grids can also be useful in the case of non-smooth solutions
when some a-priori knowledge on the location of the singularities of the solution
can be exploited to devise suitable non-equispaced meshes. Finally, we remark
that sparse grids can be seen as a simple way to parallelize pre-existing serial IGA solvers
in a straightforward fashion, which can be beneficial in many practical situations.

\medskip
\textbf{Highlights} \vspace{-5pt}
\begin{itemize}
\item A sparse grid version of classical isogeometric solvers is proposed \vspace{-8pt}
\item The proposed methodology is competitive with standard isogeometric
solvers if the solution of the PDE is sufficiently smooth \vspace{-8pt}
\item The proposed methodology can reuse pre-existing isogeometric solvers almost out-the-box and can be thought as a simple way to  
  parallelize a serial solver \vspace{-8pt}
\item Radical meshes in the parametric domain can be adopted as a remedy to improve 
  sparse grid convergence for solutions without sufficient regularity on domains with corners. \vspace{-8pt}
\end{itemize}

\end{abstract}



\end{frontmatter}

\section{Introduction}

Isogeometric analysis (IGA), which was introduced by Hughes et al. \cite{Hughes:2005,IGA-book} in 2005, consists in 
solving numerically a PDE by approximating its solution with
B-splines, Non-Uniform Rational B-Splines (NURBS), and extensions, i.e., 
with the same basis employed to parametrize the computational geometry with CAD softwares.
The method has attracted considerable attention in the engineering computing community,
not only because it could simplify the meshing process 
but also because of other interesting features, including a larger flexibility 
in the choice of the polynomial degree and regularity of the basis used to approximate the solution, 
and a more effective error vs. degrees-of-freedom ratio with respect to standard finite element
methods; see \cite{acta-IGA} and references therein.

In this paper, we focus on computational cost efficiency, and in particular 
on the fact that $d$-dimensional splines are generated
by tensorization of univariate splines, which means that the computational work
typically increases exponentially with the dimension $d$ of the problem.
Although this phenomenon is somehow acceptable for $d=2,3$, it
could nevertheless lead to excessive computational costs both in the formation
of the Galerkin matrix and in the solution of the corresponding linear system.

While a possible way out would be to embrace the tensor structure of the 
construction and exploit it to optimize as much as possible quadrature 
and linear solvers, see, e.g., \La{\cite{Gao2013,Gao2014,Gao2015,
sangalli2016isogeometric,calabro2016fast,mantzaflaris2016low,
antolin2015efficient,donatelli2017symbol,hofreither2015robust}}, 
here we consider an alternative approach and  explore the 
possibility of moving from a tensor grid to a so-called sparse grid, 
together with the combination technique.

Sparse grids emerged in the 90's as a technique to solve generic $d$-dimensional PDEs, see 
\cite{zenger91sparse,bungartz:PhD.thesis,griebel:first,Griebel.schneider.zenger:combination,b.griebel:acta,Garcke2013}, 
and have shown to be quite effective in mitigating the dependence of the computational cost on $d$,
if the solution features certain regularity assumptions, 
more demanding than the usual Sobolev ones.
More precisely, denoting by $h$ the finest mesh-size
of a discretization, a classic sparse-grid construction
with piecewise linear polynomials only needs $\mathcal{O}(h^{-1} \log(h^{-1})^{d-1})$ degrees-of-freedom
to yield an $\mathcal{O}(h^{-2} \log(h^{-1})^{d-1})$ $L^2$ approximation error,
as opposed to a standard tensor-based approximation method which would need
$\mathcal{O}(h^{-d})$ degrees-of-freedom (note the exponential dependence on $d$)
for an $\mathcal{O}(h^{-2})$ $L^2$ approximation error, i.e. 
a much larger number of degrees-of-freedom for almost the same
convergence rate. Roughly speaking, the sparse-grid technique consists in rewriting 
a tensor approximation of a function (in this work, built with tensorized splines)
as a sum of hierarchical tensor contributions, 
and then observing that under suitable regularity assumptions the more expensive contributions 
(i.e., those corresponding to fine discretizations along all directions simultaneously) 
are actually negligible and can thus be safely neglected without
compromising too much the accuracy of the approximation.
Upon some rearrangment of the sparse-grid formula,
the final form that will be used in practice (the so-called ``combination technique'')
will consist of a linear combination of a certain number of 
``small'' tensor approximations, which will recover a close-to-full order approximation.
It is important to remark that each of these tensor approximations is just a standard
PDE solve (IGA in this work), which can be computed independently: therefore, sparse grids
can also be seen as a straightforward way of parallelizing a legacy serial IGA solver.
Note in particular that the communcation between computing cores is reduced to a minimum 
if one is not interested in pointwise value of the solution everywhere in the domain but rather in 
linear functionals of the solution.

The main goal of this work is to discuss both how IGA can take advantage of a sparse-grid construction,
and vice-versa, how IGA can help broadening the field of application of the sparse-grid technique.
Indeed, IGA solvers would provide a natural and systematic way of extending the sparse-grid approach,
specifically by allowing use of basis functions with arbitrary degree and regularity, 
and general non-square domains. Note that extensions of sparse grids to these features were already proposed
in literature between the mid 1990s and the early 2000s and are mentioned in the 2004 sparse-grid survey
paper \cite{b.griebel:acta}: curvilinear domains generated
by transfinite interpolation where proposed by 
\cite{Dornseifer:curvilinear.sparse.grids,bungartz:general-deg-and-dom},
while high-order hierarchicial Lagrange polynomials were analyzed
in \cite{bungartz:general-deg-and-dom,achatz:high.order.sparse.grids}.
Other alternatives for high-order polynomial mentioned in \cite[chap. 4.5]{b.griebel:acta}
are bicubic Hermite polynomials or higher-order finite-differences/interpolets.
However, these works did not have a significant follow-up, perhaps due to 
the fact that they required advanced tools to generate non-standard polynomial bases, while
splines-based IGA solvers are nowadays widely available 
and can be easily used in a combination technique formulation in a black-box fashion.

In order to improve convergence for problems with corner and edge singularities, 
that do not possess the Sobolev regularity that is required for optimal convergence
of sparse grids, we combine sparse grids with non-uniform radical meshes,
see \cite{Babuska_book,Beirao_Cho_Sangalli}. 
This idea was briefly touched in previous sparse grids literature in \cite{garcke:graded,griebel.thurner:graded}.


\vspace{10pt}

The rest of this work is organized as follows: Section \ref{section:method} introduces the 
test case considered throughout the rest of the paper and its discretization via IGA. 
Section \ref{sec:sparse-grids} explains the details 
of the sparse-grid construction and recalls the related convergence results;
observe that these are derived in the case of a square-cubic domain, 
but hold also in the tests in Section \ref{sec:numerical-res} 
where non-square domains are considered. Section \ref{sec:numerical-res}
will also report the performance of the sparse grids in terms of accuracy
versus computational time and number of degrees-of-freedom,
over a few different geometries, and discuss
the potential advantages of the sparse-grid method over a standard tensorized
discretization. Concluding remarks and perspectives on future 
works are wrapped up in Section \ref{sec:conclusions}.


\section{Method}\label{section:method}

\subsection{Problem setting}

Let $\Omega \subset \Rset^d$, $d=2,3$ be a compact set.
In this paper, we focus on an elementary problem, i.e. the Poisson equation:
\begin{prob}\label{prob:poisson}
Find $u:\Omega \rightarrow \Rset$ such that 
\begin{equation}\label{eq:poisson}
  \begin{cases}
    -\Delta u(\xx) = f(\xx)	& \mbox{ in } \Omega \\
    u(\xx) = 0				& \mbox{ on } \partial \Omega.
  \end{cases}
\end{equation}
\end{prob}
A sufficient condition for $u \in H^k(\Omega)$, $k \geq 1$, is that $f \in H^{k-2}(\Omega)$ and $\partial \Omega \in C^k(\Rset^{d-1})$.
For non-smooth domains, for instance domains with corners, see e.g. \cite{Grisvard}.

\subsection{The isogeometric method for solving PDEs}

We now briefly recall here the fundamentals of IGA and refer to \cite{Hughes:2005,acta-IGA,BeiraoDaVeiga:2014,DeBoor} 
for a more thorough discussion.
Given a so-called parametric interval, $\hat{I}$ (typically, $\hat{I}=[0,1]$),
we introduce a knot vector, i.e, a non-decreasing vector $\Xi=[\xi_{1},\xi_{2}...,\xi_{n+p+1}]$,
with $n,p \in \Nset$, and $\xi_1, \xi_{n+p+1}$ coinciding with the extrema of $\hat{I}$; 
each $\xi_i$ is a knot and an interval $(\xi_i, \xi_{i+1})$ having non-zero length is an
element; let us further denote by number of elements as $N_{el}$. 
Observe that the elements need not have the same length:
if that is the case, we call such length mesh-size, and denote it by $h$.
A knot vector is said to be open if its first and last
knot have multiplicity  $p+1$. Observe that also 
internal nodes could have multiplicity greater than one; we define the non-decreasing
vector $Z= [\zeta_1, \ldots, \zeta_{N_{el}+1}]$ as the vector of knots of $\Xi$ without
repetitions, and let $m_i$ the multiplicity of $\zeta_i$ in $\Xi$, so that
$\sum_{i=1}^{N_{el}} m_i = n+p+1$. 

Given a knot vector $\Xi$, we define the B-splines by means of the Cox-De Boor 
recursive formula, for $i=1,\ldots,n$:
\begin{align}
& \widehat{B}_{i,0}(\xi)= 
	\begin{cases}
      1 &  \xi_{i}\leq \xi<\xi_{i+1}\\
      0 &  \textrm{otherwise,} \\
    \end{cases} 
    \label{eq:Bsp}   \\ 
& \widehat{B}_{i,p}(\xi)=
	\begin{cases}
      \dfrac{\xi-\xi_{i}}{\xi_{i+p}-\xi_{i}}\widehat{B}_{i,p-1}(\xi)+\dfrac{\xi_{i+p+1}-\xi}{\xi_{i+p+1}-\xi_{i+1}}\widehat{B}_{i+1,p-1}(\xi),
      & \xi_{i}\leq \xi<\xi_{i+p+1} \\
      0, & \textrm{otherwise,}
\end{cases} \nonumber
\end{align}
where we adopt the convention $0/0=0$; note that the basis corresponding to an open knot vector
will be interpolatory in the first and last knot. The B-splines functions are polynomials of degree $p$
and continuity $C^{p-m_i}$ at $\zeta_i$ and they form a basis of the space of splines,
\[
S_p(\Xi,\hat{I}) = \text{span}\left\{ \widehat{B}_{i,p}\,,\,\, i=1,\ldots,n \right\}.
\]
For $d=2$ we define the parametric domain $\widehat{\Omega} = \hat{I} \times \hat{I}$
(extension to the case $d=3$ is trivial). We consider two
open knot vectors $\Xi_1, \Xi_2$ with $n_1 + p_1 + 1$ and $n_2+p_2+1$ knots respectively,
the corresponding knots without repetitions $Z_1,Z_2$, and the
tensor products $\mathbf{\Xi} = \Xi_1 \otimes \Xi_2$, $\mathbf{Z}=Z_1 \otimes Z_2$; in particular, $\mathbf{Z}$ generates 
a cartesian mesh over $\widehat{\Omega}$ 
composed by $N_{el,1} \times N_{el,2}$ rectangular elements.
Taking tensor products of the univariate B-splines over $\Xi_1$ and $\Xi_2$ we obtain a basis
for the space of bi-variate splines, 
\[
S_\pp(\mathbf{\Xi},\hat{\Omega}) = \text{span}\{\widehat{B}_{\ii,\pp}, \ii \leq \nn \},
\]
where $\ii=[i_1,i_2]$, $\pp=[p_1,p_2]$, $\nn=[n_1,n_2]$, 
$\ii \leq \nn \Leftrightarrow i_1 \leq n_1, i_2 \leq n_2$,
and $\widehat{B}_{\ii,\pp}(\xi_1,\xi_2) = \widehat{B}_{i_1,p_1}(\xi_1) \widehat{B}_{i_2,p_2}(\xi_2)$. 

We assume for the sake of simplicity that the computational domain $\Omega$ 
can then be parameterized by a linear combination of B-splines with given control points $\PP_\ii \in \Rset^2$, 
\begin{equation}\label{eq:F-B-splines}
\xx \in \Omega \Leftrightarrow \xx = \FF(\xxi) = \sum_{\ii \leq \nn} \PP_{\ii} \widehat{B}_{\ii,\pp} (\xxi) \text{ for some }\xxi \in \widehat{\Omega},  
\end{equation}
where $\FF:\widehat{\Omega} \rightarrow \Omega$ is the parameterization of the geometry $\Omega$.
For the sake of exposition, we focus in this section on B-splines, but everything can be
repeated verbatim for NURBS.

According to the isogeometric principle, the spline basis is also used to approximate the solution of the 
Poisson equation. To this end, we introduce the spline space on the physical domain $\Omega$ as follows:
\begin{equation}\label{eq:spline_space_phys}
  S_{\pp}(\mathbf{\Xi},\Omega) = \text{span} \{B_{\ii,\pp} = \widehat{B}_{\ii,\pp} \circ \FF^{-1}, \ii \leq \nn  \}
\end{equation}
and then approximate the solution of Problem \ref{prob:poisson} as 
\begin{equation}\label{eq:u-IGA-approx}
  u(\xx) \approx u_\nn(\xx)  = \sum_{\ii \leq \nn} c_{\ii} B_{\ii,\pp}(\xx).  
\end{equation}
Much like with the finite element spaces, 
the isogeometric approximation of $u$ in \eqref{eq:u-IGA-approx}, $u_{\nn}$, converges to $u$ as the cardinality
of $S_\pp(\mathbf{\Xi},\Omega)$ increases. 
In this work, we focus on a $h$-refinement version of the sparse-grid technique, but 
it could be possible to devise a $p$-sparsification technique as well (or a combined $h-p$ version).
In the numerical experiments in Section \ref{sec:numerical-res}, we will consider 
isotropic degrees, i.e., $p_1=\ldots,p_d=p$.
For a fixed degree $p$, we will consider both 
$C^0$ and $C^{p-1}$ versions of the $h$-refinement (the $C^{p-1}$ case is referred to as
``$k$-refinement'' or ``$k$-method'' in the isogeometric literature).


\section{Sparse grids}\label{sec:sparse-grids}

The basic building block for a sparse-grid construction is a sequence of nested and increasingly accurate univariate 
approximation operators, which will be then first ``hiercharchized'' (i.e., an equivalent sequence of hierarchical
operators will be derived from the initial sequence) and then tensorized.

The sequence of univariate approximation spaces/operators of a generic function $f:\hat{I}\rightarrow \Rset$, 
indexed by $\alpha \in \Nset \cup \{0\}$, is constructed as follows. We first introduce 
the trivial open knot vector over the reference interval $\hat{I}$, $\Xi_0=[0,\ldots,0,1,\ldots,1]$, 
and generate the $\alpha$-th open knot vector $\Xi_\alpha$ 
by dyadic subdivision of $\Xi_0$, so that $\Xi_\alpha$ has $2^\alpha$ elements.
For notational convenience, we then denote by the short-hand $S_\alpha$
the space of B-splines/NURBS built over $\Xi_\alpha$ with fixed degree $p=1,2,3$,
i.e., $S_\alpha = S_p(\Xi_\alpha,\hat{I})$,
and $f_\alpha$ be a suitable approximation of $f$ in $S_\alpha$ (later, a Galerkin approximation). 
Finally, the univariate ``hierarchization'' of the sequence of approximation is obtained by defining the operators
\[
\Delta_\alpha (f) = f_\alpha - f_{\alpha-1},
\]
where we adopt the convention $f_{-1} = 0$. Observe that the following telescopic equality holds:
\[
f_\alpha = \sum_{k=0}^\alpha \Delta_k (f).
\]

As for the second step (the ``tensorization'') we again illustrate it in $d=2$, for ease of notation,
but the construction can be applied with straightforward modifications to any $d$. 
Thus, let now $f:\hat{I}_1 \otimes \hat{I}_2 = \widehat{\Omega}\rightarrow \Rset$ and consider
a multi-index with non-negative components $\aalpha= [\alpha_1,\, \alpha_2] \in (\Nset \cup \{0\})^2$. Furthermore,
let $\mathbf{\Xi}_{\aalpha} = \Xi_{\alpha_1} \times \Xi_{\alpha_2}$, and the corresponding 
B-splines space be $S_{\aalpha}$.
Extending the univariate notation, we now let $f_{\aalpha}$ be
the approximation of $f$ in $S_{\aalpha}$, and in particular
we denote by $f_{[\alpha_1, \infty]}, f_{[\infty, \alpha_2]}$ the semi-discrete approximations
of $f$ along directions $\xi_1$ and $\xi_2$. 
We then consider hierarchical operators along each direction individually, i.e., let
\[
\Delta^1_{\alpha_1} f = f_{[\alpha_1,\infty]} - f_{[\alpha_1-1,\infty]},
\]
be the difference of two consecutive semi-discrete approximations of $f$ along direction $\xi_1$, 
and analogously 
\[
\Delta^2_{\alpha_2} f = f_{[\infty,\alpha_2]} - f_{[\infty,\alpha_2-1]},
\]
and we introduce the multidimesional hierachical operator as the 
tensor product of the univariate $\Delta_\alpha$ operators, i.e., 
\begin{align}
  \Delta_{\aalpha}(f) 
  & = \Delta^1_{\alpha_1} [\Delta^2_{\alpha_2} (f)] \label{eq:Delta-op} \\
  & = \Delta^1_{\alpha_1} [f_{[\infty,\alpha_2]} - f_{[\infty,\alpha_2-1]}] \nonumber \\
  & = f_{[\alpha_1,\alpha_2]} - f_{[\alpha_1-1,\alpha_2]} - f_{[\alpha_1,\alpha_2-1]} + f_{[\alpha_1-1,\alpha_2-1]}, \label{eq:Delta-combi-tec}       
\end{align}
where the last form is usually referred to as the combination technique.
Observe that by telescopy 
\begin{equation}\label{eq:TP-telescopic}
  f_{\aalpha} = \sum_{\bbeta \leq \aalpha} \Delta_{\bbeta}(f).         
\end{equation}

Upon introducing this hierarchical decomposition of the approximation of $f$, the crucial observation 
is that for $f$ regular enough the norm of each component $\Delta_{\bbeta}(f)$
can be expected to be decreasing with respect to $|\bbeta| = \sum_{\ell=1}^d \beta_\ell$. 
If that is the case, many of the terms in \eqref{eq:TP-telescopic}
can be dismissed from the approximation without compromising too much its accuracy. At the same time however,
due to the dyadic subdivision approach used to generate the univariate spaces $S_\alpha$ to be tensorized,
it can be easily seen from \eqref{eq:Delta-combi-tec} that the number of degrees-of-freedom of $\Delta_{\bbeta}$ 
grows exponentially in $|\bbeta| = \sum_{\ell=1}^d\beta_\ell$. Thus, introducing 
\begin{equation}\label{eq:TD-approx}
  f_{J} = \sum_{|\bbeta| \leq J} \Delta_{\bbeta}(f), \quad \text{for some } J \in \Nset \cup \{0\},       
\end{equation}
we are left with an approximation which gives up moderately on the accuracy while, in fact, significantly reducing the number
of degrees-of-freedom. 
Observe that in principle, such a computation can be already performed by considering the sparse-grid approximation in the form of 
\eqref{eq:TD-approx}; however, this requires using discretizations of hierarchical type, which may not be straightforward.
Instead, by further expanding each $\Delta_{\bbeta}$ in \eqref{eq:TD-approx}
as in \eqref{eq:Delta-combi-tec}, we obtain the following expression of the sparse-grid approximation 
\La{(see \cite[Lemma 1]{wasi.wozniak:cost.bounds} for the derivation of the expression)} which is more
amenable to computation, because it only employs standard solvers:
\begin{align}\label{eq:TD-combi-tec}
  f_{J} 
  & = \sum_{\kk \in \{0,1\}^d, |\bbeta + \kk| \leq J}(-1)^{|\kk|} f_{\bbeta}\\
  & = \sum_{J-d+1 \leq |\bbeta| \leq J} (-1)^{J-|\bbeta|} \binom{d-1}{J-|\bbeta|} f_{\bbeta} \\[6pt]
  & = 
  \begin{cases}
    \displaystyle \sum_{|\bbeta|=J} f_{\bbeta} -   \sum_{|\bbeta|=J-1} f_{\bbeta} & \text{if } d=2, \\[12pt]
    \displaystyle \sum_{|\bbeta|=J} f_{\bbeta} - 2 \sum_{|\bbeta|=J-1} f_{\bbeta} + \sum_{|\bbeta|=J-2} f_{\bbeta} & \text{if } d=3,
  \end{cases}
\end{align}
where the derivation of the second equality can be found in, e.g., \cite{wasi.wozniak:cost.bounds}.
We refer to \eqref{eq:TD-combi-tec} as the combination technique form of $f_J$, see, e.g.,
\cite{Griebel.schneider.zenger:combination,Bungartz.Griebel.Roschke.ea:pointwise.conv,Hegland:combination,
griebel.harbrecht:combi-conv,Garcke2013,reisinger:analysis,reisinger:efficient}.
We remark that the number of components to be computed in the combination technique 
\eqref{eq:TD-combi-tec} is substantially smaller than the number of hierarchical components in \eqref{eq:TD-approx}:
indeed, the set in \eqref{eq:TD-approx} is a discrete simplex while in \eqref{eq:TD-combi-tec} only
the uppermost $d$ layers of such simplex need to be computed. More precisely:
\begin{itemize}
\item for $d=2$ one has
$\mathrm{card}\big(\{|\bbeta| \leq J\}\big) = \frac{J(J+1)}{2} = \mathcal{O}(J^2)$ while 
$\mathrm{card}\big(\{|\bbeta| = J\}\big) + \mathrm{card}\big(\{|\bbeta| = J-1\}\big) = 2J-1 =  \mathcal{O}(J)$;
\item for $d=3$ one has
$\mathrm{card}\big(\{|\bbeta| \leq J\}\big) = \frac{J(J+1)(J+2)}{6} =  \mathcal{O}(J^3)$ 
while $\mathrm{card}\big(\{|\bbeta| = J\}\big) + \mathrm{card}\big(\{|\bbeta| = J-1\}\big) + \mathrm{card}\big(\{|\bbeta| = J-2\}\big) = 
\frac{2}{3}J(J-1)+1 =  \mathcal{O}(J^2)$. 
\end{itemize}
In practice, computing $f_J$ entails 
solving a number of ``moderate'' tensor approximations of $f$ where the $d$ directions
are never simultaneously refined to full accuracy, compared to the full tensor approximation $f_{\aalpha}$
(see also Example \ref{ex:combi} detailed next). 
The quantity $h$ appearing in the estimates below is defined as $h = 2^{-J}$
and is proportional to the finest mesh-size employed in the sparse-grid discretization
(which, we recall again, is never reached ``simultaneously'' along each direction).
A combination technique formula uses $\mathcal{O}(J^{d-1} 2^{J}) = \mathcal{O}(h^{-1} |\log_2 h|^{d-1})$ 
degrees of freedom, see, e.g., \cite{bungartz:general-deg-and-dom,achatz:high.order.sparse.grids},  
compared to the corresponding fully resolved grid that employs 
$\mathcal{O}(2^{Jd}) = \mathcal{O}(h^{-d})$ degrees of freedom. Crucially, note that the number of degrees of freedom of the sparse
grid depends in a milder way on $d$.
\begin{figure}
  \centering
  \includegraphics[width=0.24\linewidth]{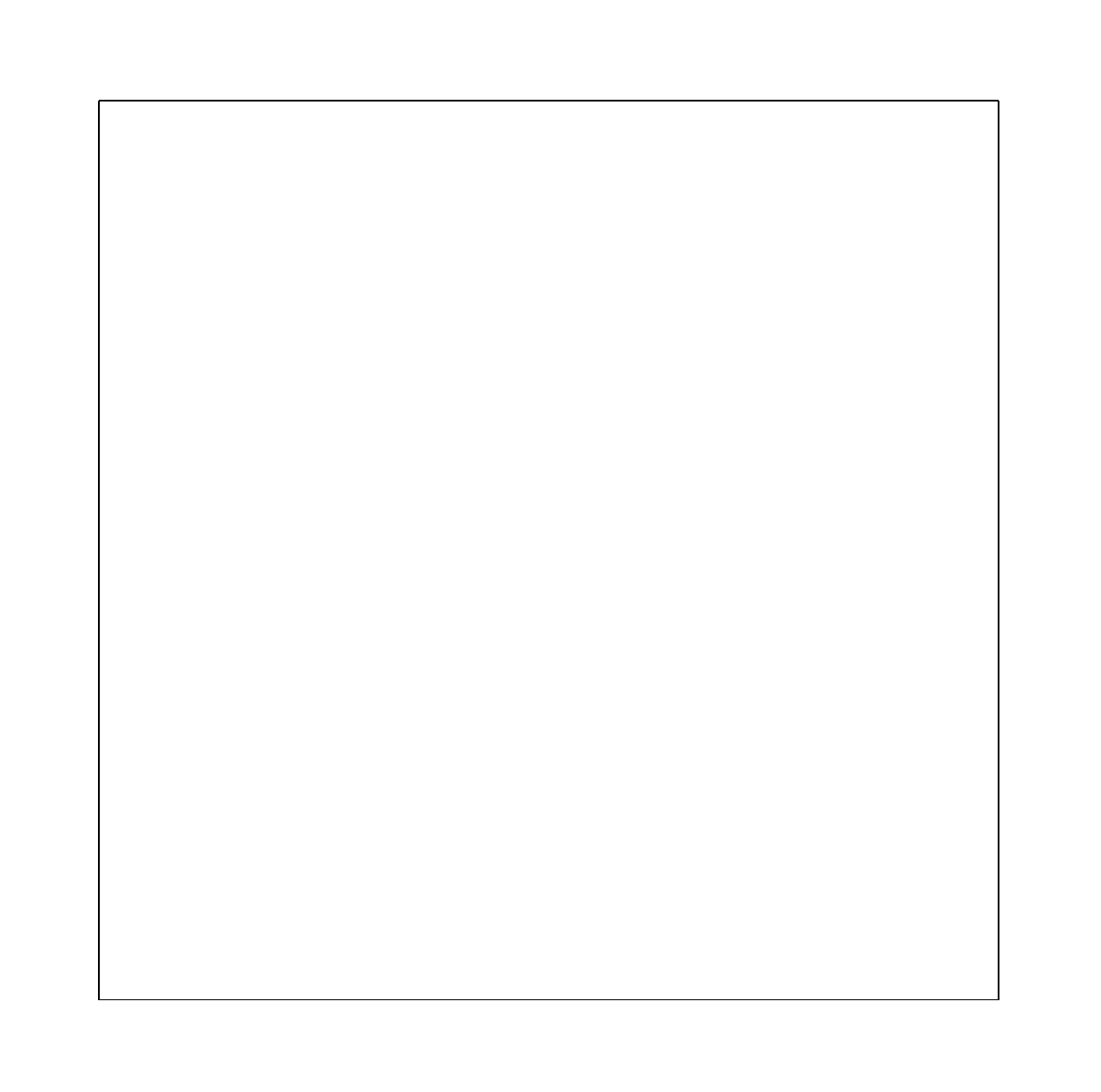}
  \includegraphics[width=0.24\linewidth]{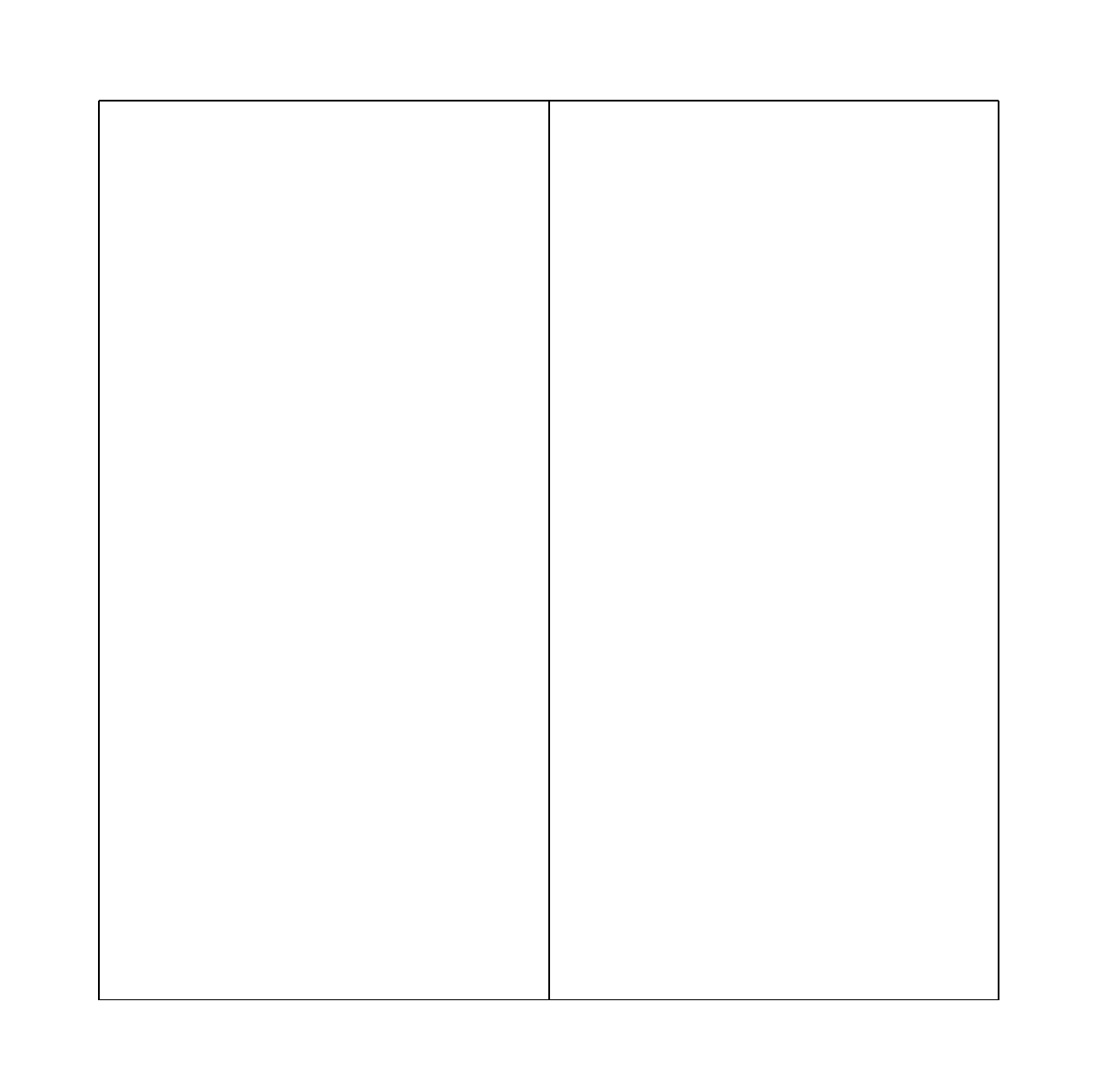}
  \includegraphics[width=0.24\linewidth]{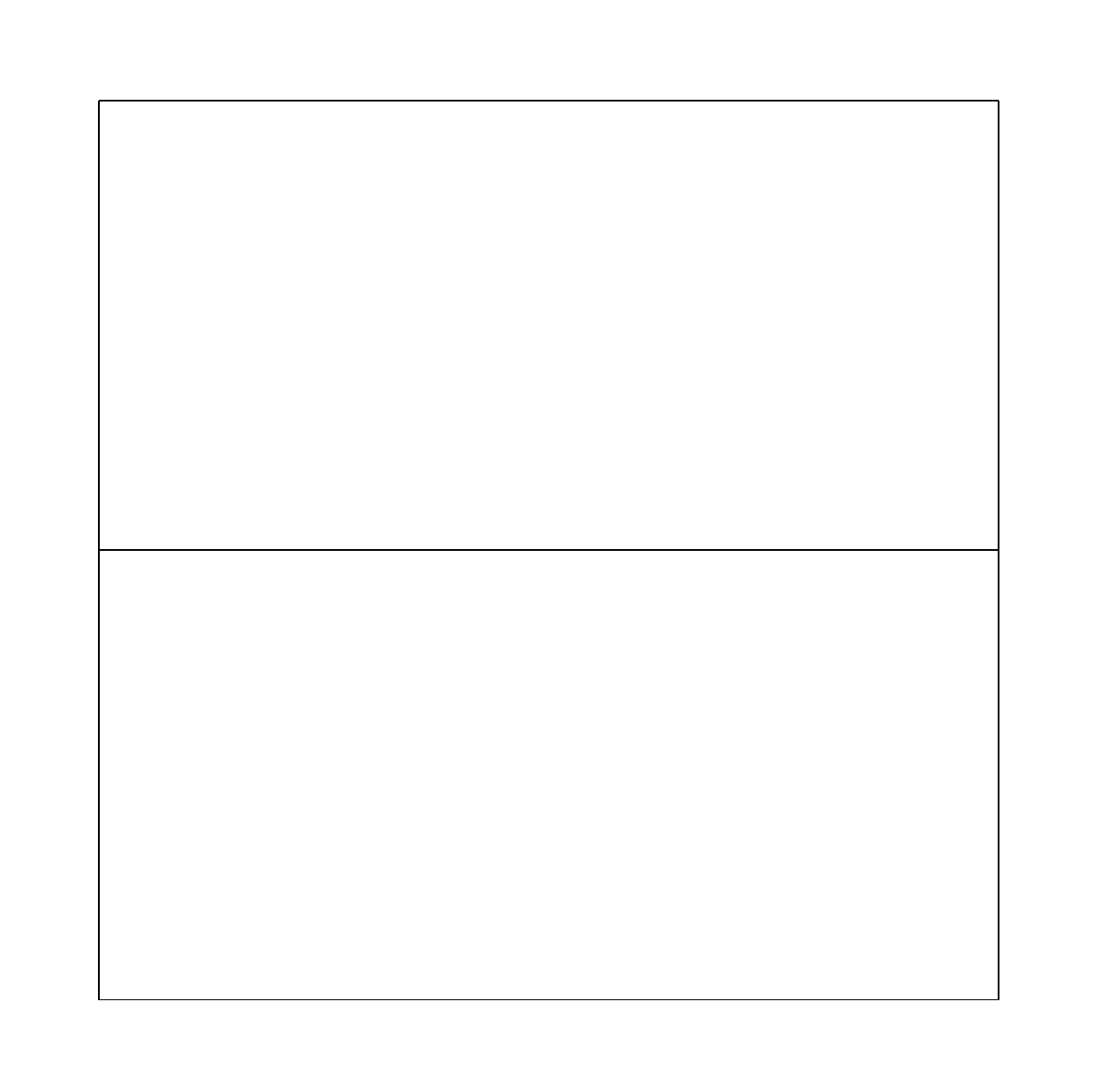}
  \includegraphics[width=0.24\linewidth]{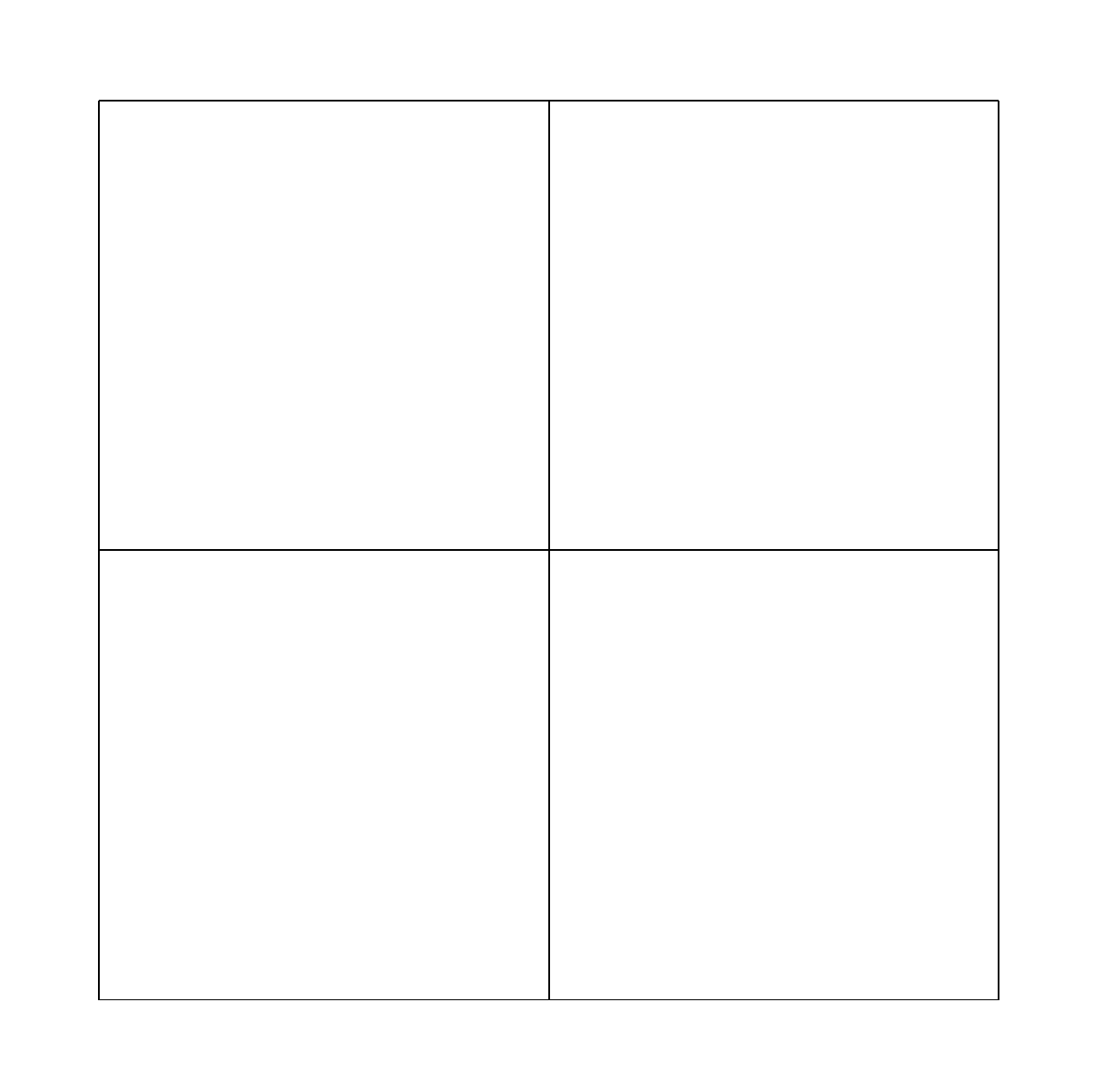}
  \caption{Bezier meshes for Example \ref{ex:combi}. From left to right: $f_{[0\,0]}$, $f_{[1\,0]}$, $f_{[0\,1]}$, $f_{[1\,1]}$.
  The first three are to be linearly combined as in Equation \eqref{eq:combi-ex}
  to obtain the combination technique, to be contrasted with the last 
  full approximation grid.}
  \label{fig:sparse-vs-tensor}
\end{figure}
\begin{example}\label{ex:combi}
In this example we consider $\hat{I}=[0,1]$, the knot-vector $\Xi_0=[0,0,1,1]$,
and $\widehat{\Omega}=[0,1]^2 = \hat{I} \otimes \hat{I}$.
We then consider an approximation of $f:\widehat{\Omega}\rightarrow \Rset$ with a sparse grid using $J=1$, and below we show the
calculus to convert the sparse-grid representation of $f_J$, \eqref{eq:TD-approx}, to 
its combination technique equivalent, \eqref{eq:TD-combi-tec}. 
  \begin{align}
    f_{J} = \sum_{\beta_1 + \beta_2 \leq J=1} \Delta_{\bbeta}(f)
    & = \Delta_{[1,0]}(f) + \Delta_{[0,1]}(f) + \Delta_{[0,0]}(f) \nonumber \\[-10pt]
    &= f_{[1,0]} - f_{[0,0]} 					- \underbrace{f_{[1,-1]}}_{0} + \underbrace{f_{[0,-1]}}_{0} \nonumber \\
    &+ f_{[0,1]}	- \underbrace{f_{[-1,1]}}_{0} - f_{[0,0]} 		+ \underbrace{f_{[-1,0]}}_{0} \nonumber \\
    &+ f_{[0,0]} - \underbrace{f_{[-1,0]}}_{0} - \underbrace{f_{[0,-1]}}_{0}	+ \underbrace{f_{[-1,-1]}}_{0} \nonumber\\
    & =f_{[1,0]} +  f_{[0,1]} - f_{[0,0]}. \label{eq:combi-ex}
  \end{align}
Observe also that by adding the operator $\Delta_{[1,1]}(f)$ to the sparse-grid representation above,
we obtain the full tensor approximation $f_{[1,1]}$, in agreement with the telescopic formula \eqref{eq:TP-telescopic}:\La{
  \begin{align}
    f_{J} + \Delta_{[1,1]}(f)
     = \sum_{[\beta_1,\beta_2] \leq [1,1]} \Delta_{\bbeta}(f)
    & = \Delta_{[1,1]}(f) + f_{[1,0]} +  f_{[0,1]} - f_{[0,0]} \nonumber \\[-10pt]
    &= f_{[1,1]} - f_{[0,1]} - f_{[1,0]}+ f_{[0,0]} + f_{[1,0]} +  f_{[0,1]} - f_{[0,0]} \nonumber \\
    &= f_{[1,1]}.
  \end{align}}
The results are illustrated in Figure \ref{fig:sparse-vs-tensor}.
The meshes of the three components of $f_{J=1}$ in \eqref{eq:combi-ex}
are the three left-most meshes, while 
the right-most mesh is the one of the full tensor grid approximation $f_{[1,1]}$.
It is evident how the full resolution $h=1/2$ is reached only in one direction at a time
for the combination technique.
\end{example}

Finally, we discuss the convergence of the combination technique,
following closely \cite{griebel.harbrecht:combi-conv,griebel.harbrecht:tensor-spaces,griebel.harbrecht:L-fold};
see also \cite{Bungartz.Griebel.Roschke.ea:pointwise.conv,reisinger:analysis,
pflaum-zhou:combi_conv,Bungartz:two.proofs}.
Note that, although the sparse-grid representation \eqref{eq:TD-approx} and 
the combination technique \eqref{eq:TD-combi-tec} are formally two different representations
of the same approximation, the strategies to prove convergence of the two techniques are different, 
see, e.g., \cite{b.griebel:acta} and \cite{griebel.harbrecht:combi-conv}. 
To state the convergence results, we need 
to introduce the so-called mixed regularity Sobolev spaces.


First, let $H^{l}(\hat{I})$ for $l \in \Rset_+ \setminus \Nset$ be the fractional order Sobolev space,
extending the definition of a standard Sobolev space $H^{l}(\hat{I})$ for integer $l$.
Then, borrowing notation from \cite{griebel.harbrecht:combi-conv}, 
we define the Sobolev spaces of mixed derivatives on $\widehat{\Omega}=[0,1]^d$ as 
\[
H^{r_1,\ldots,r_d}_{mix}(\widehat{\Omega}) = H^{r_1}(\hat{I}) \otimes \cdots \otimes H^{r_d}(\hat{I}),
\]
and finally the spaces
\[
\mcH^{r_1,\ldots,r_d}(\widehat{\Omega}) = 
H^{r_1+1,r_2,\ldots,r_d}_{mix}(\widehat{\Omega}) \cap 
H^{r_1,r_2+1, r_3 \ldots,r_d}_{mix}(\widehat{\Omega})\cap \ldots \cap
H^{r_1,r_2,\ldots,r_{d-1}, r_d+1}_{mix}(\widehat{\Omega}).
\]
As will be detailed below, the estimate of the $H^1$ norm of the error will be valid
provided that $u$ belongs to certain $\mcH^{r_1,\ldots,r_d}$ spaces, while 
the estimate of the $L^2$ norm of the error will require $H^{r_1,\ldots,r_d}_{mix}$ regularity.
However, it is usually not trivial to perform a sharp mixed regularity analysis of the solution of a PDE.
Thus, in the following we will also state the convergence results
in terms of the standard Sobolev spaces. To this end, observe that the following inclusions trivially hold
\begin{align}
&H^{k}(\widehat{\Omega}) \subset H^{k/d,k/d,\ldots,k/d}_{mix}(\widehat{\Omega})  \label{eq:mixed-stand-incl}\\
&H^{1+k}(\widehat{\Omega}) \subset \mcH^{k/d,k/d,\ldots,k/d}(\widehat{\Omega}).
\nonumber
\end{align}

In case of $C^0$ finite elements (though the case of higher continuity is similar),
the convergence results that we report in the following have been proved for non-tensor operators 
on rectangular domains, i.e., on $\widehat{\Omega}$, 
and are valid under certain additional technical assumptions 
(see \cite{griebel.harbrecht:combi-conv} for details
on the proof of the $H^1$ convergence, while the proof of the $L^2$ convergence 
is detailed in  \cite{griebel.harbrecht:tensor-spaces}).
Given the exploratory nature of this work, we will not try to prove the same kind of result
in generic domains $\Omega = \mathbf{F}(\widehat{\Omega})$, and instead we will content ourselves with
verifying numerically that these convergence rates are valid also in this latter framework.

\paragraph{$H^1$ norm}
For $0<s\leq p$ there holds
\begin{equation}\label{eq:H1-conv-wrt-mix}
  \norm{u-u_J}{H^1(\widehat{\Omega})} \leq h^s J^{(d-1)/2} \norm{u}{\mcH^{s,\ldots,s}(\widehat{\Omega})}.  
\end{equation}
Observe that using the inclusions in \eqref{eq:mixed-stand-incl}, it holds also that 
\begin{equation}\label{eq:H1-conv}
\norm{u-u_J}{H^1(\widehat{\Omega})} \leq h^s J^{(d-1)/2} \norm{u}{H^{1+ds}(\widehat{\Omega})}. 
\end{equation}
Comparing the last equation with the convergence estimate of standard tensor methods
\begin{equation}
  \tag{\ref{eq:H1-conv}-bis}\label{eq:H1-conv-tensor}
  \norm{u-u_{[J,\ldots,J]}}{H^1(\widehat{\Omega})} \leq h^s \norm{u}{H^{1+s}},
\end{equation}
we see in a more quantitative way the fact that sparse grids reach the same convergence rate
than standard tensor methods (up to a logarithmic term) 
by employing less degrees of freedom, at the expense of higher regularity requirements.

\paragraph{$L^2$ norm}
For $0 < s \leq r=p+1$ there holds
\begin{equation}\label{eq:L2-conv-wrt-mix}
  \norm{u-u_J}{L^2(\widehat{\Omega})} \leq h^s J^{(d-1)/2} \norm{u}{H^{s,\ldots,s}_{mix}\widehat{\Omega}},  
\end{equation}
and again from \eqref{eq:mixed-stand-incl},
\begin{equation}\label{eq:L2-conv}
  \norm{u-u_J}{L^2(\widehat{\Omega})} \leq h^s J^{(d-1)/2} \norm{u}{H^{sd}(\widehat{\Omega})}, 
\end{equation}
to be compared to the standard tensor convergence 
\begin{equation}
  \tag{\ref{eq:L2-conv}-bis}\label{eq:L2-conv-tensor}
  \norm{u-u_{[J,\ldots,J]}}{L^2(\widehat{\Omega})} \leq 2^{-J s} \norm{u}{H^{s}},
\end{equation}
leading to same conclusions as in the $H^1$ case.

We summarize the expected convergence rates in Table \ref{tab:sparse-grids-rates}.
Observe that, with a slight abuse, in Table \ref{tab:sparse-grids-rates} and in the rest of this work 
we use the expression ``sparse-grid convergence \emph{rate}'' to indicate $s$, 
i.e. the exponent of $h$ at the right-hand side of Equations 
\eqref{eq:H1-conv-wrt-mix} to \eqref{eq:L2-conv},
neglecting the factor $J^{(d-1)/2}$ which depends logarithmically on $h$.
Moreover, here and in the rest of the manuscript, 
we will denote by ``SG-IGA'' the sparse-grid version of isogeometric analysis.

\begin{table}[tbp]
\centering
{\scriptsize
\begin{tabular}{ccccccc}
  \multicolumn{7}{c}{$ \textbf{$H^1$ error, $d=2$}$} \\ \hline
			& $u \in H^2$ 	& $u \in H^3$ 	& $u \in H^4$ 	& $u \in H^5$ 	& $u \in H^6$ 	& $u \in H^7$	\\
			& ($ u \in \mcH^{1/2,1/2}$) 
 			& ($ u \in \mcH^{1,1}$)	
 			& ($ u \in \mcH^{3/2,3/2}$)	
 			& ($ u \in \mcH^{2,2}$)	
 			& ($ u \in \mcH^{5/2,5/2}$)	
        	& ($ u \in \mcH^{3,3}$) \\ \hline
$p=1$, SG-IGA& 	1/2			&		1		&		1		&		1		&		1		&		1		\\
$p=1$, IGA	& 	1			&		1		&		1		&		1		&		1		&		1		\\ \hline
$p=2$, SG-IGA& 	1/2			&		1		&		3/2		&		2		&		2		&		2		\\
$p=2$, IGA	& 	1			&		2		&		2		&		2		&		2		&		2		\\ \hline
$p=3$, SG-IGA& 	1/2			&		1		&		3/2		&		2		&		5/2		&		3		\\
$p=3$, IGA	& 	1			&		2		&		3		&		3		&		3		&		3		\\
\end{tabular}
\vskip 10pt
\begin{tabular}{ccccccc}
  \multicolumn{7}{c}{$ \textbf{$L^2$ error, $d=2$}$} \\ \hline
			& $u \in H^2$ 	& $u \in H^3$ 	& $u \in H^4$ 	& $u \in H^5$ 	& $u \in H^6$ 	& $u \in H^7$	\\
			& ($ u \in H^{1,1}$) 
 			& ($ u \in H^{3/2,3/2}$)	
 			& ($ u \in H^{2,2}$)	
 			& ($ u \in H^{5/2,5/2}$)	
 			& ($ u \in H^{3,3}$)	
        	& ($ u \in H^{7/2,7/2}$) \\ \hline
$p=1$, SG-IGA& 	1			&		3/2		&		2		&		2		&		2		&		2		\\
$p=1$, IGA	& 	2			&		2		&		2		&		2		&		2		&		2		\\ \hline
$p=2$, SG-IGA& 	1			&		3/2		&		2		&		5/2		&		3		&		3		\\
$p=2$, IGA	& 	2			&		3		&		3		&		3		&		3		&		3		\\ \hline
$p=3$, SG-IGA& 	1			&		3/2		&		2		&		5/2		&		3		&		7/2		\\
$p=3$, IGA	& 	2			&		3		&		4		&		4		&		4		&		4		\\
\end{tabular}
\vskip 10pt
\begin{tabular}{cccccccc} 
  \multicolumn{7}{c}{$ \textbf{$H^1$ error, $d=3$}$} \\ \hline
			& $u \in H^2$ 	& $u \in H^3$ 	& $u \in H^4$ 	& $u \in H^5$ 	& $u \in H^6$ 	& $u \in H^7$	\\
			& ($ u \in \mcH^{1/3,1/3,1/3}$)\!\!\!\!\!\!\!\!\! 
 			& ($ u \in \mcH^{2/3,2/3,2/3}$)\!\!\!\!\!\!\!\!\!	
 			& ($ u \in \mcH^{1,1,1}$)\!\!\!\!\!\!\!\!\!	
 			& ($ u \in \mcH^{4/3,4/3,4/3}$)\!\!\!\!\!\!\!\!\!	
 			& ($ u \in \mcH^{5/3,5/3,5/3}$)\!\!\!\!\!\!\!\!\!
        	& ($ u \in \mcH^{2,2,2}$) \\ \hline
$p=1$, SG-IGA& 	1/3			&		2/3		&		1		&		1		&		1		&		1		\\
$p=1$, IGA	& 	1			&		1		&		1		&		1		&		1		&		1		\\ \hline
$p=2$, SG-IGA& 	1/3			&		2/3		&		1		&		4/3		&		5/3		&		2		\\
$p=2$, IGA	& 	1			&		2		&		2		&		2		&		2		&		2		\\ \hline
$p=3$, SG-IGA& 	1/3			&		2/3		&		1		&		4/3		&		5/3		&		2		\\
$p=3$, IGA	& 	1			&		2		&		3		&		3		&		3		&		3		\\
\end{tabular}
\vskip 10pt
\begin{tabular}{ccccccc}
  \multicolumn{7}{c}{$ \textbf{$L^2$ error, $d=3$}$} \\ \hline
			& $u \in H^2$ 	& $u \in H^3$ 	& $u \in H^4$ 	& $u \in H^5$ 	& $u \in H^6$ 	& $u \in H^7$	\\
			& ($ u \in H^{2/3,2/3,2/3}$)\!\!\!\!\!\!\!\!\!
 			& ($ u \in H^{1,1,1}$)\!\!\!\!\!\!\!\!\! 	
 			& ($ u \in H^{4/3,4/3,4/3}$)\!\!\!\!\!\!\!\!\! 	
 			& ($ u \in H^{5/3,5/3,5/3}$)\!\!\!\!\!\!\!\!\! 	
 			& ($ u \in H^{2,2,2}$)\!\!\!\!\!\!\!\!\! 	
        	& ($ u \in H^{7/3,7/3,7/3}$) \\ \hline
$p=1$, SG-IGA& 	2/3			&		1		&		4/3		&		5/3		&		2		&		2		\\
$p=1$, IGA	& 	2			&		2		&		2		&		2		&		2		&		2		\\ \hline
$p=2$, SG-IGA& 	2/3			&		1		&		4/3		&		5/3		&		2		&		7/3		\\
$p=2$, IGA	& 	2			&		3		&		3		&		3		&		3		&		3		\\ \hline
$p=3$, SG-IGA& 	2/3			&		1		&		4/3		&		5/3		&		2		&		7/3		\\
$p=3$, IGA	& 	2			&		3		&		4		&		4		&		4		&		4		\\
\end{tabular}}
\caption{Provable lower bounds on the $H^1$ and $L^2$ convergence rates for sparse-grid (SG-IGA) and tensor approximation (IGA) and $d=2,3$.}\label{tab:sparse-grids-rates}
\end{table}


\section{Numerical tests}\label{sec:numerical-res}

\begin{figure}[tp]
  \centering
  \includegraphics[width=0.25\textwidth]{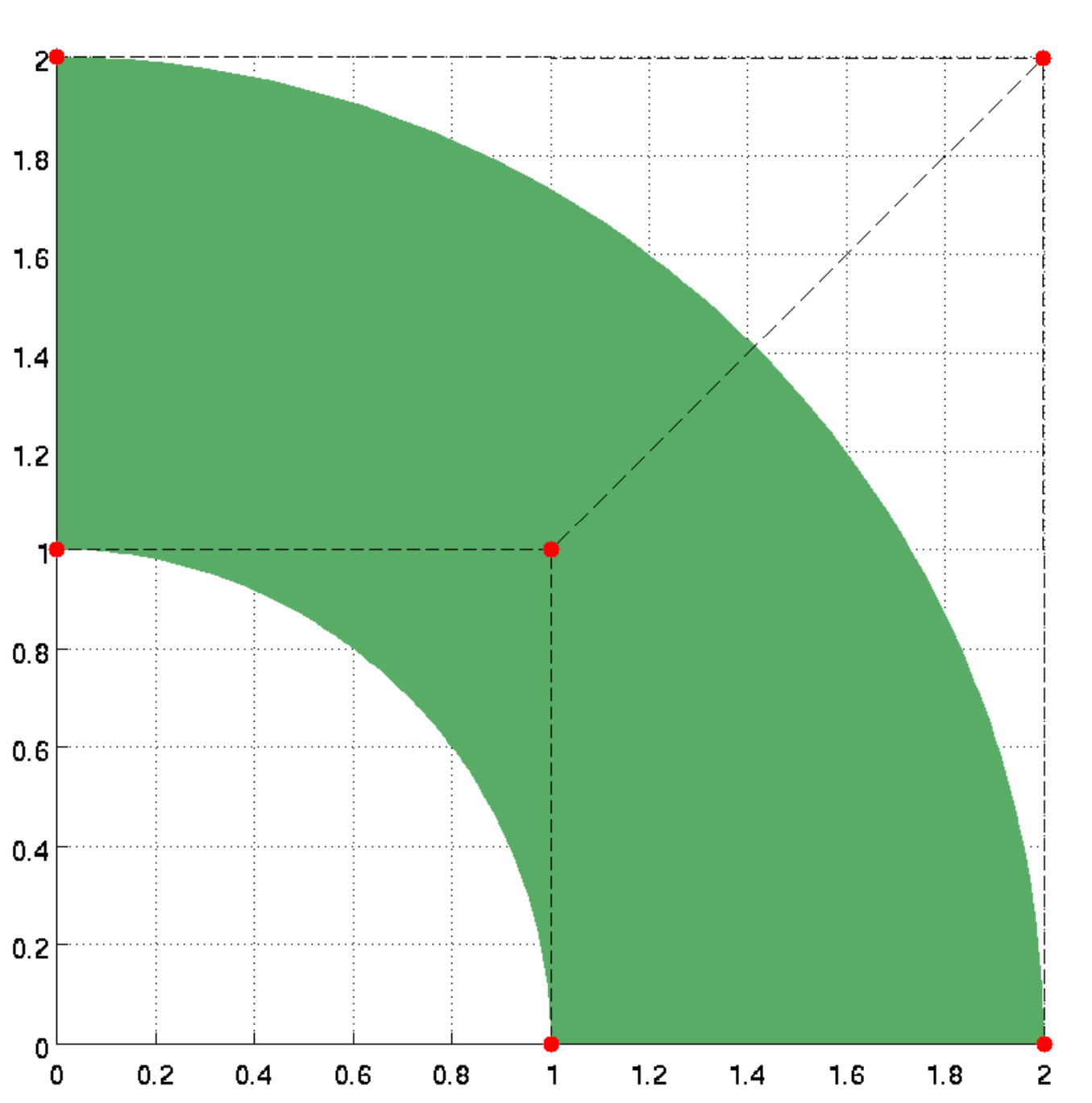}
  \includegraphics[width=0.33\textwidth]{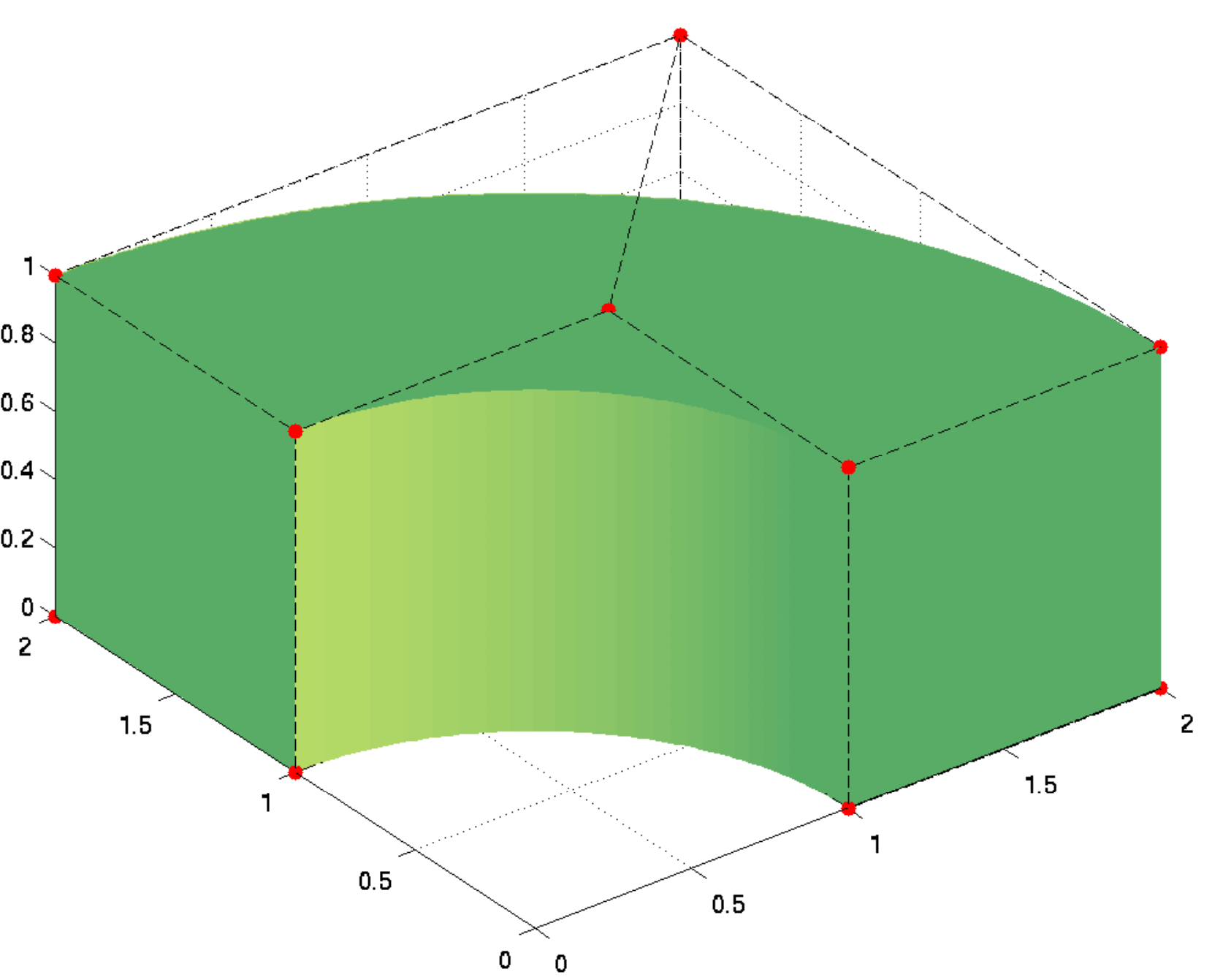}
  \includegraphics[width=0.25\textwidth]{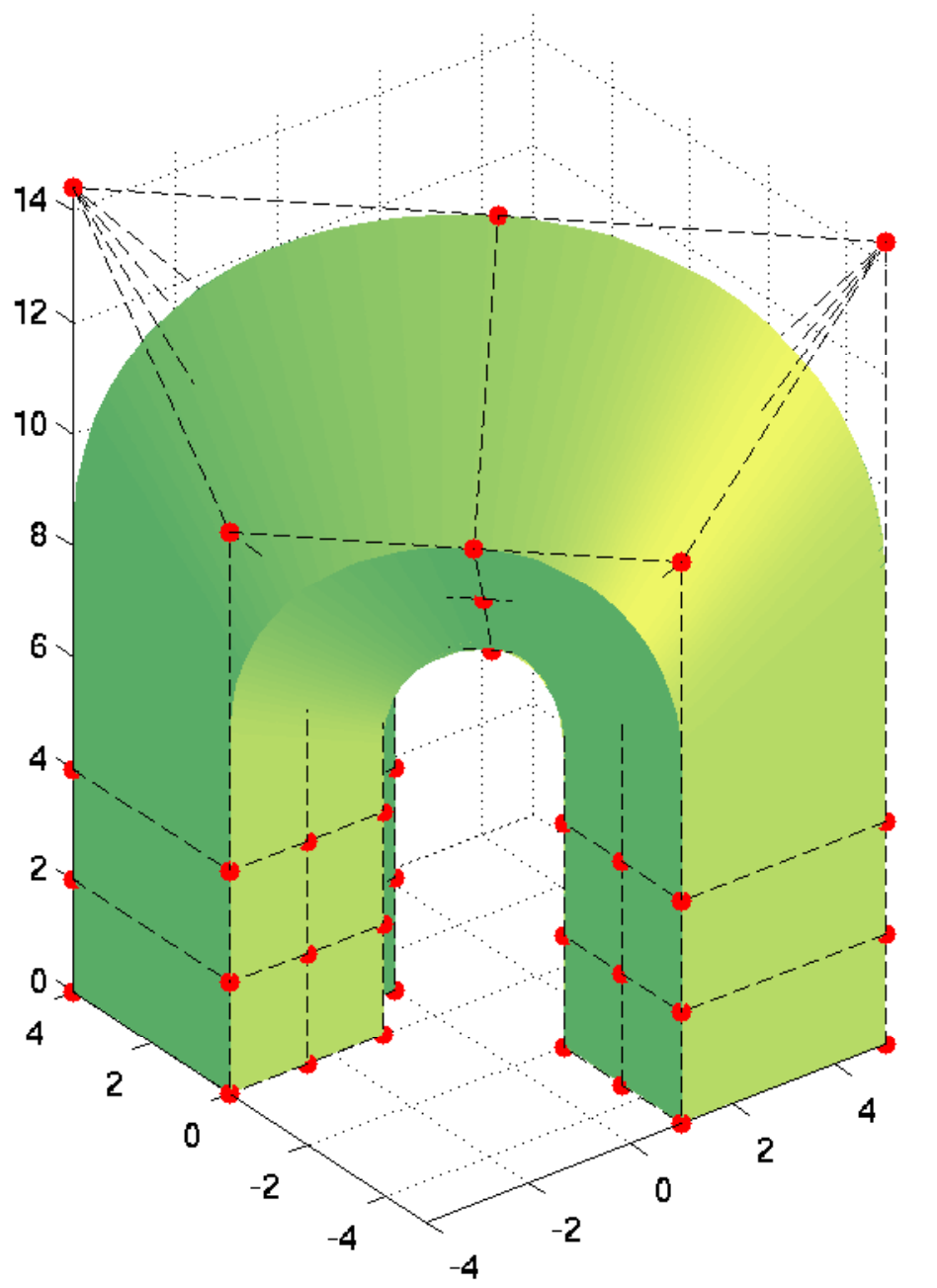}
  \caption{Domains for the numerical tests with control points (in red). 
    Left: Quarter of annulus, $d=2$; center: quarter of annulus $d=3$; right: horseshoe.}
  \label{fig:domains}
\end{figure}

In what follows, we consider a few variations of Problem \ref{prob:poisson}, which
highlight the features of SG-IGA. More specifically, two different geometries are considered, 
the classical quarter of annulus in $d=2,3$ and a horseshoe-shaped domain, see Figure \ref{fig:domains}.
We verify numerically that SG-IGA attains for problems with regular and low-regular solutions the theoretical convergence rates 
presented in Table \ref{tab:sparse-grids-rates}, and compare SG-IGA to standard IGA in terms of computational cost.  
All tests have been performed using the Matlab/Octave IGA package GeoPDEs \cite{geopdesv3},
on a Linux workstation with 12 i7 cores at 3.30 GHz (although all tests reported have been executed in serial unless explicitely mentioned), 64 GB RAM and Matlab 2015a 8.5.

\subsection{A matter of costs}  

\La{In our experiments we adopt GeoPDEs 3.0
\cite{geopdesv3}, which is an isogeometric solver based on a finite
element architecture. Denoting by  $N_{\mathsf{EL}} $ and $N_{\mathsf{DOF}} $ the
number of mesh elements and degrees-of-freedom, respectively, the
standard matrix formation requires $ {O(N_{\mathsf{EL}} p^{3d})}$
FLOPs  for degree $p$ splines, any continuity (that is, $
{O(N_{\mathsf{DOF}} p^{2d})}$ FLOPs  for $C^0$  Bernstein polynomials and  $
{O(N_{\mathsf{DOF}} p^{3d})}$ FLOPs  for $C^{p-1}$ splines). The linear
solver cost for preconditioned iterative solvers is proportional to
$N_{\mathsf{DOF}} $ but robustness with respect to $p$ is difficult or
costly to achieve with standard finite element preconditioners. For
small systems, a direct solver (that we adopt) is faster.
Recent approaches have significantly improved computational efficiency,
especially for $C^{p-1}$ continuity (see
e.g. \cite{calabro2016fast,mantzaflaris2016low,sangalli2016isogeometric,hofreither2015robust,sangalli2017matrix}). }

\La{In this work, we relate the error  both to the computational time (with GeoPDEs
based routines)  and to the number of degrees-of-freedom: 
time measures the ``actual cost'' of our implementation 
(mainly dependent on the element number and polynomial degree, independent on the spline regularity), 
while degrees-of-freedom represent the ``best'' 
possible computational cost in an ideal setting.}

When considering degrees-of-freedom as a measure of the computational cost, 
we remark that the components in the combination technique formula of SG-IGA are
solved separately, and so, we find it of interest to look at both the
largest component's degrees-of-freedom, and the sum of all degrees-of-freedom. 
In fact, when fully exploiting parallelization the component
with the largest degrees-of-freedom is the computational bottleneck in
SG-IGA, as the combination technique is ``embarrassingly parallel,''
see, e.g., \cite{griebel:first,griebel:parallerl92,heene:efficient,Heene:scalable}, 
and one can exploit that existing IGA solvers can be used out-of-the-box.  

The combination technique has two evident drawbacks: first, determining the optimal
balance of its components over the cores at hand is a combinatoric
problem, which is non-trivial to solve (see 
\cite{heene:efficient,Heene:scalable} for more details on an efficient
implementation of a combination technique strategy, 
\La{as well as \cite{garcke_hegland_nielsen_2006} for jointly optimizing both the distribution of componentes over the
available cores and the use of multiple cores, i.e., employing a parallel solver, 
for computing a single component}). 
Second, the approximation of the solution on the entire domain is not
available unless all grids have been combined together, e.g.,
interpolated on a reference grid, which will have in general a
non-negligible computational cost. Such a ``post-process'' step can be
avoided if one is only interested in a linear functional of the
solution.

As for the computational time, we show timings for both serial and parallel execution over $C$ cores.
The time for serial execution will be the true computational time while the 
time for parallel executions will not be the true one but instead an 
``optimized time'': after having clocked the computational time for each tensor
during the serial execution, we sort the tensors needed for each SG-IGA level
in decreasing computational time, and assign them to the $C$ cores at disposal
in this order, in such a way that as soon as a core is free it takes up the largest tensor
not assigned yet.
The rationale is to show the ``ideal'' behaviour of the combination technique by 
factoring out implementation suboptimalities.
Furthermore, the computational time is measured as the sum of setup time 
(meshing and preliminary operations, including matrix and right hand side
assembly) and solve time, i.e., in our results we do not consider the 
time it takes to evaluate the solution on the reference grid used to compute the
error, which in a suboptimal implementation may dominate any other 
cost, for both IGA and SG-IGA. 

\subsection{Quarter of annulus domains, regular solution}

With reference to Equation \eqref{eq:poisson}, we first look at the quarter-of-annulus benchmark,
in its two- and three-dimensional versions. 
In the following, $x,y,z$ will as usual denote the components of the $d$-dimensional vector $\xx \in \Omega$. 
In this first set of experiments, we choose $f(\xx)$ such that 
\begin{align*}
  u(x,y) &= -(x^2+y^2-1)(x^2+y^2-4)xy^2 					& \mbox{ for } d=2, \\              
  u(x,y,z)& = -(x^2+y^2-1)(x^2+y^2-4)xy^2 \sin(\pi z)  	& \mbox{ for } d=3.  
\end{align*}


\begin{table}[tp]
  \centering
\begin{tabular}{l|rr|rr}
						& \multicolumn{2}{c|}{$d=2$}	& \multicolumn{2}{c}{$d=3$}					\\
						& $H^1$ error 	& $L^2$ error 	& 	$H^1$ error & $L^2$ error 	\\ \hline
    $C^0, p=2$, SG-IGA	&	2.29(2)		&	3.03(3)		&	2.38(2)		&	3.12(3)		\\
    $C^0, p=2$, IGA		&	2.10(2)		&	2.96(3)		&	2.05(2)		&	2.97(3)		\\
    $C^1, p=2$, SG-IGA	&	2.16(2)		&	3.04(3)		&	2.23(2)		&	3.08(3)		\\
    $C^1, p=2$, IGA		&	2.10(2)		&	3.00(3)		&	2.03(2)		&	3.08(3)		\\
    $C^0, p=3$, SG-IGA	&	3.21(3)		&	3.94(4)		&	3.30(3)		&	4.00(4)		\\
    $C^0, p=3$, IGA		&	3.08(3)		&	3.98(4)		&	3.06(3)		&	4.01(4)		\\
    $C^1, p=3$, SG-IGA	&	3.05(3)		&	3.93(4)		&	3.13(3)		&	3.93(4)		\\
    $C^1, p=3$, IGA		&	2.91(3)		&	3.98(4)		&	2.89(3)		&	3.85(4)		\\
    $C^2, p=3$, SG-IGA	&	3.06(3)		&	4.15(4)		&	3.15(3)		&	4.26(4) 	\\
    $C^2, p=3$, IGA		&	2.93(3)		&	4.08(4)		&	2.98(3)		&	4.13(4)		\\
\end{tabular}
  \caption{Convergence of SG-IGA and IGA methods with respect to the sparse-grid level, $J$. 
    The number in parenthesis indicates the expected rate from theory, cf. Table \ref{tab:sparse-grids-rates}.} 
  \label{tab:quarter-of-annulus-reg-rates}

\end{table}

We start the discussion by examining the convergence of the SG-IGA and IGA approximation errors in $L^2$ and $H^1$ norm 
with respect to the sparse-grid level $J$, cf. Equations \eqref{eq:H1-conv-wrt-mix} to \eqref{eq:L2-conv}. The measured numerical rates 
(recall the ``informal definition'' of sparse-grid rate we gave while discussing Table \ref{tab:sparse-grids-rates},
at the end of Section \ref{sec:sparse-grids})
are reported in Table \ref{tab:quarter-of-annulus-reg-rates} and should be compared to 
Table \ref{tab:sparse-grids-rates}, which lists the theoretical convergence rates. 
\La{Numerical rates are computed by using a least squares fit (Matlab command \texttt{polyfit})
of the quantities in equations \eqref{eq:H1-conv} and \eqref{eq:L2-conv}}. 
For convenience, we complement Table \ref{tab:quarter-of-annulus-reg-rates} by indicating in parenthesis, 
next to each measured rate, the corresponding theoretical value taken from Table \ref{tab:sparse-grids-rates}. 

\begin{figure}[tp]
  \centering
  \includegraphics[width=\ConvSize\linewidth]{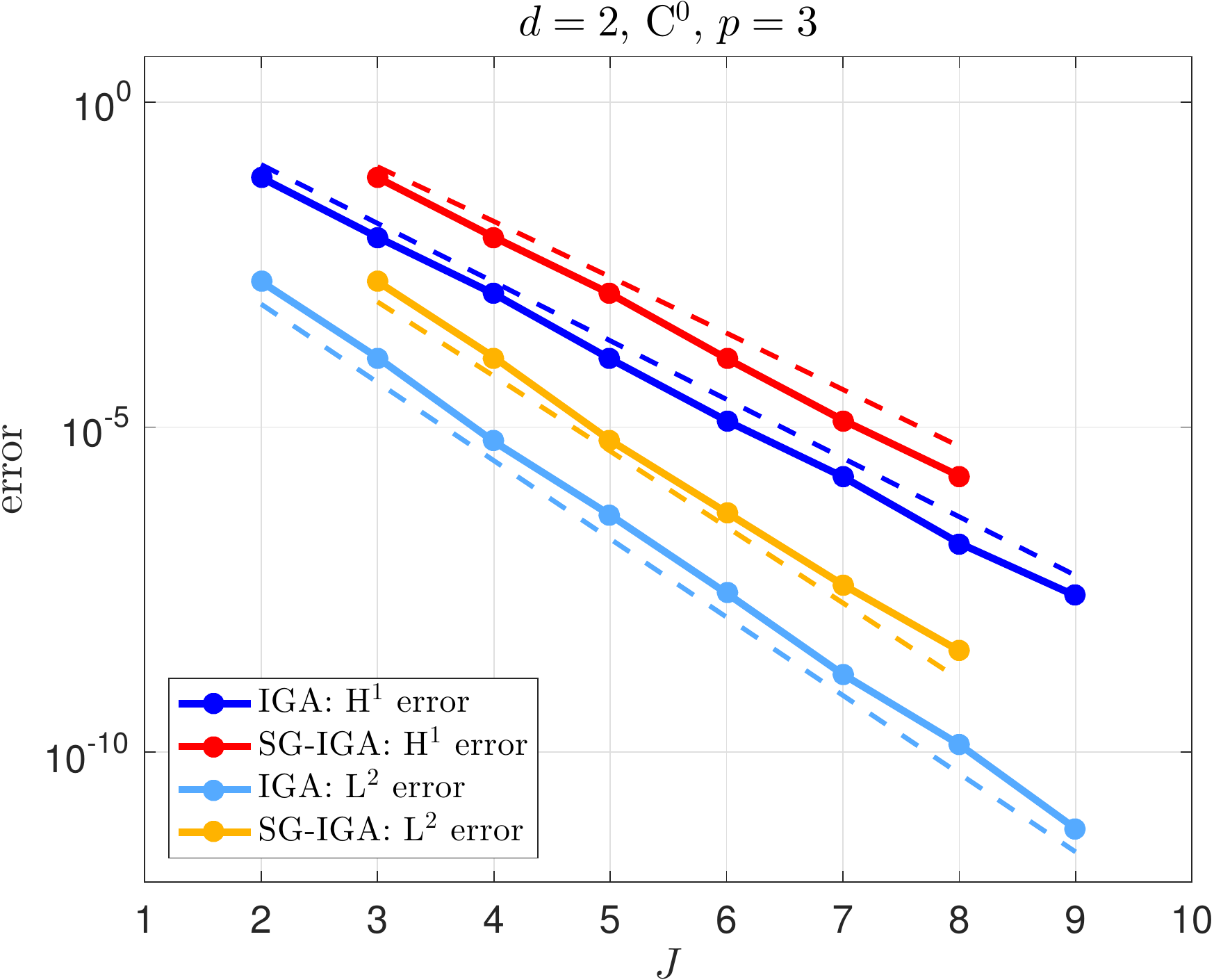}
  \includegraphics[width=\ConvSize\linewidth]{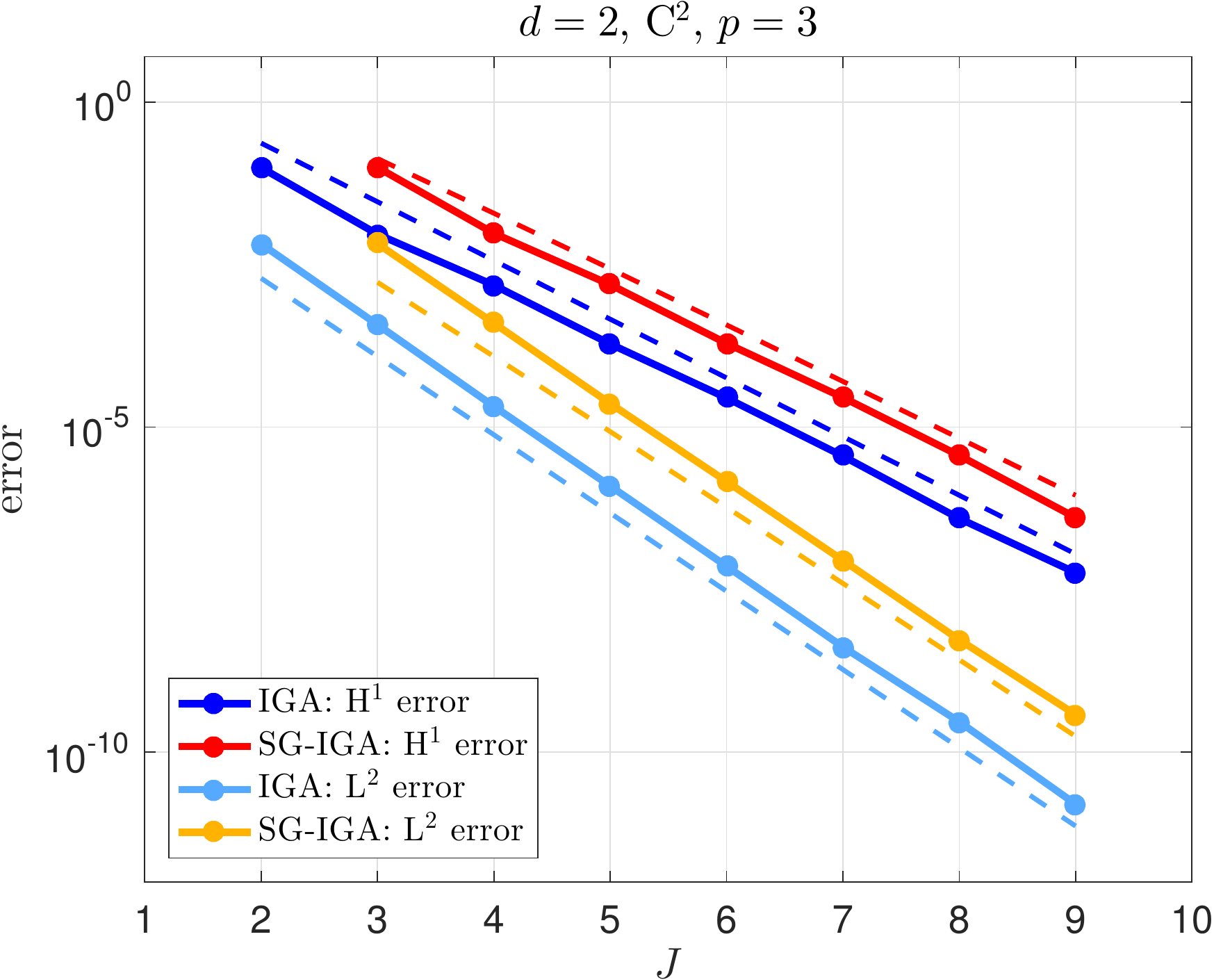} \\[\SpaceSize]
  \includegraphics[width=\ConvSize\linewidth]{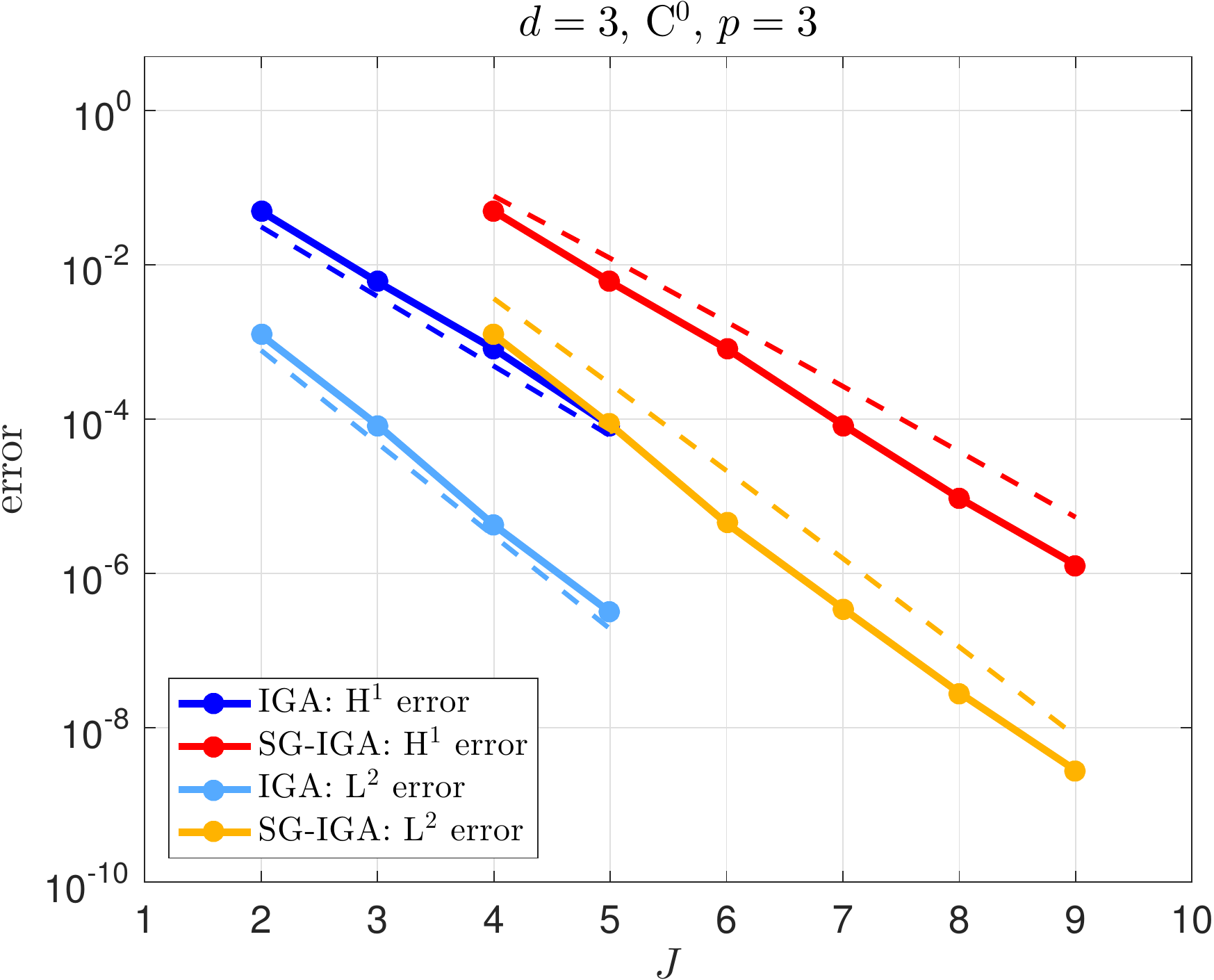}
  \includegraphics[width=\ConvSize\linewidth]{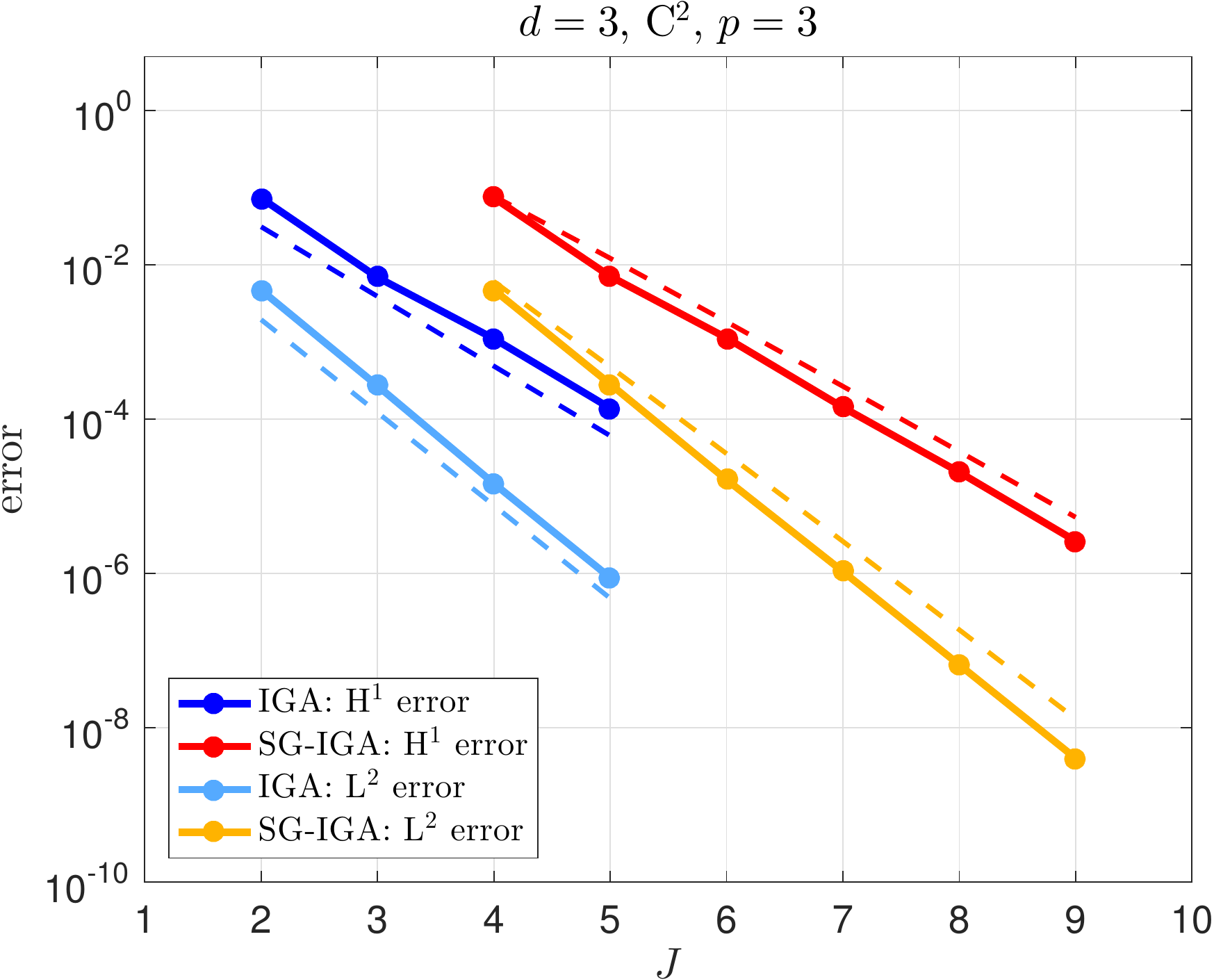}
  \caption{Error versus grid level for the problem on a quarter of annulus domain with regular solution,
    for $p=3$ and different values of $d$ (top row: $d=2$, bottom row: $d=3$)
    and of the regularity of the B-splines basis \La{(left column: $C^0$, right column: $C^2$)}. 
    The solid lines are the errors computed against the analytical solution, 
    and the dashed lines show the corresponding theoretical rates.}\label{fig:err-vs-w-reg-ring}
\end{figure}

The observed SG-IGA convergence rates agree well with the theoretical ones given in the parentheses.
In some cases, the resulting rates are even greater than expected: this is consistent with 
the observation in \cite{griebel.harbrecht:combi-conv}
that reports that sparse methods  can achieve the same convergence behaviour as standard tensor methods, (i.e., 
the estimate without the correction term $J^{(d-1)/2}$ in Equations \eqref{eq:H1-conv-wrt-mix} to \eqref{eq:L2-conv}
holds true) when the function has some additional regularity. 
Some of the convergence figures from which the rates in Table
\ref{tab:quarter-of-annulus-reg-rates} are deduced
are shown in Figure \ref{fig:err-vs-w-reg-ring}.  
The solid lines are the $H^1$ and $L^2$ errors
computed with respect to the analytical solution, and the dashed lines
are the corresponding theoretical slopes given in Table \ref{tab:sparse-grids-rates}.

\begin{figure}[tp]
  \centering
  \includegraphics[width=\ConvSize\linewidth]{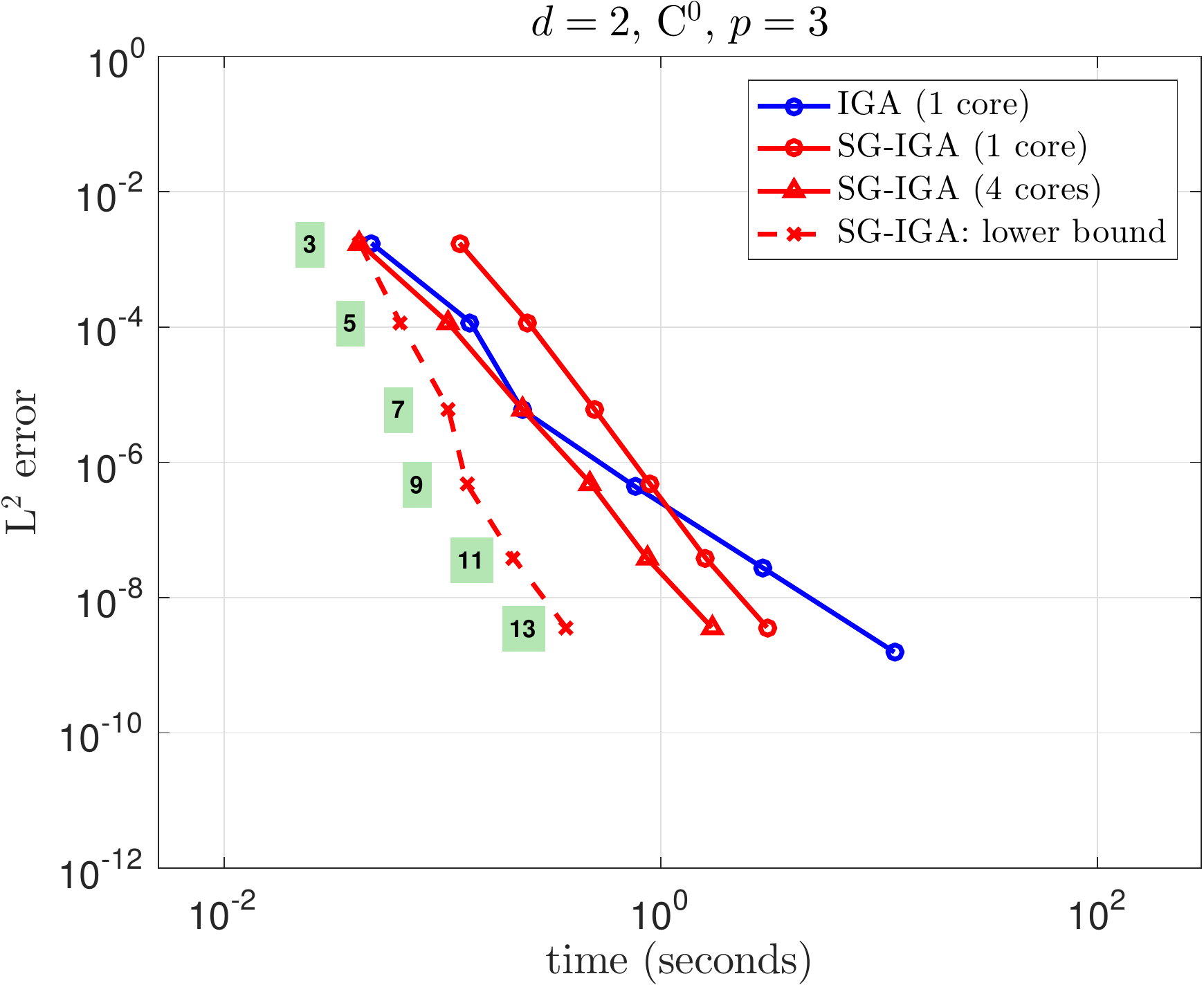}
  \includegraphics[width=\ConvSize\linewidth]{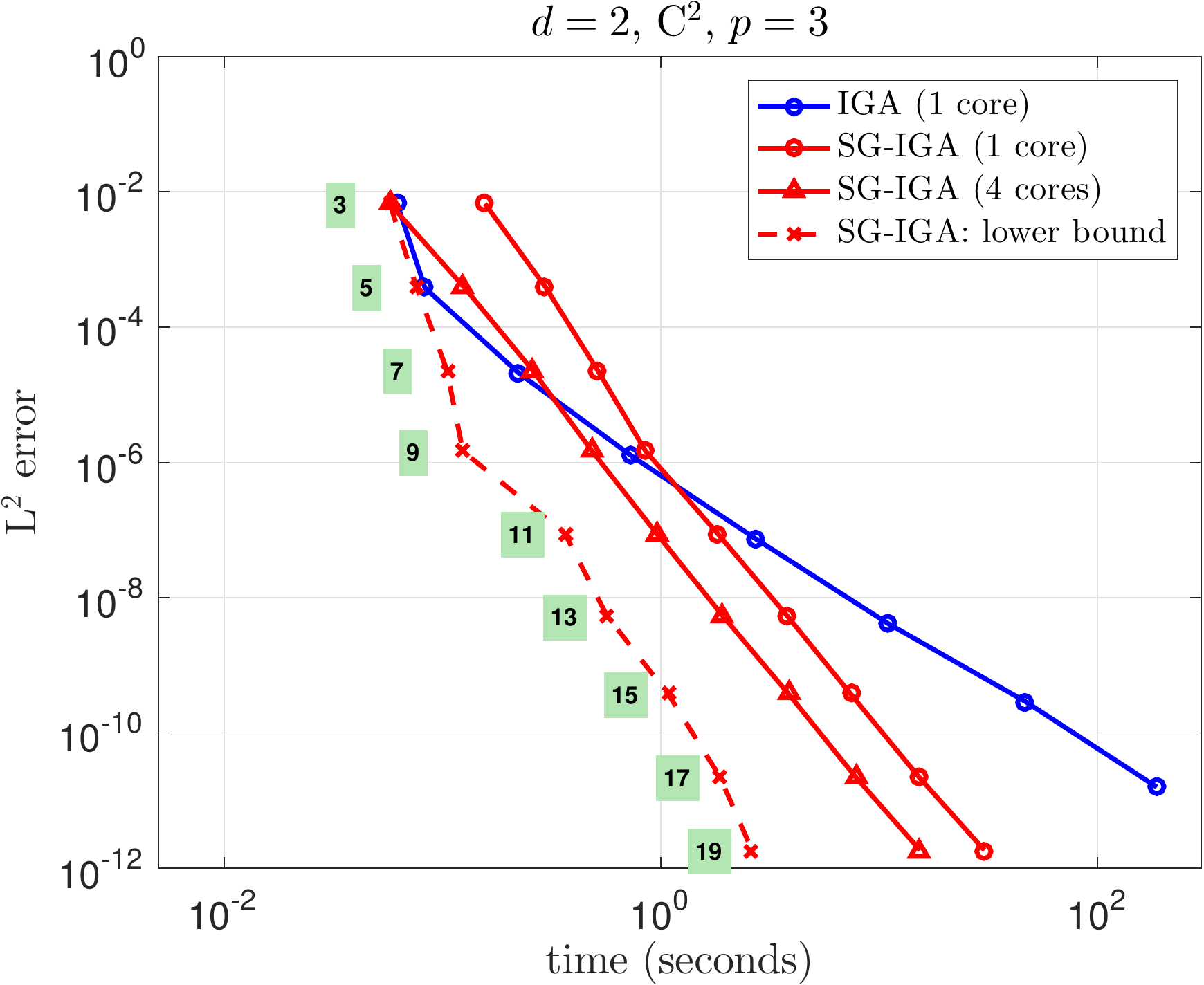} \\[\SpaceSize]
  \includegraphics[width=\ConvSize\linewidth]{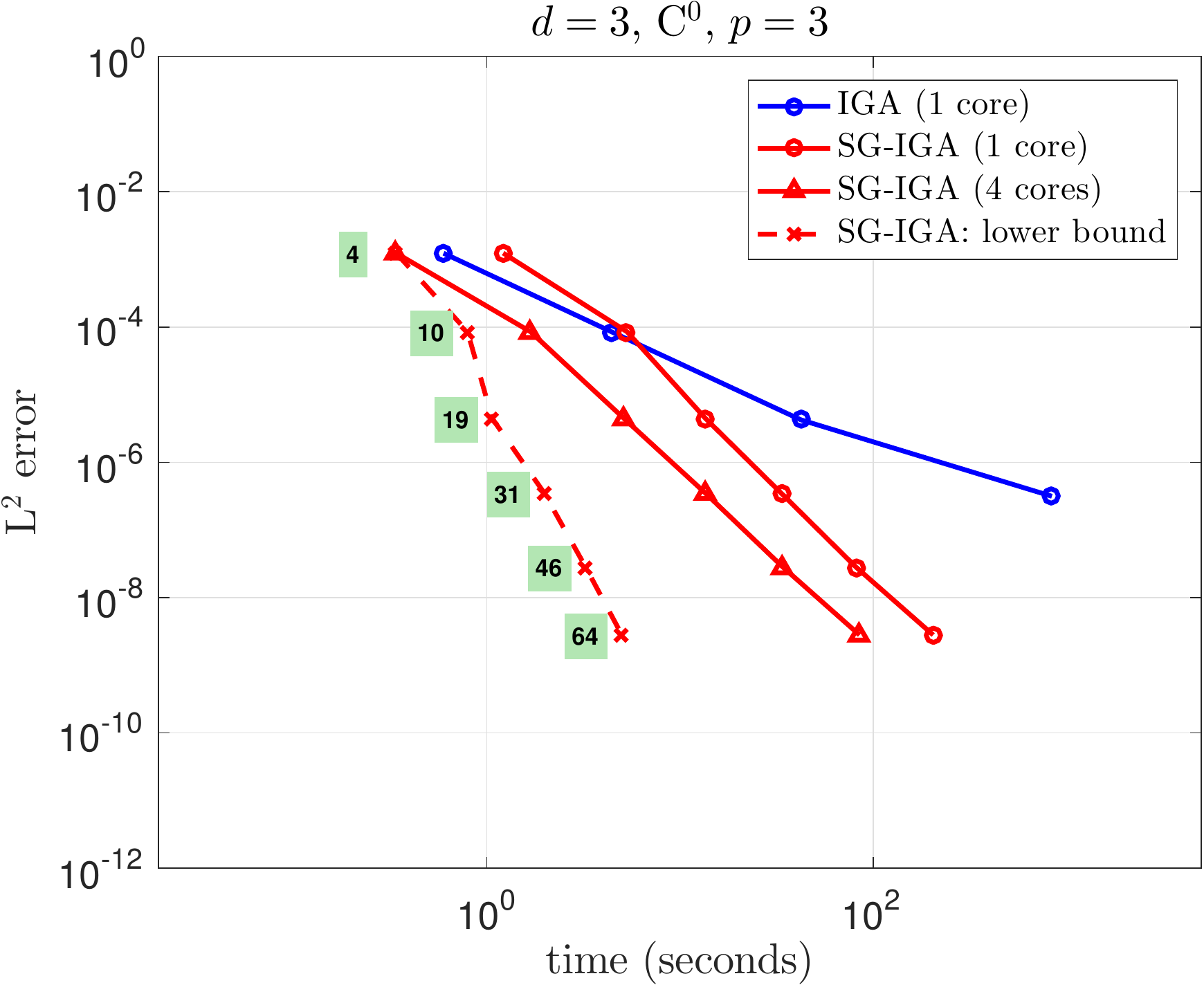}
  \includegraphics[width=\ConvSize\linewidth]{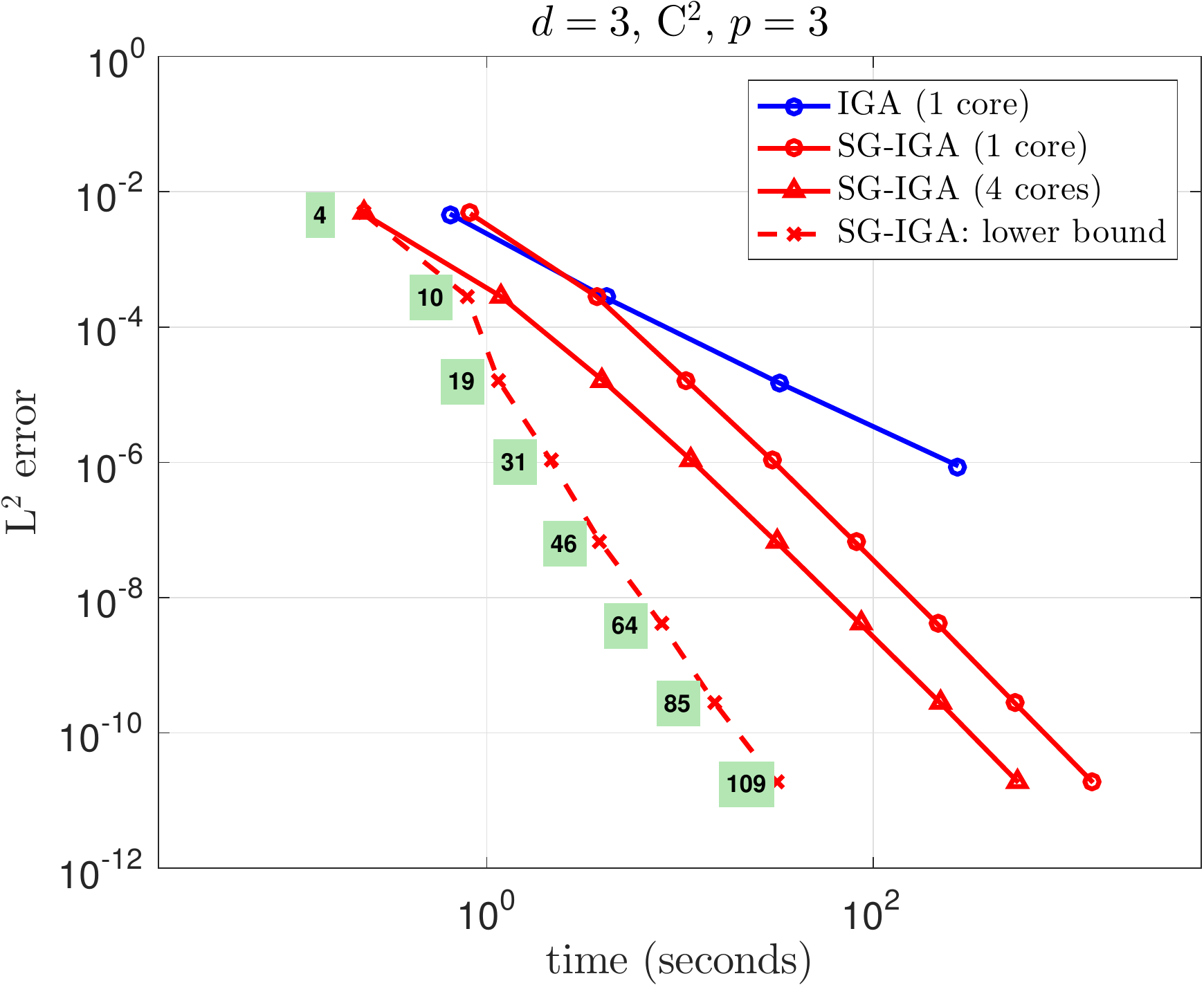}
  \caption{$L^2$ error vs. time for the problem on a quarter of annulus domain with regular solution.
    Here we fix $p=3$ and change $d$ (top row: $d=2$, bottom row: $d=3$) and the regularity of the B-splines basis
    \La{(left column: $C^0$, right column: $C^2$)}. The dashed line is the lower bound that 
    can be achieved if the number of available cores is at least equal to the number of components 
    of the combination technique for each level, given by the numbers in green boxes.}\label{fig:err-vs-time-reg-ring} 

  \bigskip

  \includegraphics[width=\ConvSize\linewidth]{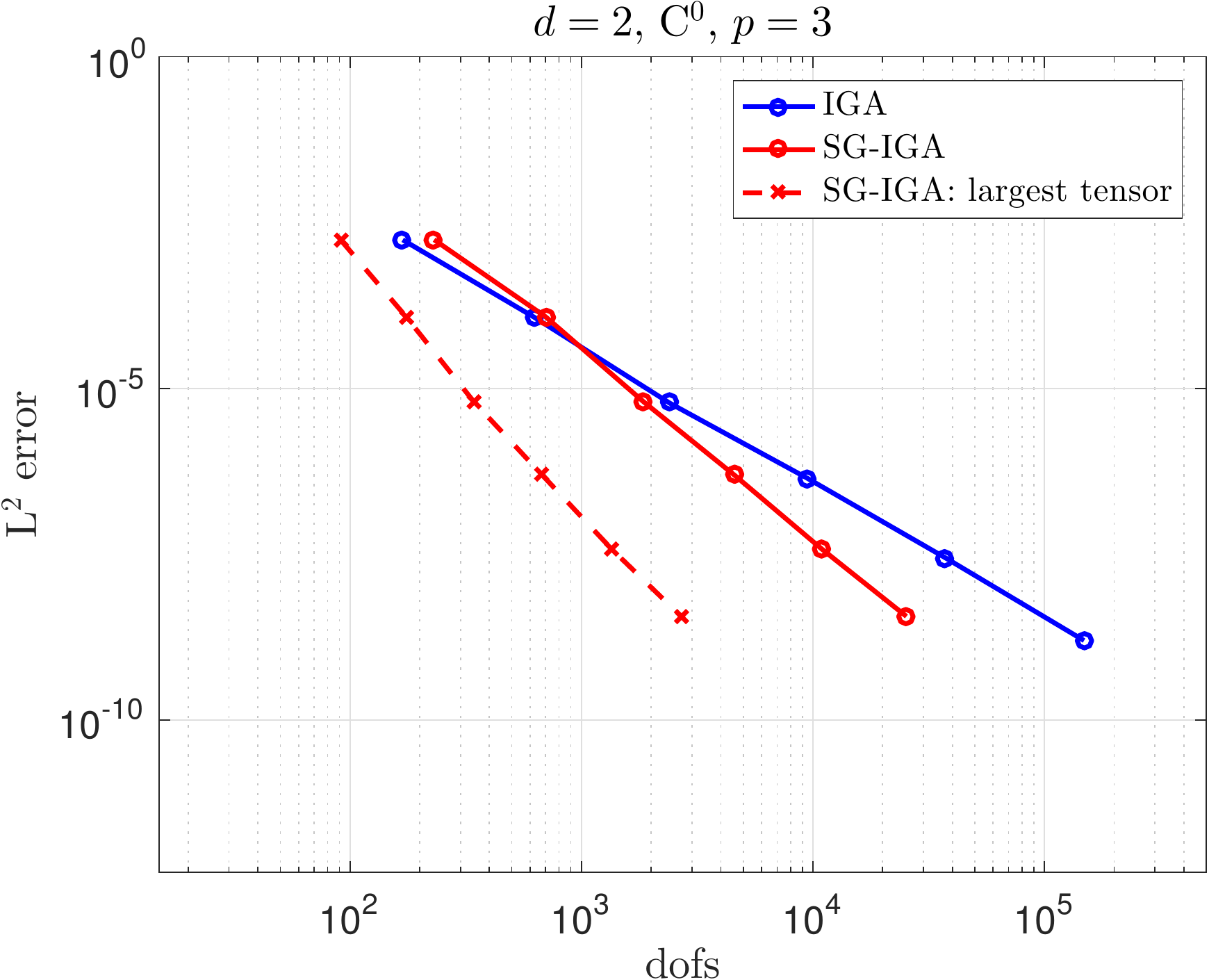}
  \includegraphics[width=\ConvSize\linewidth]{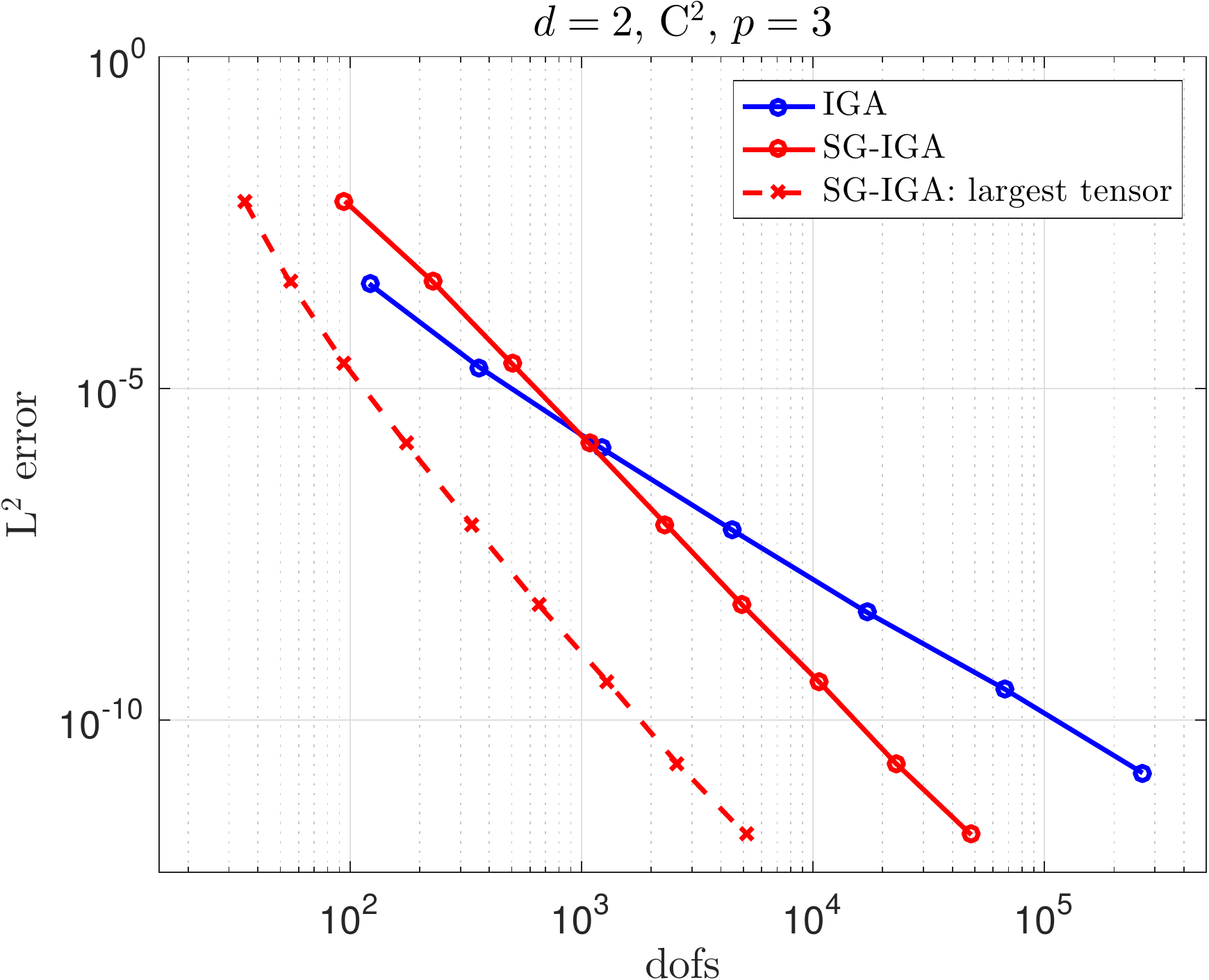} \\[\SpaceSize]
  \includegraphics[width=\ConvSize\linewidth]{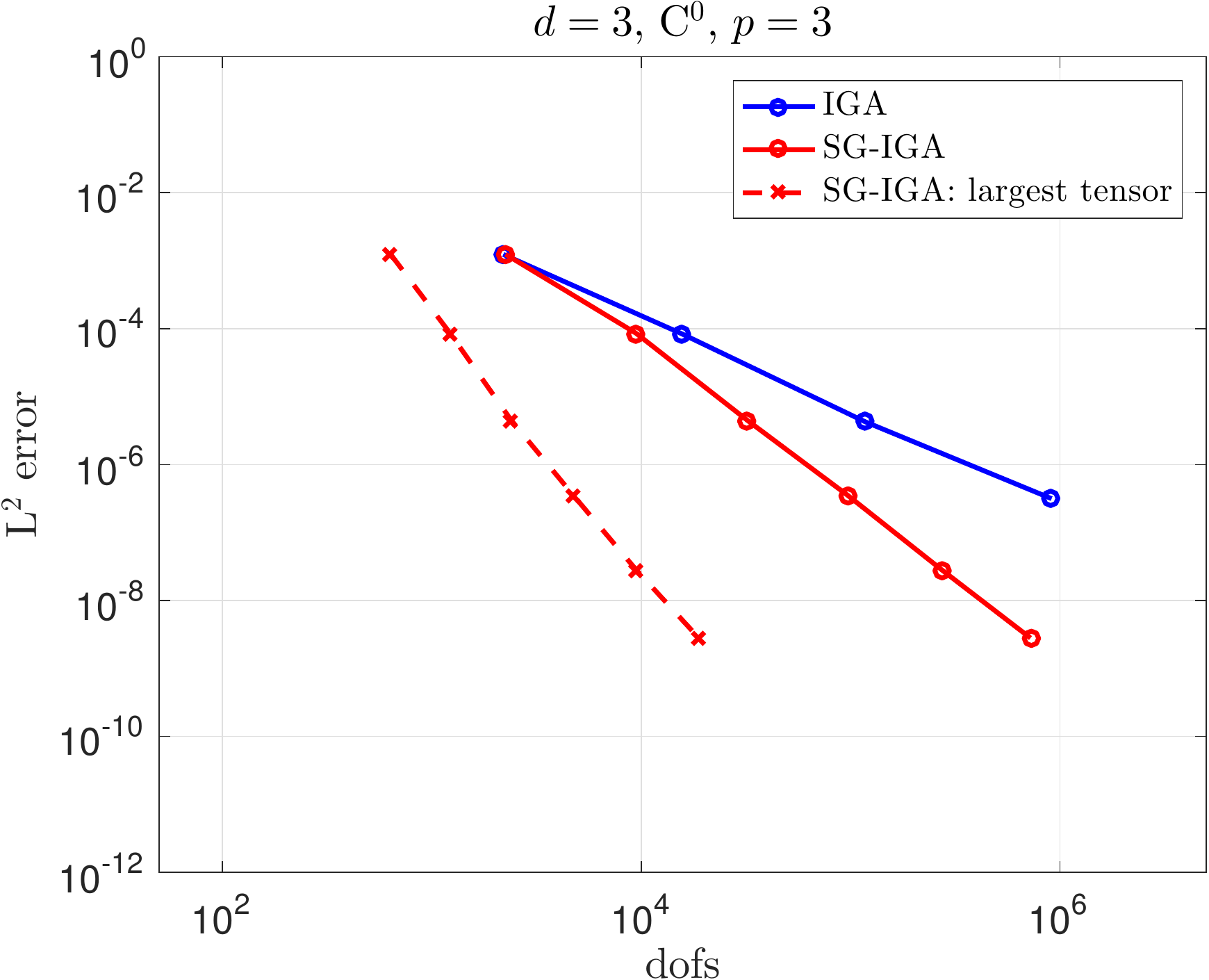}
  \includegraphics[width=\ConvSize\linewidth]{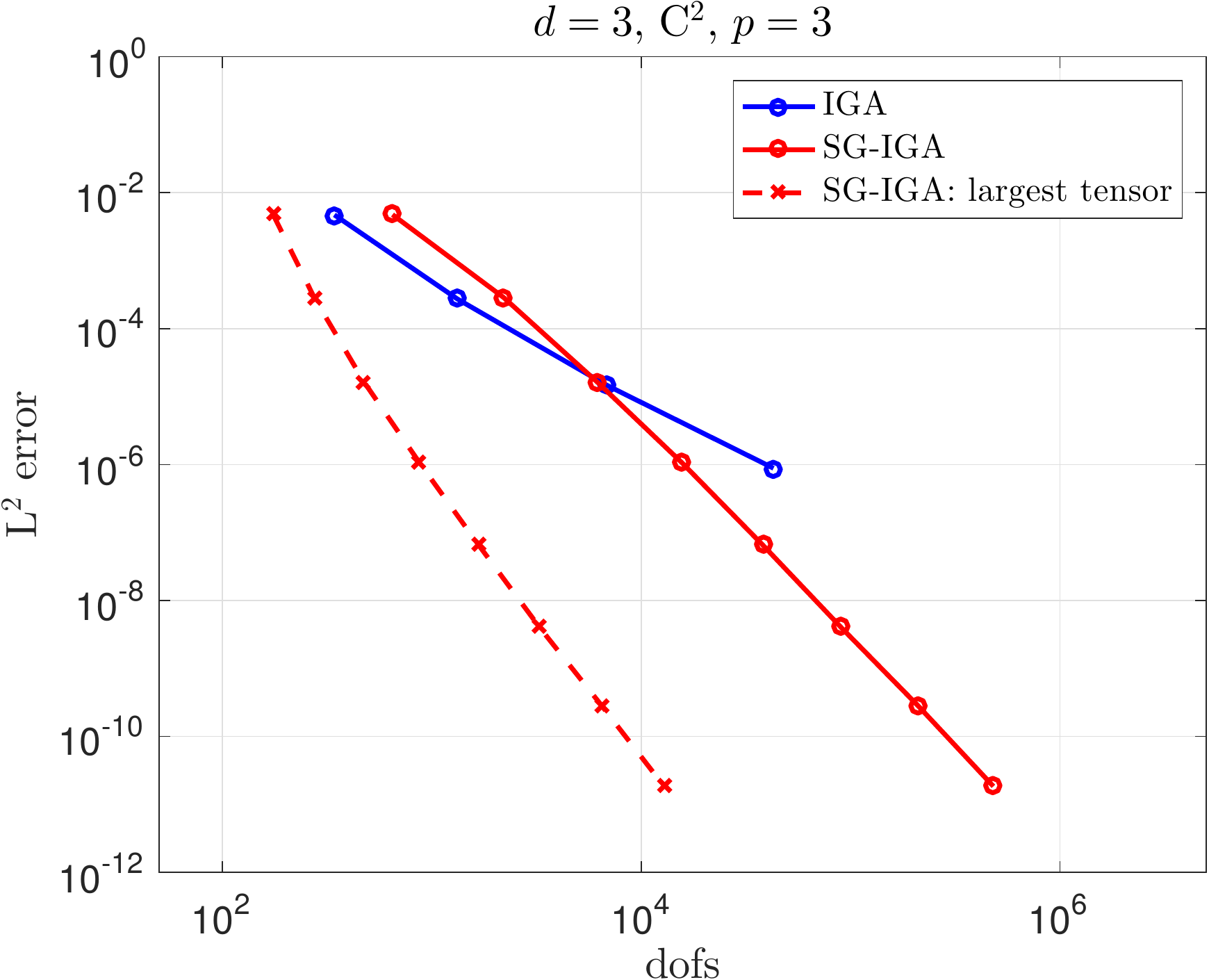}

  \caption{$L^2$ error vs. dofs for the problem on a quarter of annulus domain with regular solution.
    Here we fix $p=3$ and change $d$ (top row: $d=2$, bottom row: $d=3$) and the regularity of the B-splines basis
    \La{(left column: $C^0$, right column: $C^2$)}. 
    In each figure, the dashed line represents the error of SG-IGA vs. the number
    of the degrees-of-freedom for the largest tensor grid in the combination technique.
  }
  \label{fig:err-vs-dofs-reg-ring}

\end{figure}

The error versus computational time is shown in Figure \ref{fig:err-vs-time-reg-ring}:
here we only show the case $p=3$ with minimal and maximal regularity for brevity.
We can see that SG-IGA is an effective alternative to the standard IGA for this problem: 
for $d=2$, this is especially true in the multi-core situation when 4 or more cores are used,
while in the case $d=3$ even using one core gives clear advantages. 
The dashed line indicates the lower bound on the computational cost that can be achieved if the 
number of available cores is at least equal to the number of components of the combination technique 
to be computed (such number is shown to the left in green boxes).
The gain in computational effort is evident also if we consider 
the number of degrees-of-freedom as a computational cost indicator, 
cf. Figure \ref{fig:err-vs-dofs-reg-ring}. 
In contrast to the dashed line in the error vs. time, the dashed line here is the number of the 
degrees-of-freedom for the largest tensor grid in the combination technique.

\begin{figure}[tp]
  \centering
  \includegraphics[width=\OptSetSize\linewidth]{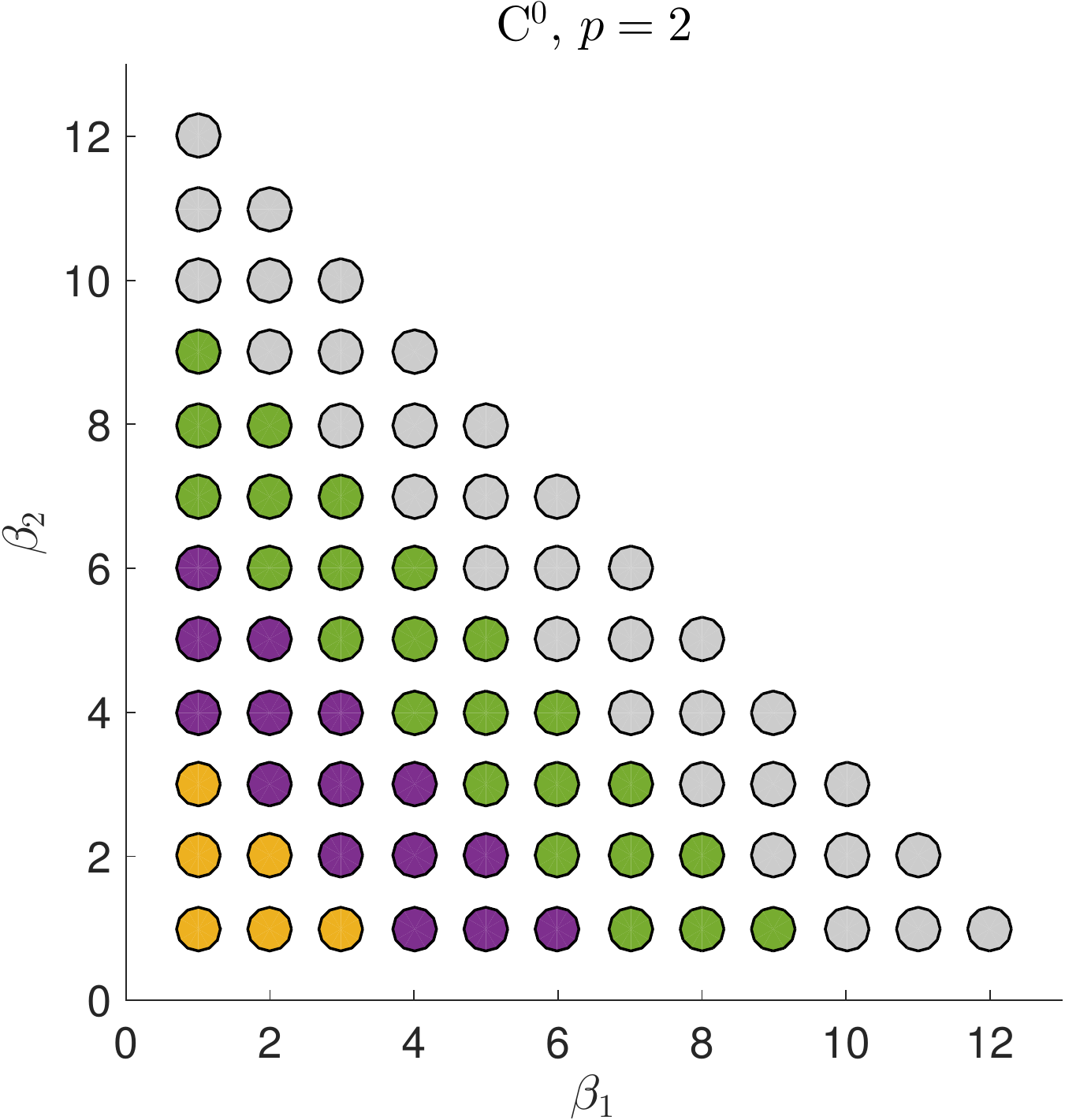}
  \includegraphics[width=\OptSetSize\linewidth]{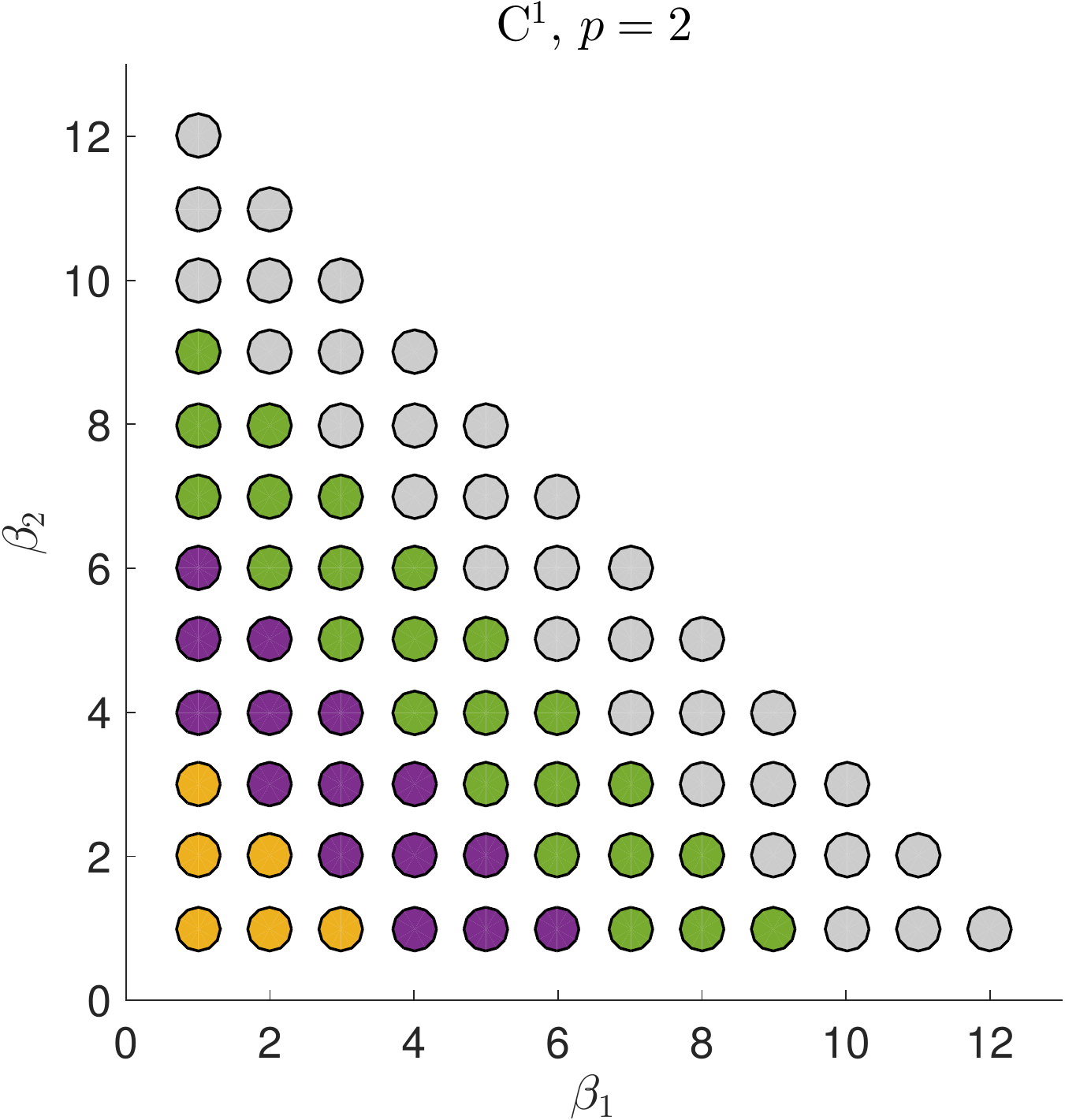} \\[\SpaceSize]
  \includegraphics[width=\OptSetSize\linewidth]{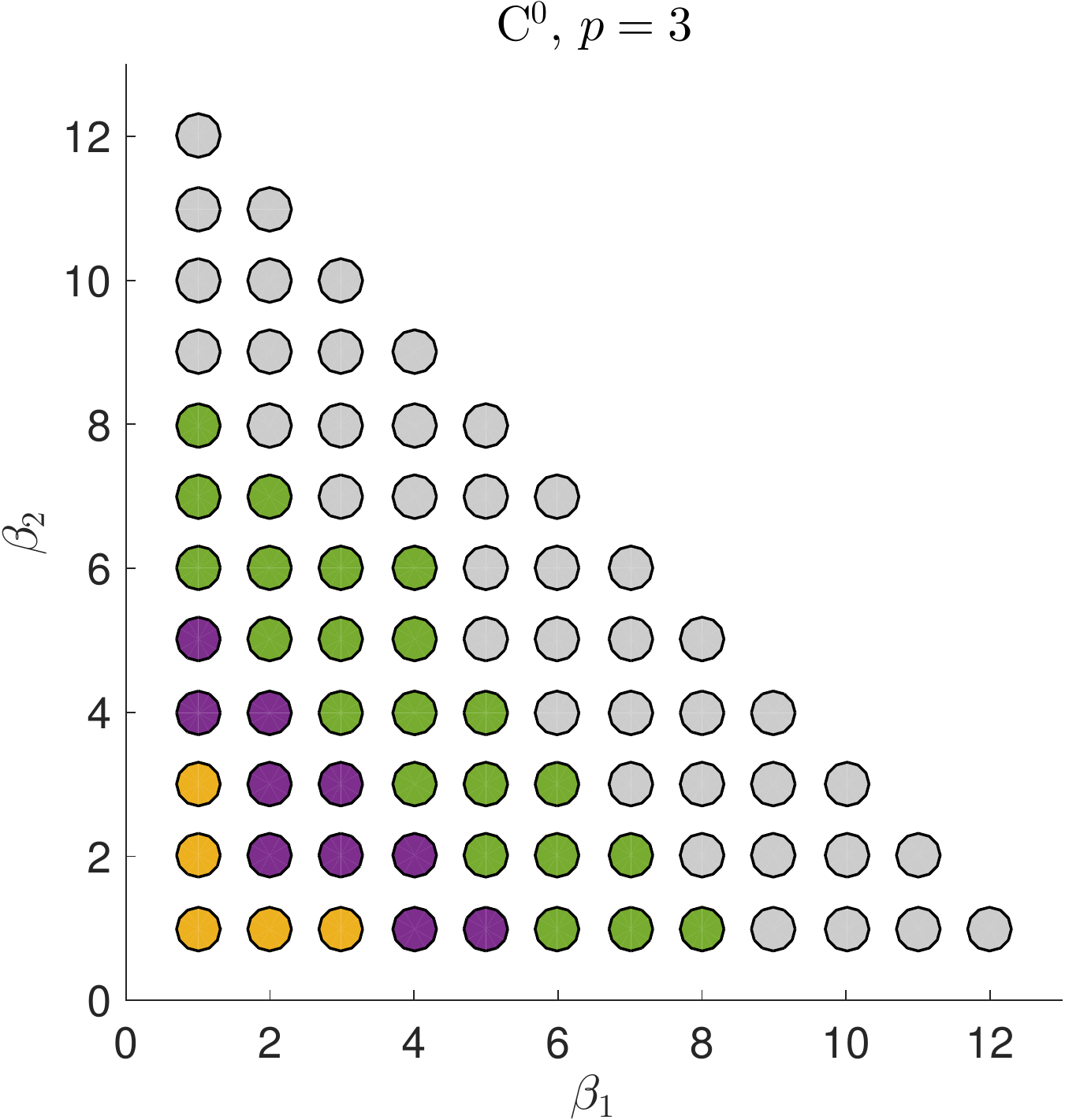}
  \includegraphics[width=\OptSetSize\linewidth]{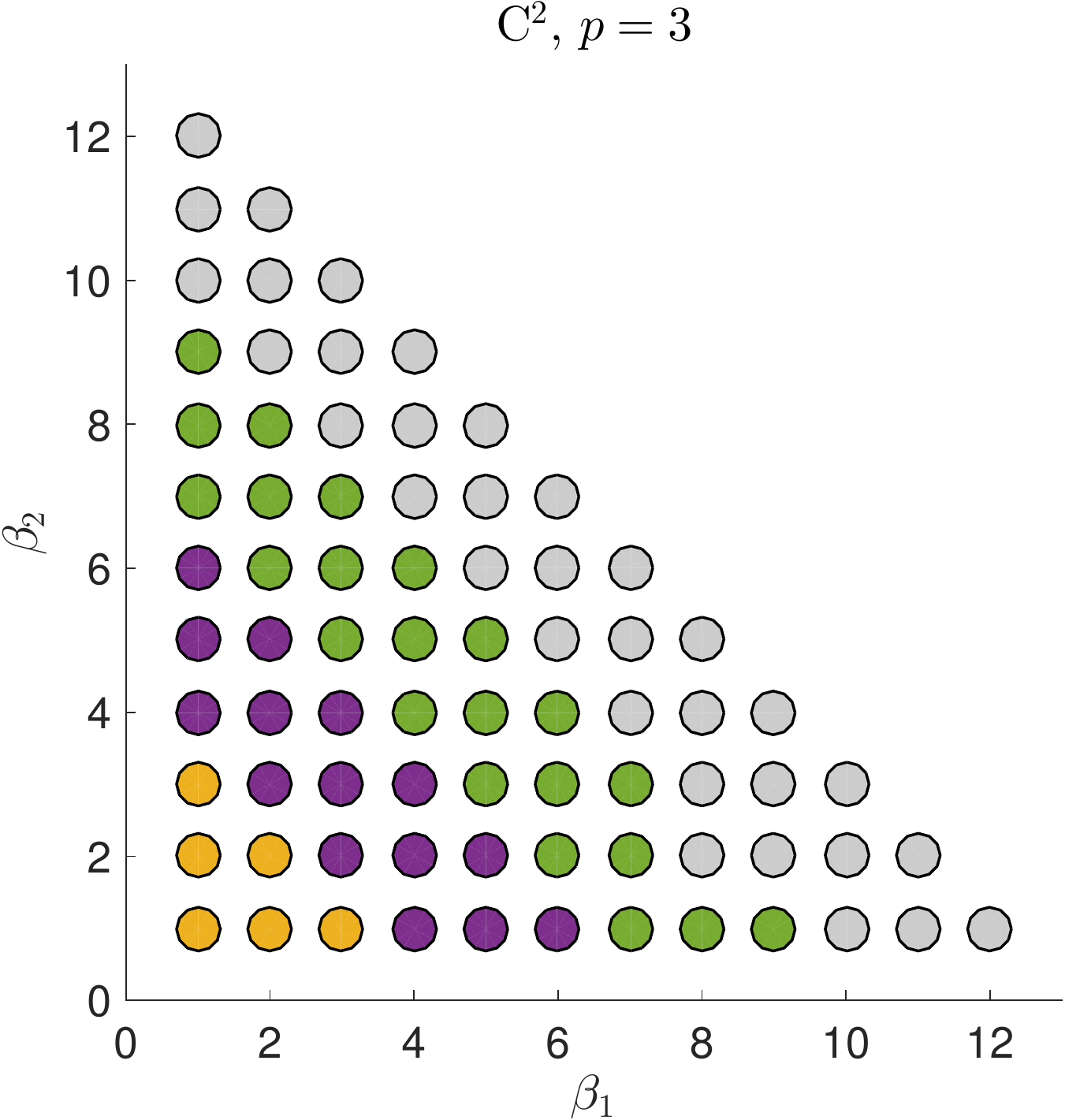} \\[\SpaceSize]
  \includegraphics[width=\OptSetSize\linewidth]{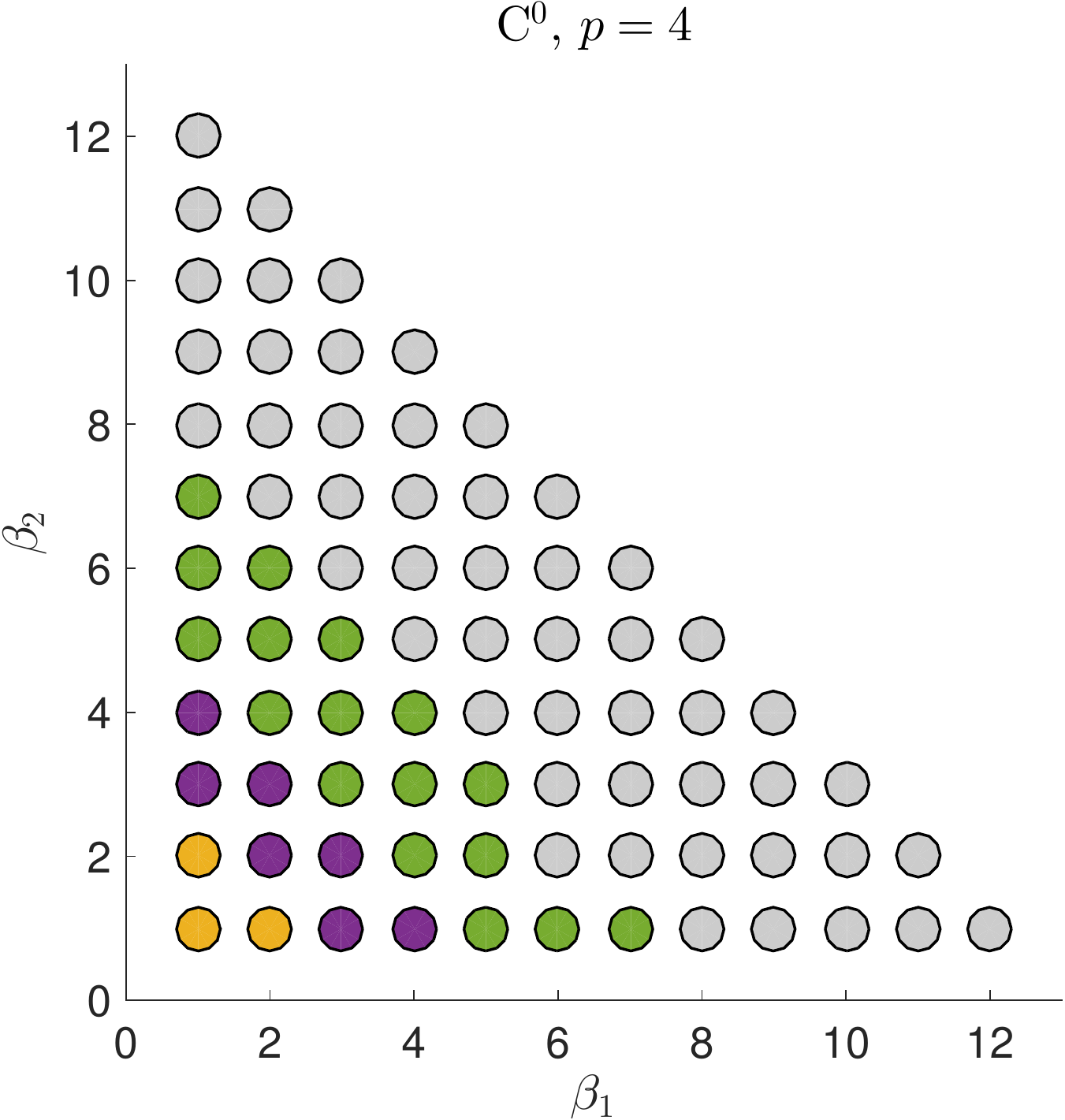}
  \includegraphics[width=\OptSetSize\linewidth]{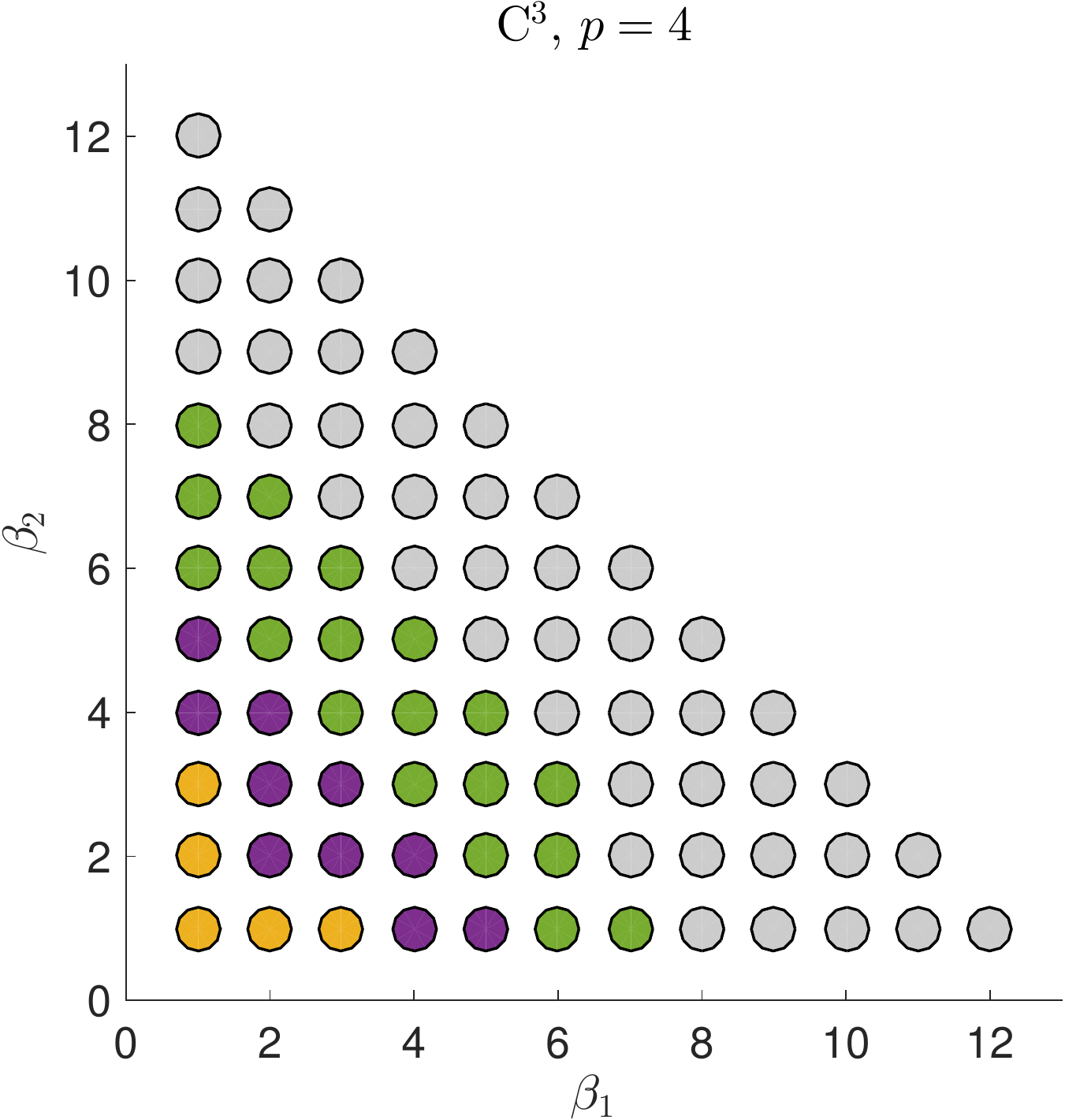} 
  \caption{optimal sets for the quarter-of-annulus problem with regular solution. 
    From top to bottom row: $p=2,3,4$. Left column: $C^0$ B-spline basis; Right column: maximal continuity
    B-splines basis.
}\label{fig:opt-sets-ring-reg}
\end{figure}

Finally, 
we give for this test problem a computational validation of the truncation procedure
that led to defining the sparse-grid approximation \eqref{eq:TD-approx}, which we now revisit. 
To obtain \eqref{eq:TD-approx}, we first introduced the decomposition of a full tensor IGA approximation $u_{[J,\ldots,J]}$ as
\[
u_{[J,\ldots,J]} = \sum_{\bbeta \in \mathcal{I}_{J,\text{IGA}}} \Delta_{\bbeta}(u), \quad 
\mathcal{I}_{J,\text{IGA}} = \{\bbeta \in \Rset^d : \max_{\kappa=1,\ldots,d} \beta_\kappa \leq J\},
\]
cf. Equations \eqref{eq:Delta-op} to \eqref{eq:TP-telescopic}, and then ``pruned'' the set $\mathcal{I}_{J,\text{IGA}}$ to obtain 
\[
u_{J} = \sum_{\bbeta \in \mathcal{I}_{J,\text{SG-IGA}}} \Delta_{\bbeta}(u), \quad 
\mathcal{I}_{J,\text{SG-IGA}} = \Big\{\bbeta \in \Rset^d : \sum_{\kappa=1}^d \beta_\kappa \leq J\Big\},
\]
under the assumption that the magnitude of $\|\Delta_{\bbeta}(u)\|_{L^2(\Omega)}$  
would be decreasing with respect to $\sum_{\kappa=1}^d \beta_\kappa$.
One may therefore set up an optimization problem, and look for the best set $\mathcal{I}$, i.e., the set that delivers the best 
approximation of $u$ for a given computational budget, 
\begin{align*}
& \min_{\mathcal{I} \subset \Nset^d} \|u-u_{\mathcal{I}}\|_{L^2(\Omega)}, \\
& \text{such that}  \quad u_{\mathcal{I}} = \sum_{\bbeta \in \mathcal{I}} \Delta_{\bbeta}(u) \\ 
& \quad \quad \quad \quad \sum_{\bbeta \in \mathcal{I}} \text{degrees-of-freedom}(\Delta_{\bbeta}(u)) \leq K.
\end{align*}
This is a classic ``knapsack'' optimization problem \cite{martello:knapsack} and can be approximately solved by the so-called
Dantzig method: 
\begin{enumerate}
\item define a ``revenue'' and a ``cost'' for each item that can be picked;
\item define the ``profit'' for each item, by taking the ratio of ``revenue'' over ``cost'';
\item sort items according to profit in decreasing order;
\item pick up items according to the ordering just introduced, until the sum of their costs exceeds the budget constraint.
\end{enumerate}
In our case, the cost of a $\Delta_{\bbeta}(u)$ would be the number of degrees-of-freedom 
and the revenue would be its $L^2$ norm, which quantifies how much the 
sparse-grid approximation $u_{\mathcal{I}}$ would improve if $\Delta_{\bbeta}(u)$ was added to the sparse-grid approximation. 
If one had a-priori estimates on the decay/growth of cost and revenue, one could then devise the ``optimal''
sparse-grid approximation for a fixed cost, see \cite{b.griebel:acta,griebel.knapek:optimized,nobile.eal:optimal-sparse-grids} for more details.
Observe that we can compute cost and revenue for every $\Delta_{\bbeta}(u)$ with $\bbeta$ in a sufficiently large set,
cf. Equation \eqref{eq:Delta-op} and \eqref{eq:Delta-combi-tec}. 
The results of this computation for the two-dimensional quarter-of-annulus problem
are reported in Figure \ref{fig:opt-sets-ring-reg}, where each dot represents 
a multi-index $\bbeta$ in the plane $(\beta_1, \beta_2)$, and dots are colored 
according to the value of their associated profit. 
It is clearly visible that indices with $\beta_1+\beta_2 = constant$ have profit of equal magnitude, which means
that they should be picked together in the approximation, thus obtaining exactly the sparse-grid expression
\eqref{eq:TD-approx}. 
This will no longer be the case for the next numerical test.

\subsection{Quarter of annulus domains, low-regular solution}


\begin{table}[tp]
  \centering
\begin{tabular}{l|rr|rr}
					& \multicolumn{2}{c|}{$d=2$}	& \multicolumn{2}{c}{$d=3$}	\\
					& $H^1$ error 	& $L^2$ error	& $H^1$ error 	& $L^2$ error\\
$C^0, p=2$, SG-IGA	&	1.35(1)		&	1.72(1.5)	&	1.30(2/3)	&	1.62(1)	\\
$C^0, p=2$, IGA		&	2.11(2)		&	2.71(3)		&	1.84(2)		&	2.70(3)	\\
$C^1, p=2$, SG-IGA	&	1.20(1)		&	1.72(1.5)	&	1.22(2/3)	&	1.63(1)	\\
$C^1, p=2$, IGA  	&	2.02(2)		&	2.82(3)		&	1.81(2)		&	2.88(3)	\\
$C^0, p=3$, SG-IGA	&	1.26(1)		&	1.74(1.5)	&	1.30(2/3)	&	1.71(1)	\\
$C^0, p=3$, IGA		&	2.21(2)		&	3.00(3)		&	2.21(2)		& 	2.89(3)	\\
$C^1, p=3$, SG-IGA	&	1.25(1)		&	1.78(1.5)	&	1.28(2/3)	&	1.75(1)	\\
$C^1, p=3$, IGA		&	2.14(2)		&	3.05(3)		&	2.15(2)		&	2.90(3) \\
$C^2, p=3$, SG-IGA	&	1.19(1)		&	1.58(1.5)	&	1.20(2/3)	&	1.55(1)	\\
$C^2, p=3$, IGA		&	2.21(2)		&	2.92(3)		&	2.12(2)		&	2.81(3) \\
\end{tabular}  
  \caption{Convergence of SG-IGA and IGA methods with respect to the sparse-grid level, $J$ for the low-regular problem. 
    The number in parenthesis indicates the expected rate from theory, cf. Table \ref{tab:sparse-grids-rates}.} 
  \label{tab:quarter-of-annulus-non-reg-rates}
\end{table}

In this second set of experiments we consider again the quarter-of-annulus domain, 
but we choose as forcing term in Equation \eqref{eq:poisson}
the function $f(\xx)=1$. This implies that the solution, $u$, will have limited regularity
due to corner and edge singularities, i.e., $u \in H^{3-\epsilon}(\Omega)$ 
for $\epsilon>0$ but $u \not \in H^{3}(\Omega)$, see \cite{dauge2006elliptic}; 
therefore, we expect SG-IGA to converge at a lower rate than for the previous problem.

Table \ref{tab:quarter-of-annulus-non-reg-rates} reports the measured convergence rate with 
respect to the sparse-grid level as well as the lower bound for the error estimates
\eqref{eq:H1-conv}-\eqref{eq:H1-conv-tensor} and \eqref{eq:L2-conv}-\eqref{eq:L2-conv-tensor}
In the case $d=2$, the measured rates are in close agreement with the lower bounds: 
this is also consistent with what was observed in the previous problem. In the case $d=3$ instead
the convergence of SG-IGA is significantly better than the lower bounds
\eqref{eq:H1-conv} and \eqref{eq:L2-conv}.
This may be related to a mixed Sobolev regularity of the solution 
higher than $H^{1,1,1}_{mix}(\widehat{\Omega})$ or $\mcH^{2/3,2/3,2/3}(\widehat{\Omega})$.
In any case, there is a significant difference between the convergence
rate of SG-IGA and IGA for this problem.



\begin{figure}[tp]
  \centering
  \includegraphics[width=\OptSetSize\linewidth]{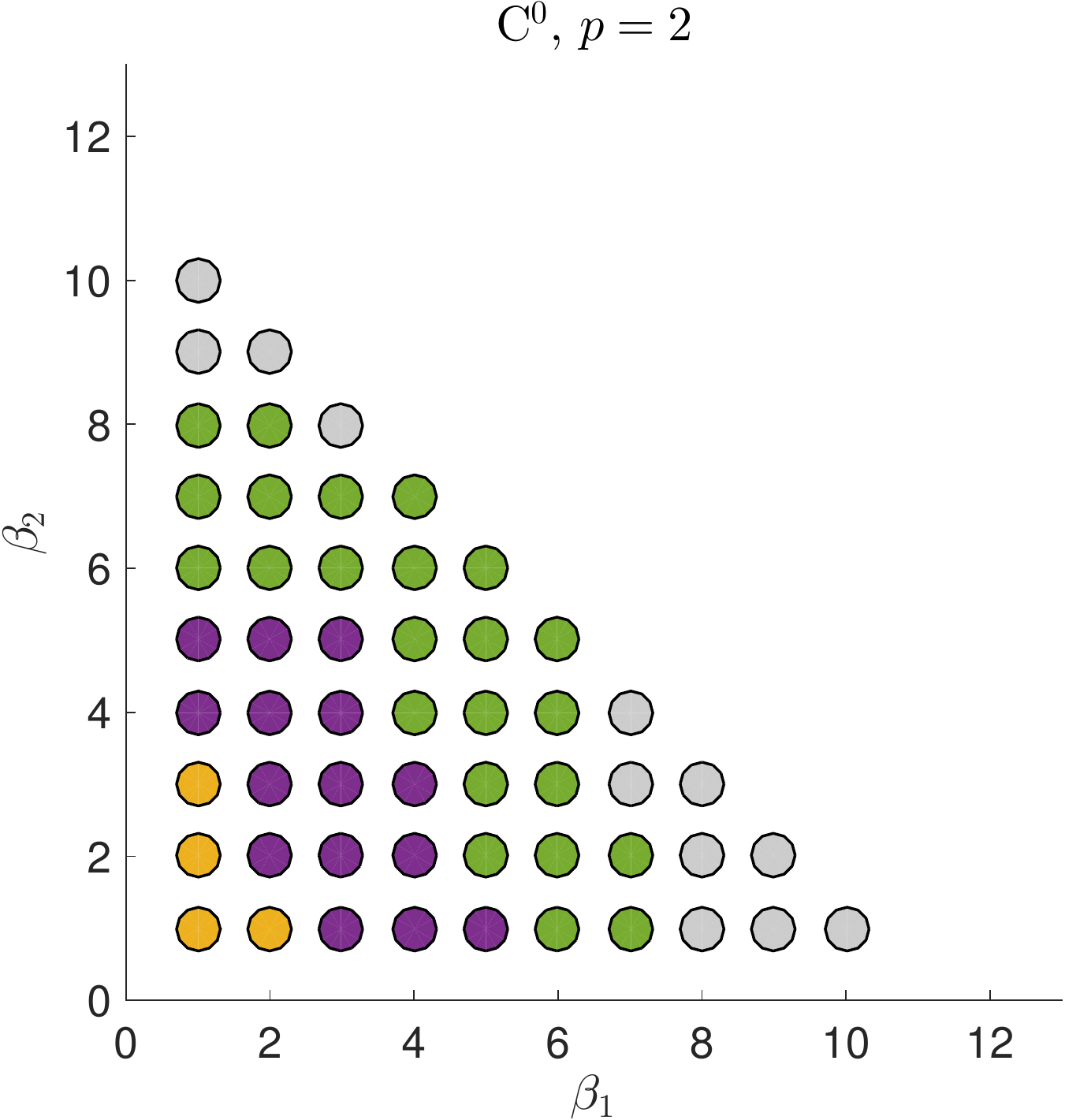}
  \includegraphics[width=\OptSetSize\linewidth]{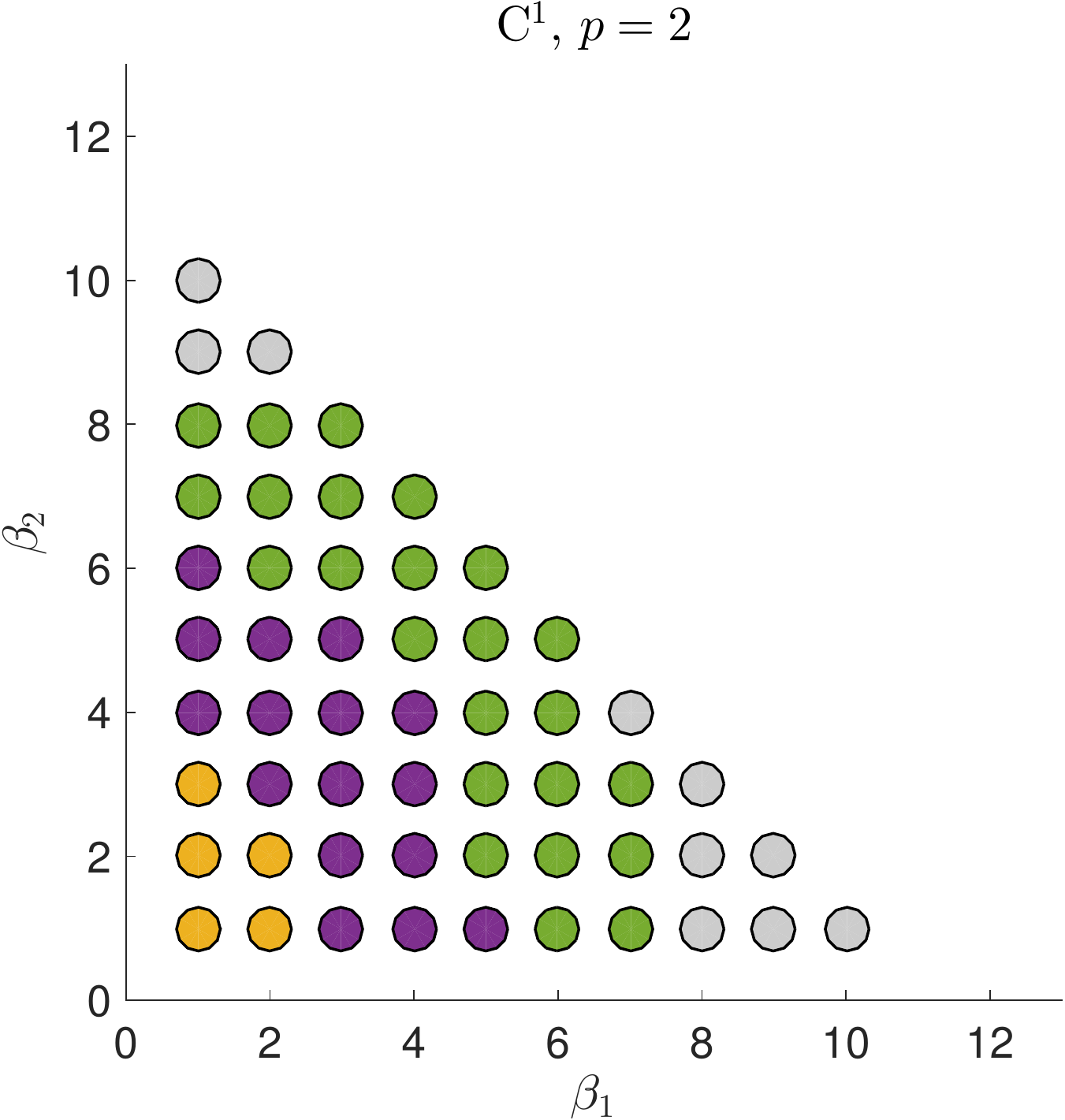} \\[\SpaceSize]
  \includegraphics[width=\OptSetSize\linewidth]{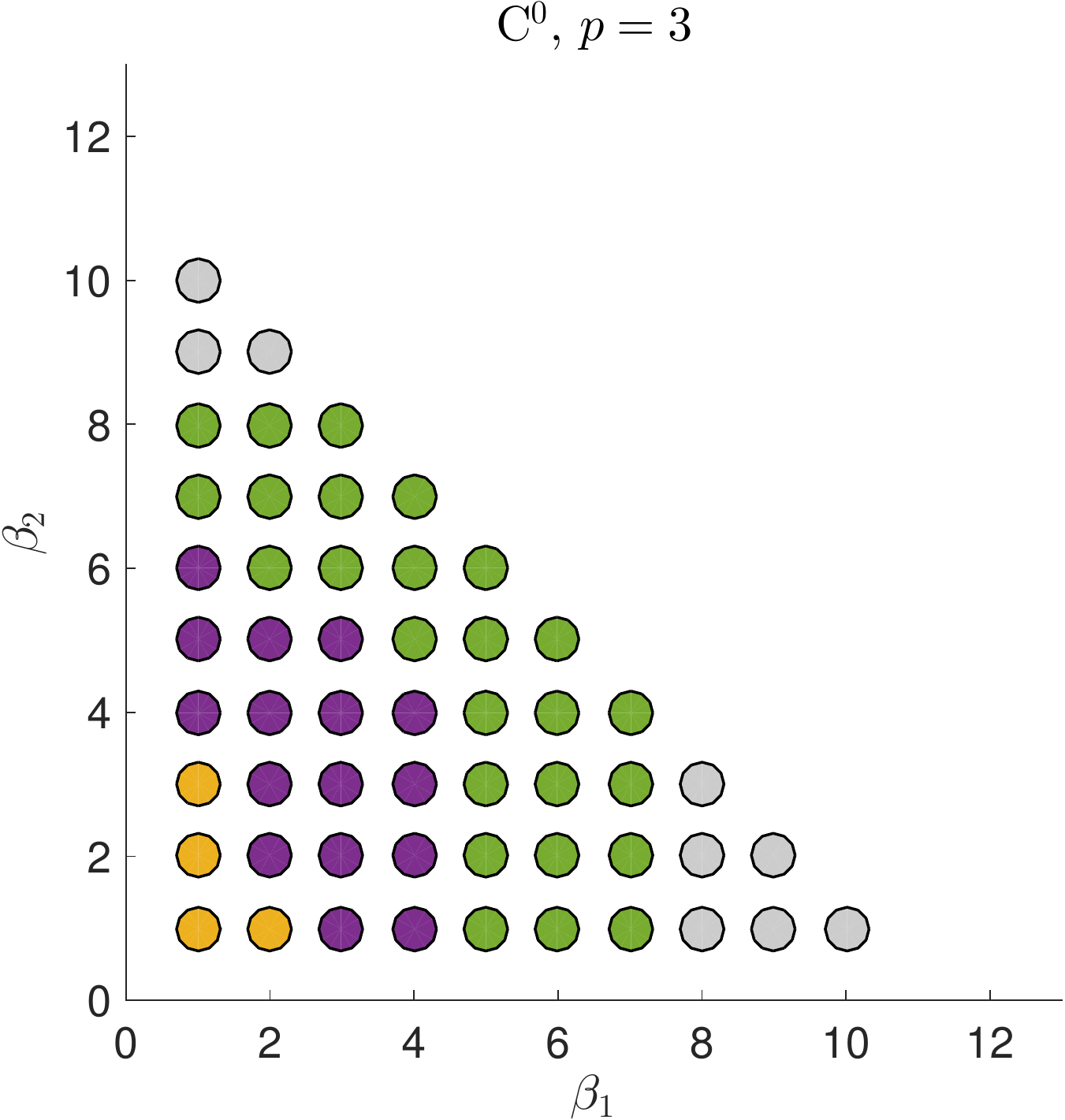}
  \includegraphics[width=\OptSetSize\linewidth]{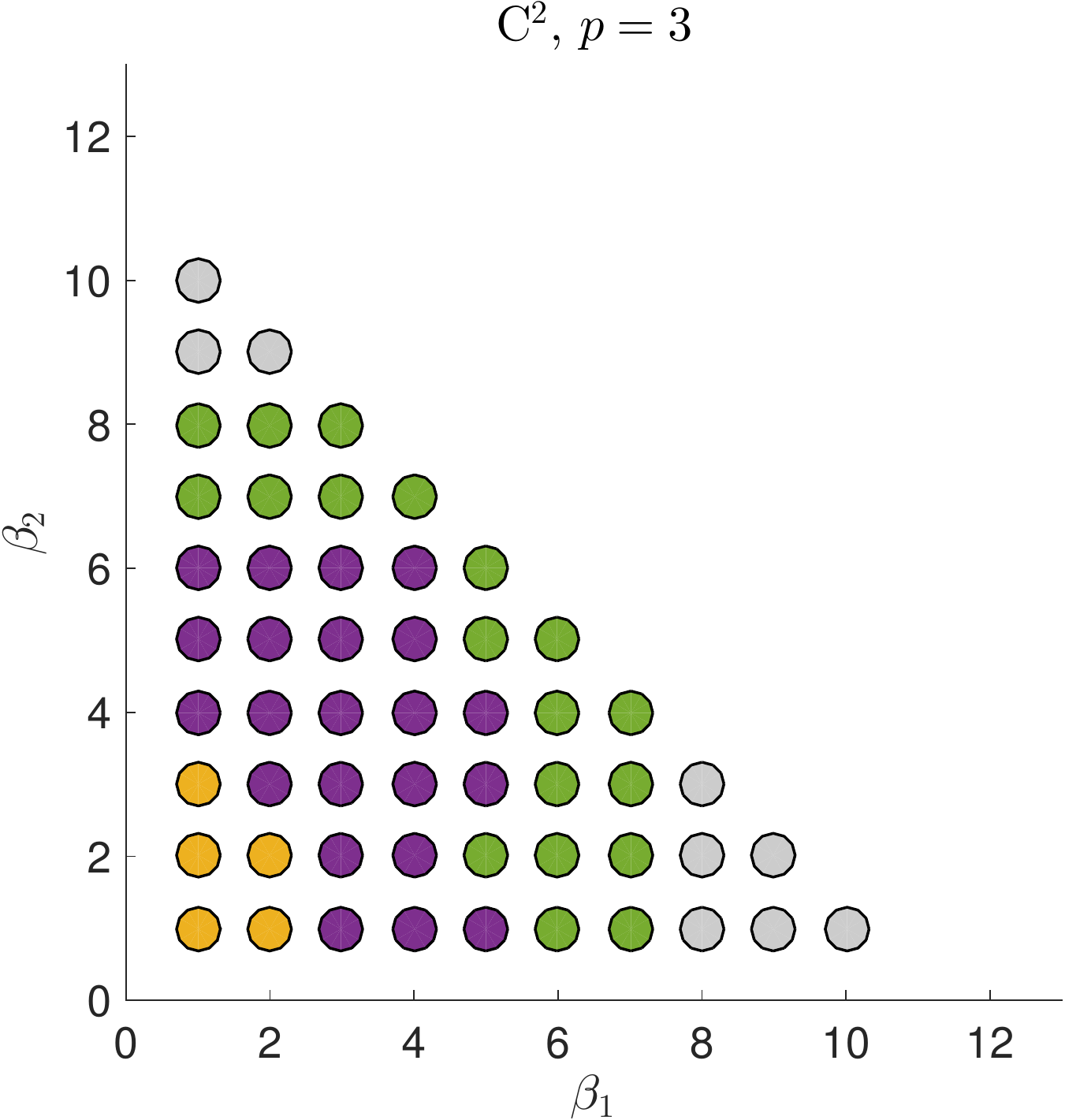} \\[\SpaceSize]
  \includegraphics[width=\OptSetSize\linewidth]{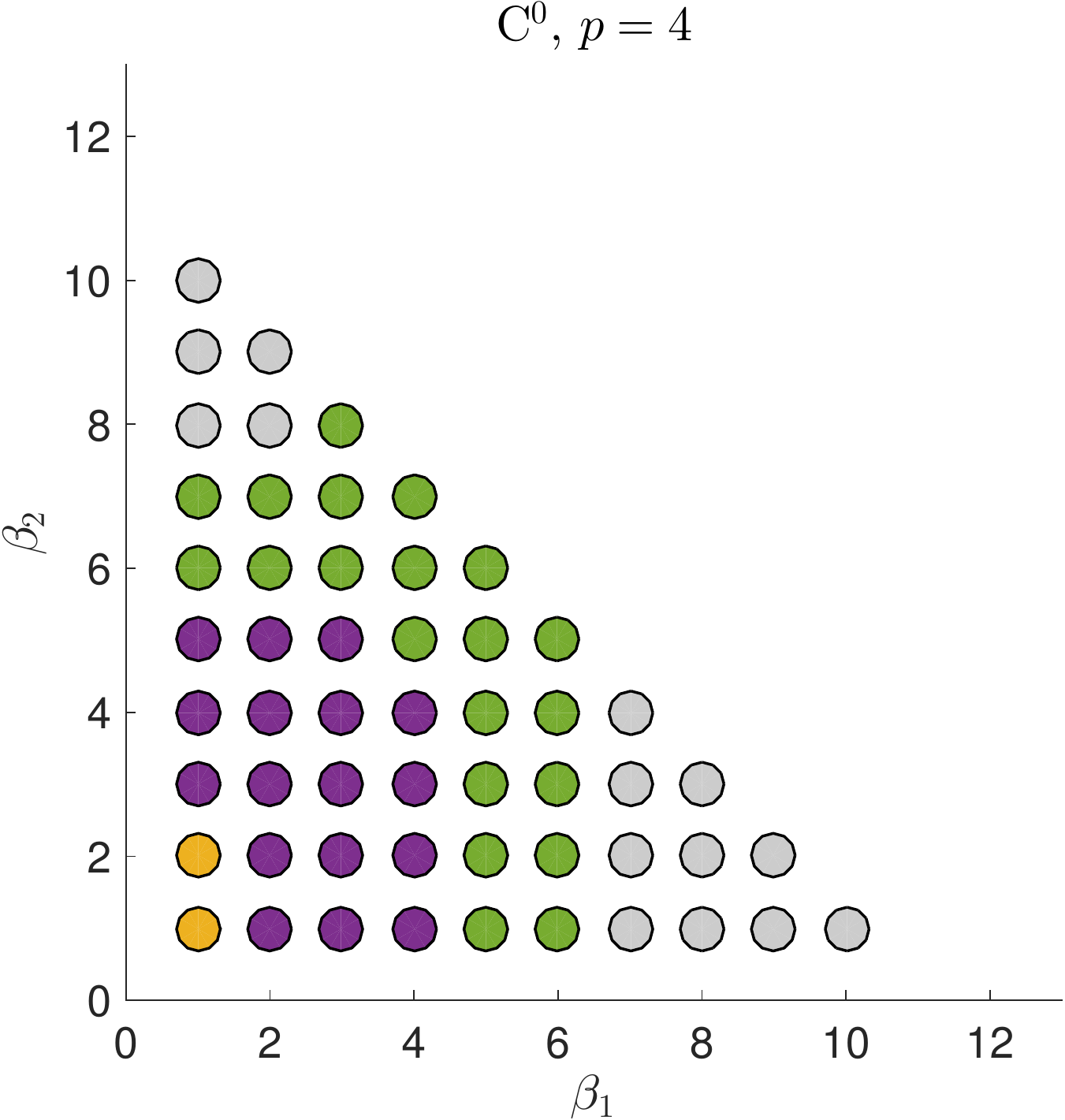}
  \includegraphics[width=\OptSetSize\linewidth]{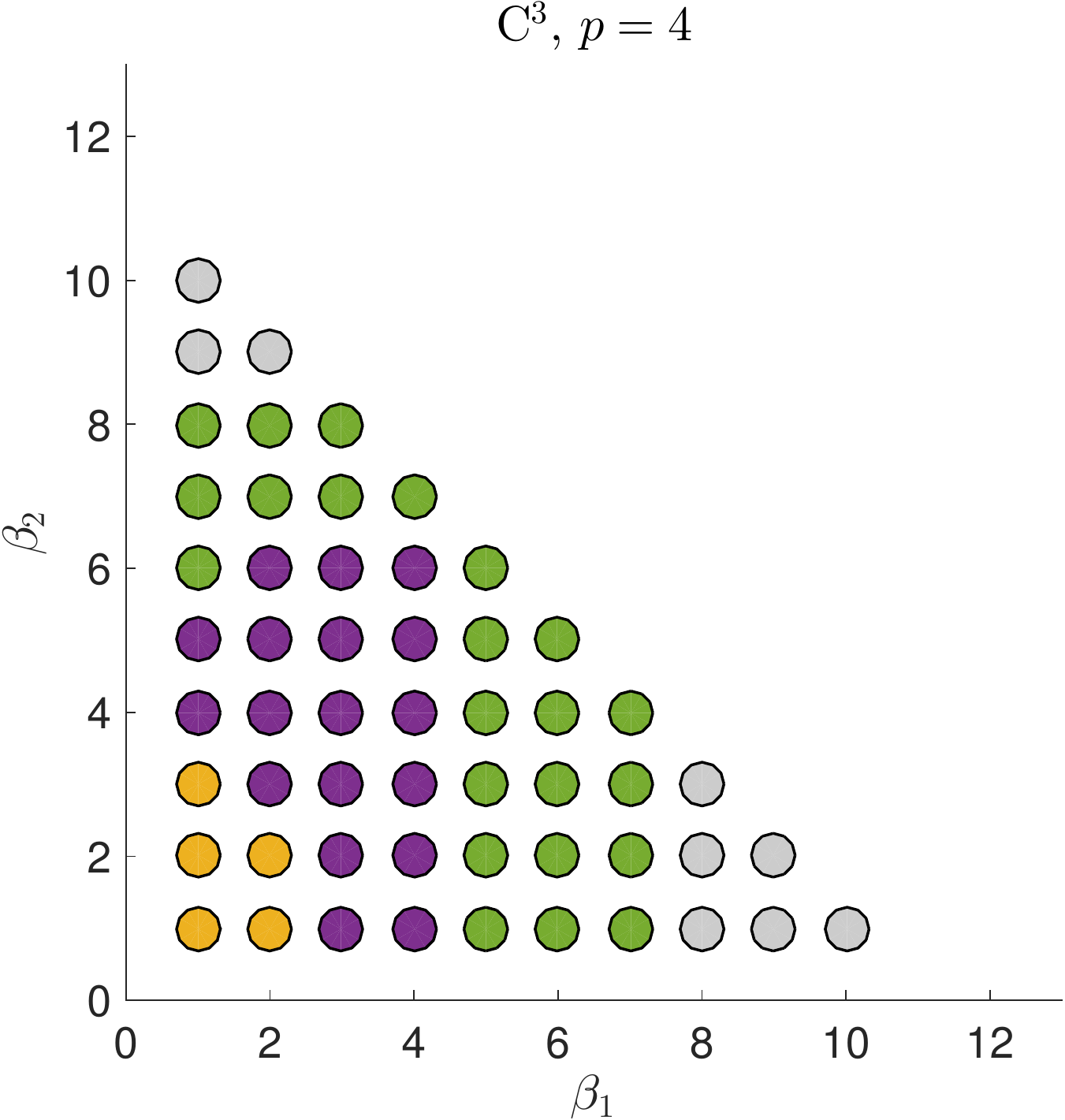}
\caption{Optimal sets for the quarter-of-annulus problem with low-regular solution.
    From top to bottom row: $p=2,3,4$. Left column: $C^0$ B-spline basis; Right column: maximal continuity
    B-splines basis.}\label{fig:opt-sets-ring-nonreg}
\end{figure}

If we compute and look at the profits for this problem,
we can indeed observe that the optimal sets are closer to a rectangle than a triangle, especially for $p=3,4$,
i.e., the method of choice for this problem should be the standard (i.e., tensor-based) IGA,
see Figure \ref{fig:opt-sets-ring-nonreg}.
It is also interesting to note that the profits decay more slowly along
the second parametric direction (i.e., the angular direction in the physical domain),
which therefore requires a finer discretization to be properly resolved, as one could have anticipated.
The exact shape of the optimal set is however difficult to describe
with a closed-form formula so that a \La{dimension-}adaptive scheme that iteratively chooses
the best operator $\Delta_{\bbeta}$ to be added to the sparse grid
should be devised if one wants to take advantage of these partial sets, see e.g. \cite{b.griebel:acta}.



\begin{figure}[tp]
  \centering
  \includegraphics[width=\ConvSize\linewidth]{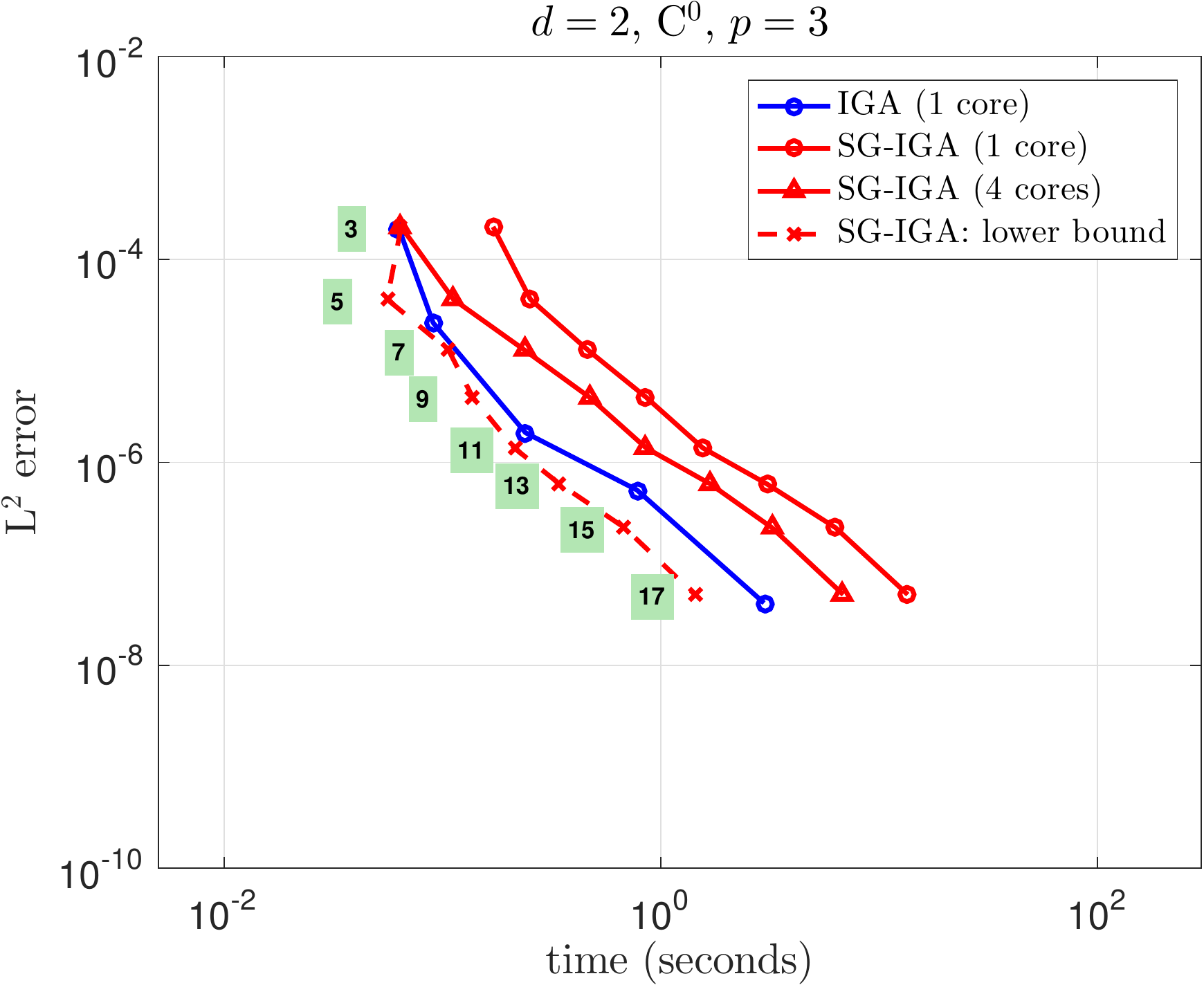}
  \includegraphics[width=\ConvSize\linewidth]{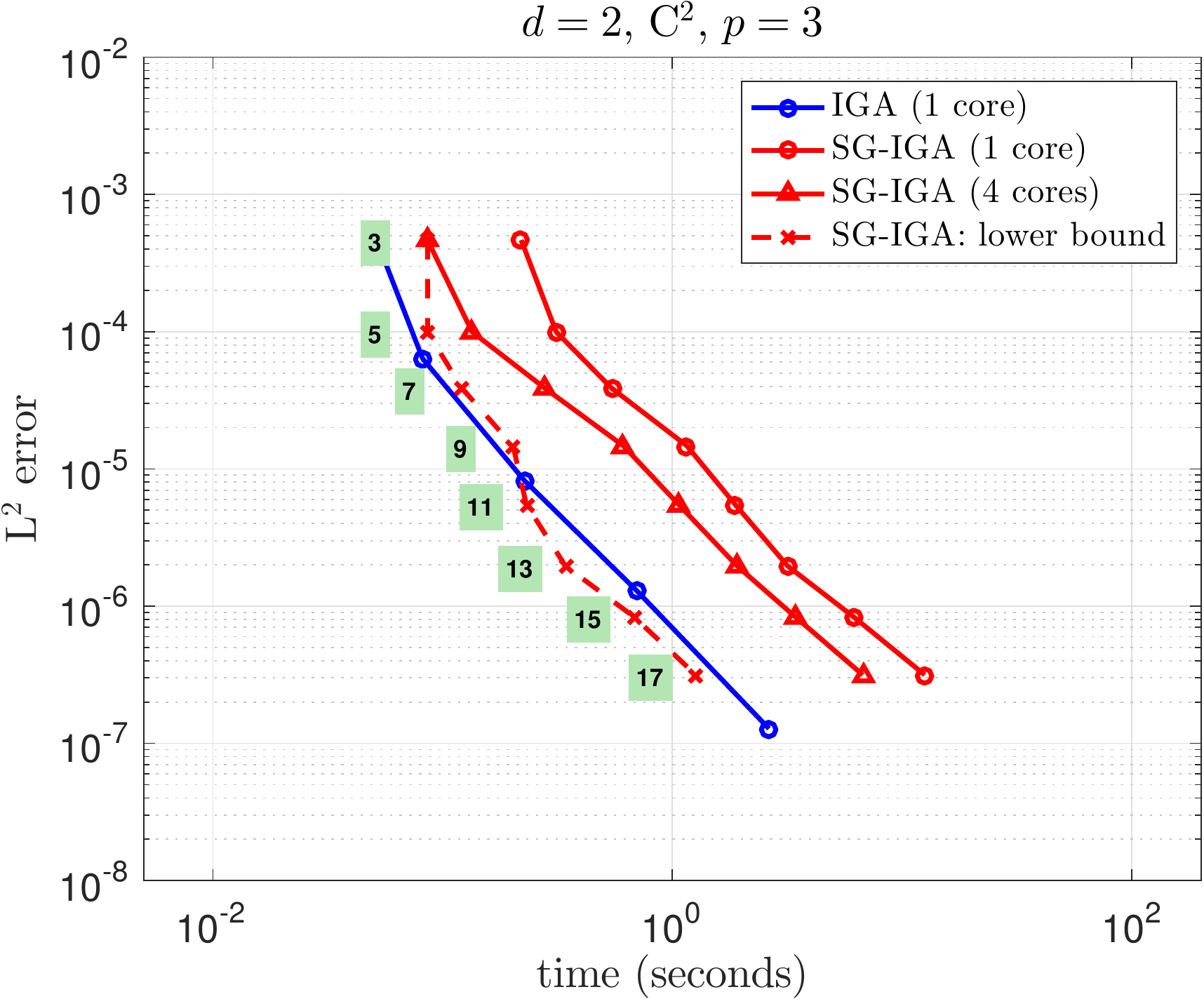} \\[\SpaceSize]
  \includegraphics[width=\ConvSize\linewidth]{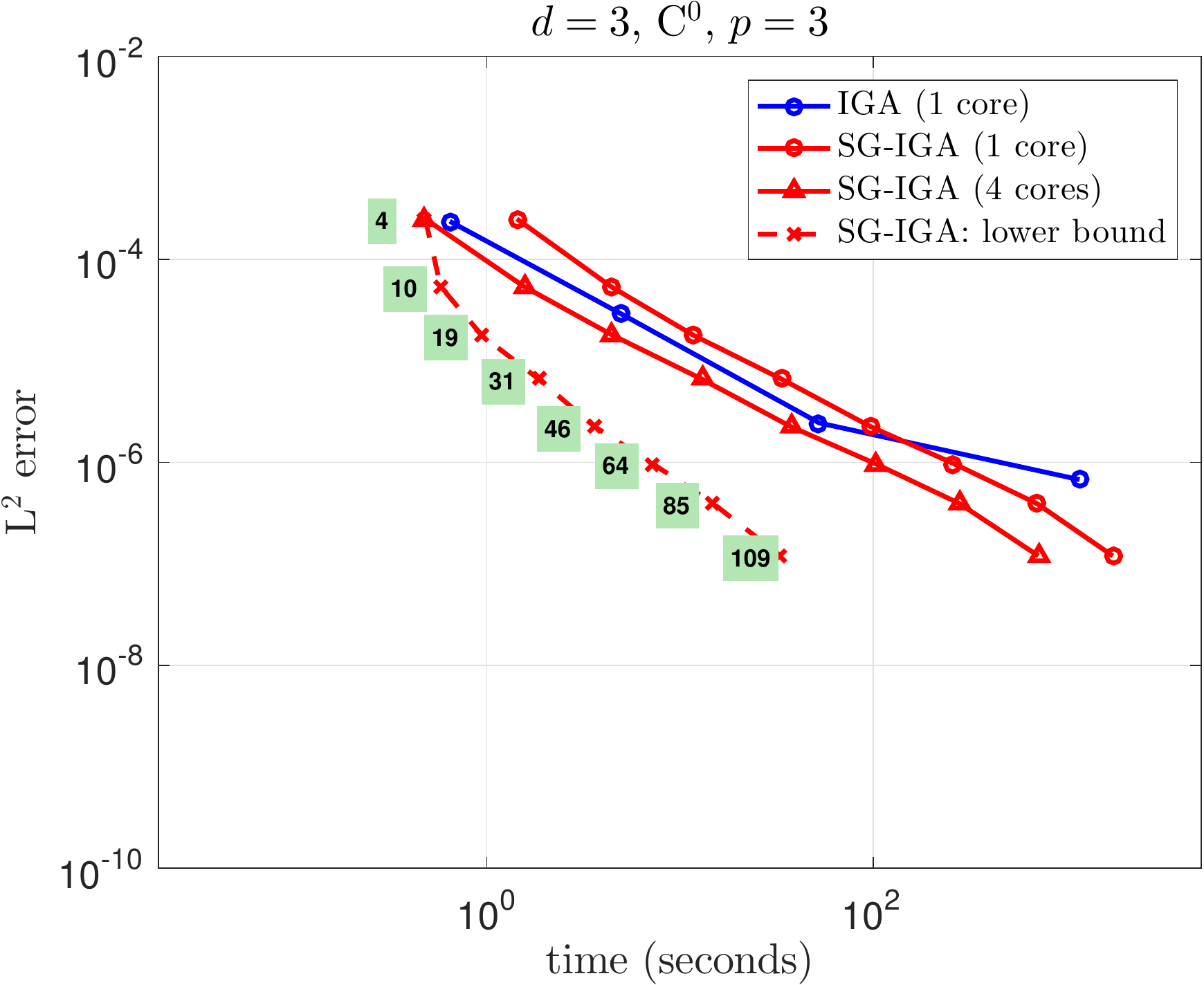}
  \includegraphics[width=\ConvSize\linewidth]{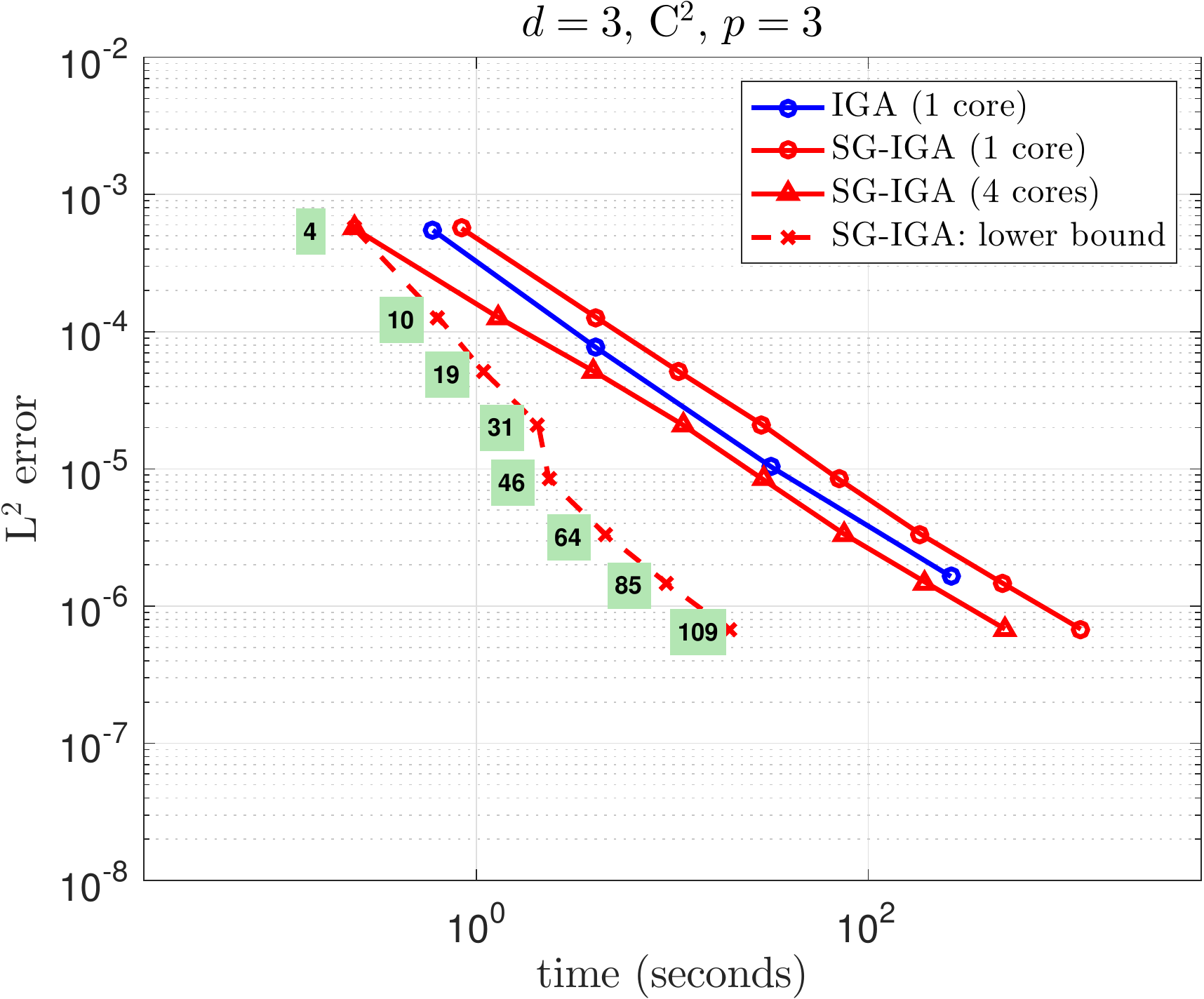}
  \caption{$L^2$ error vs. time for the quarter-of-annulus problem with low-regular solution.
    Here we fix $p=3$ and change $d$ (top row: $d=2$, bottom row: $d=3$)
    and the regularity of the B-splines basis \La{(left column: $C^0$, right column: $C^2$)}. 
    The dashed line is the lower bound that 
    can be achieved if the number of available cores is at least equal to the number of components 
    of the combination technique for each level, given by the numbers in green boxes.}\label{fig:err-vs-time-non-reg-ring}

  \bigskip

  \includegraphics[width=\ConvSize\linewidth]{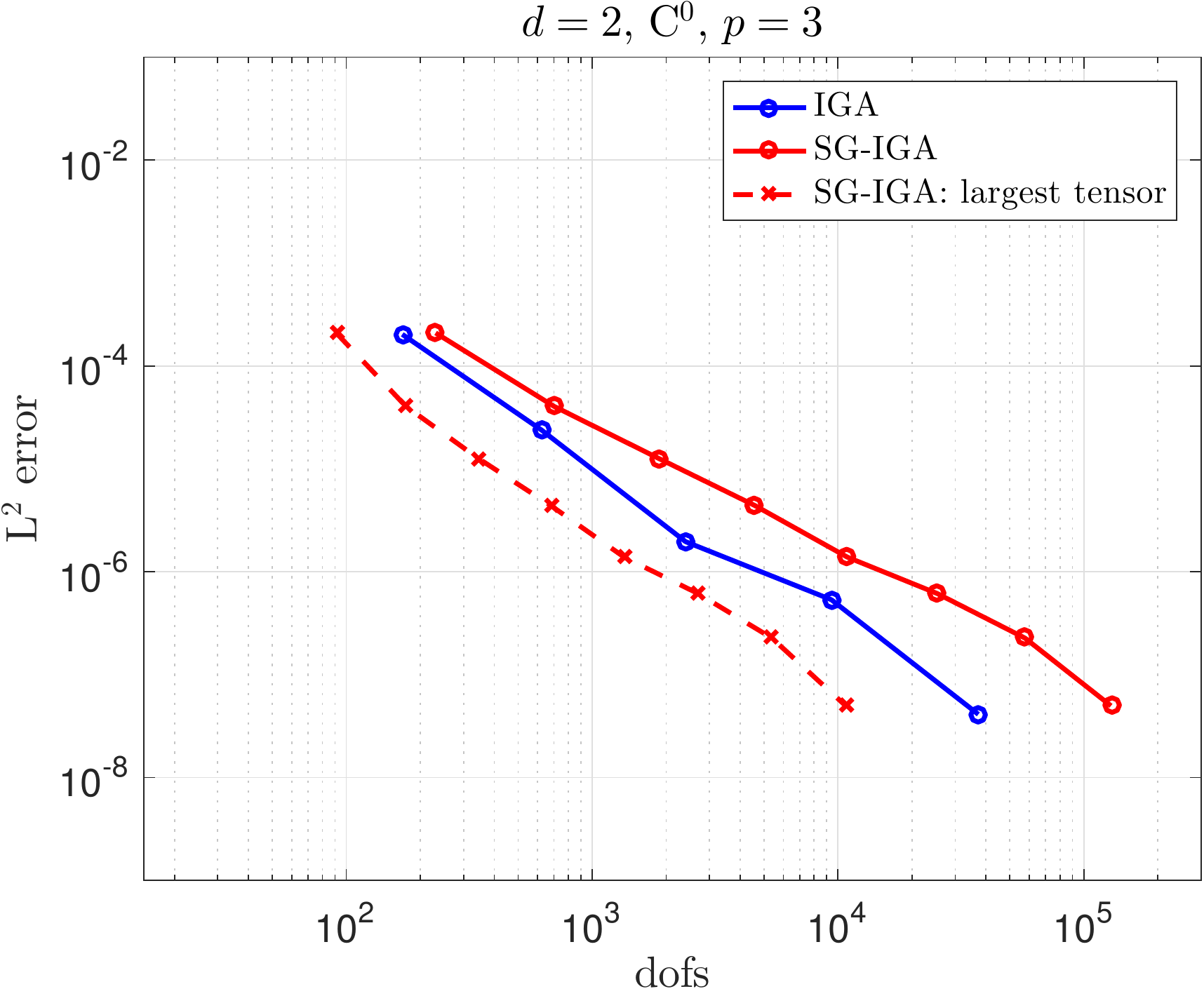}
  \includegraphics[width=\ConvSize\linewidth]{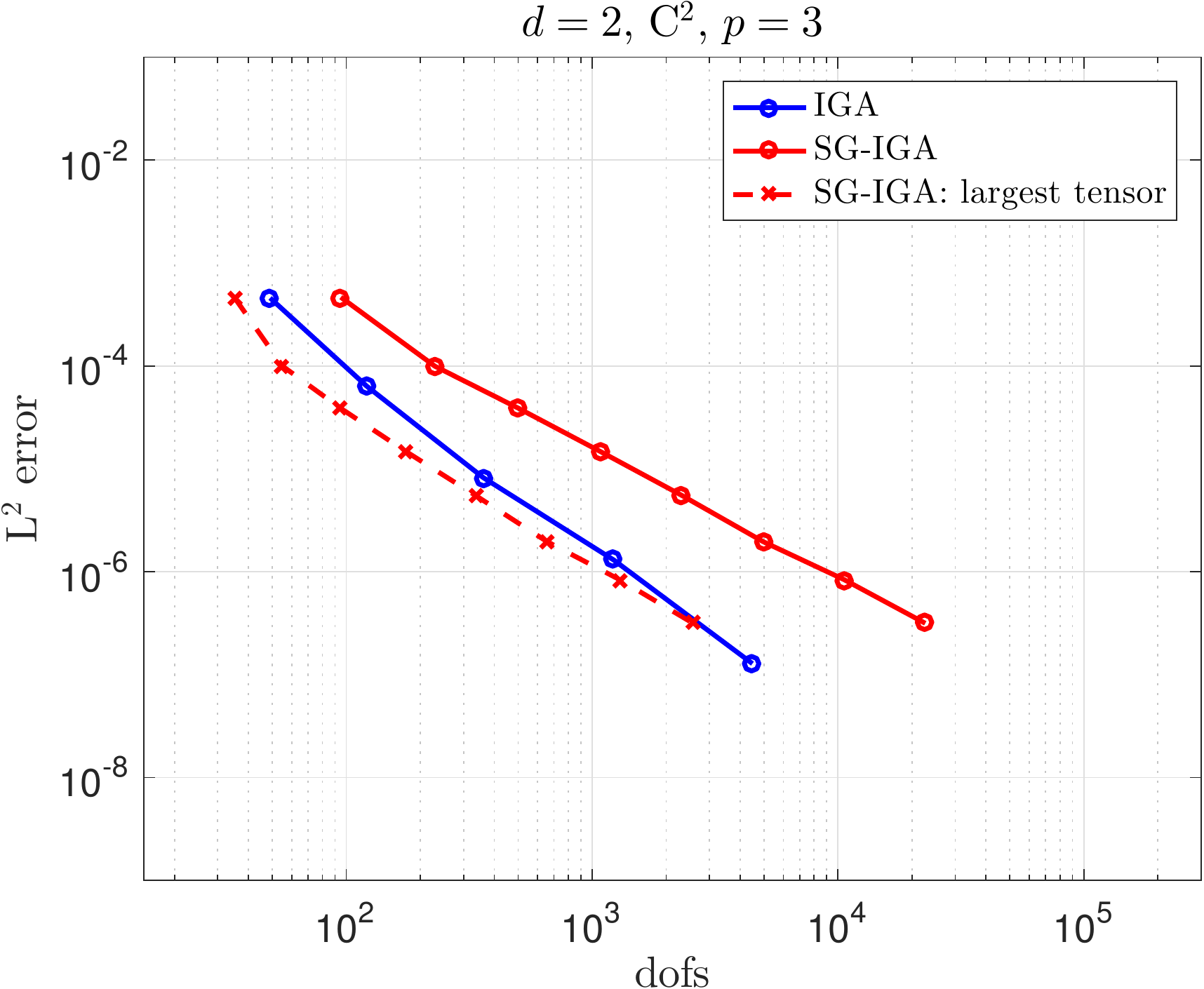} \\[\SpaceSize]
  \includegraphics[width=\ConvSize\linewidth]{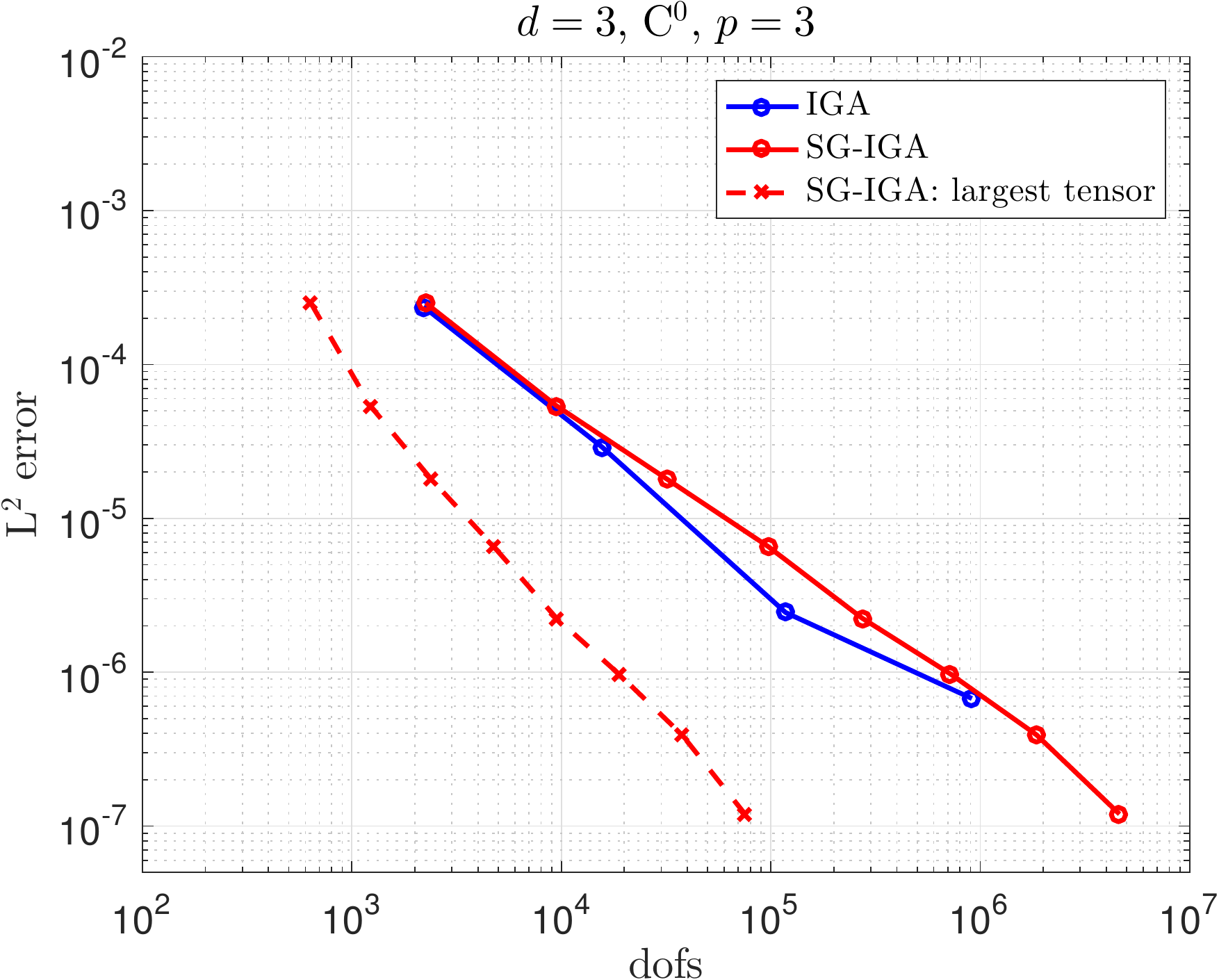}
  \includegraphics[width=\ConvSize\linewidth]{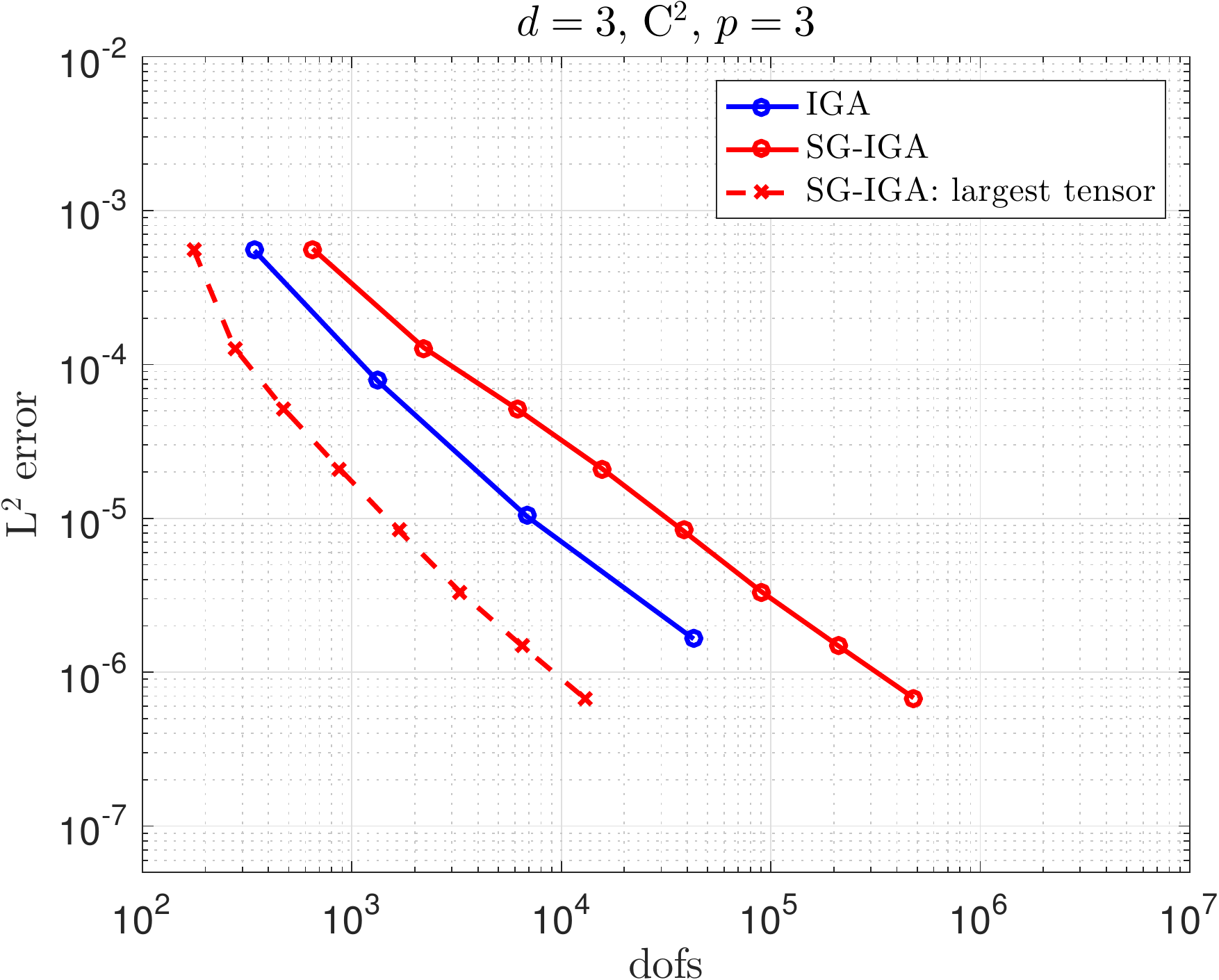}
  \caption{$L^2$ error vs. dofs for the quarter-of-annulus problem with low-regular solution. 
    Here we fix $p=3$ and change $d$ (top row: $d=2$, bottom row: $d=3$) and the regularity of the B-splines basis
    \La{(left column: $C^0$, right column: $C^2$)}.
    In each figure, the dashed line represents the error of SG-IGA vs. the number
    of the degrees-of-freedom for the largest tensor grid in the combination technique.}\label{fig:err-vs-dofs-non-reg-ring}
\end{figure}

Figures \ref{fig:err-vs-time-non-reg-ring} and \ref{fig:err-vs-dofs-non-reg-ring}
show the decay of the error with respect to time and degrees-of-freedom, again focusing on the case $p=3$ only.
These results confirm the underperformance of SG-IGA with respect to IGA in this case:
compared to the problem with regular solution, cf. Figures \ref{fig:err-vs-time-reg-ring} and \ref{fig:err-vs-dofs-reg-ring},
the rate of SG-IGA with respect to both time and degrees-of-freedom is now nearly identical to the rate of IGA
and a larger number of cores is needed to obtain SG-IGA computational times close to IGA times
(more than 4 cores in $d=2$, more than 1 core in $d=3$).


\subsection{Radical meshes in the parametric domain}

\begin{figure}[tp]
  \centering
  \includegraphics[width=0.32\linewidth]{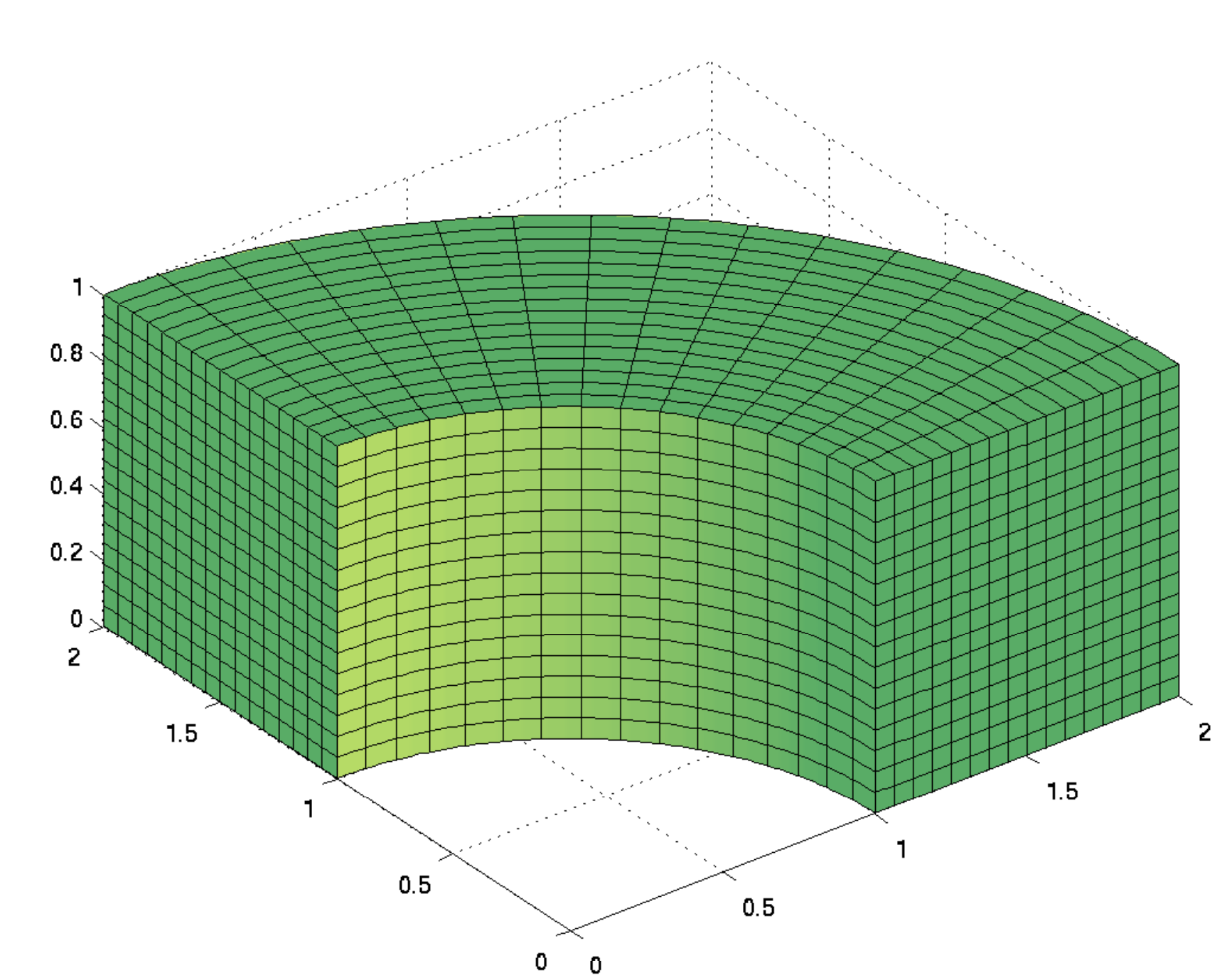}
  \includegraphics[width=0.32\linewidth]{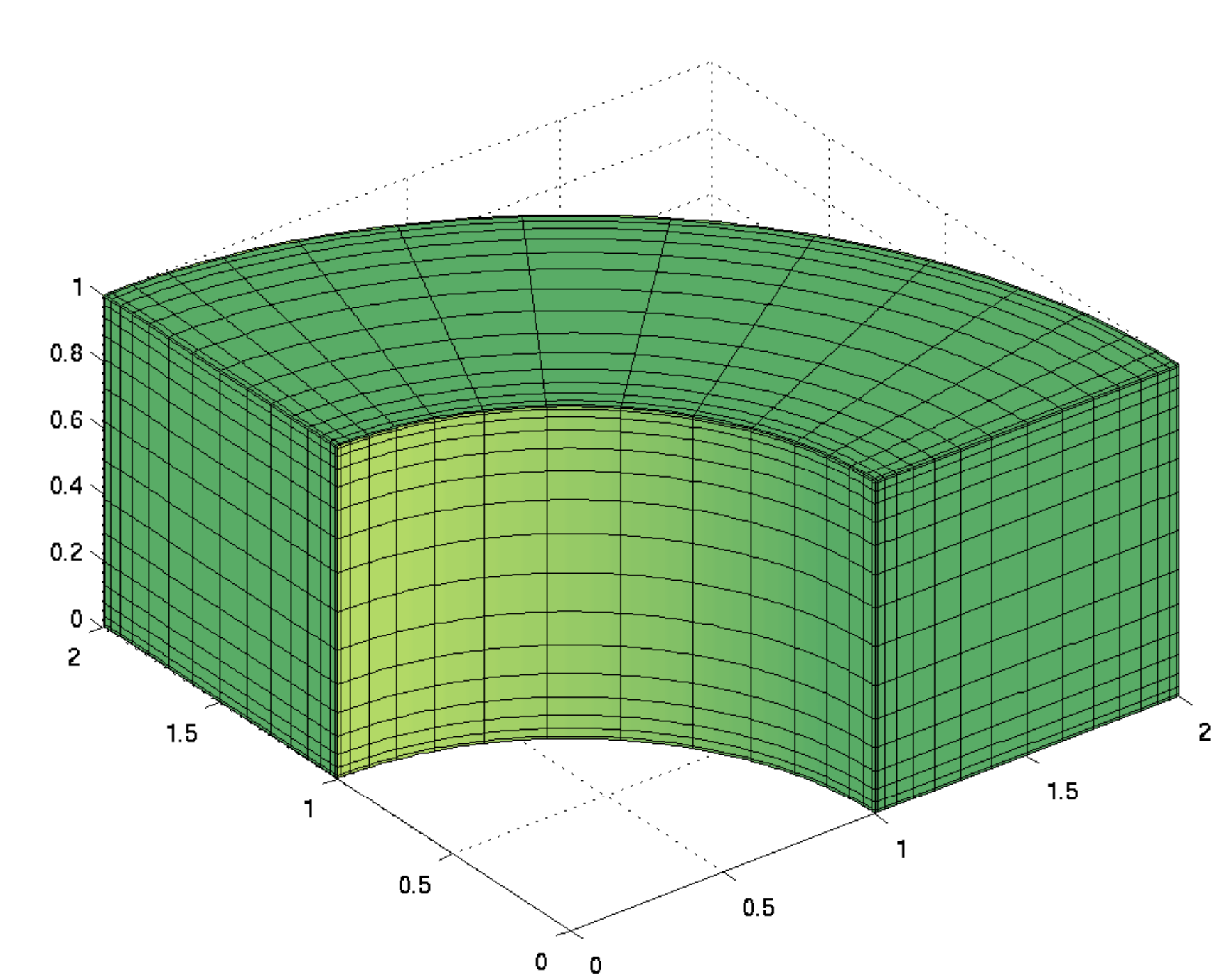}
  \includegraphics[width=0.32\linewidth]{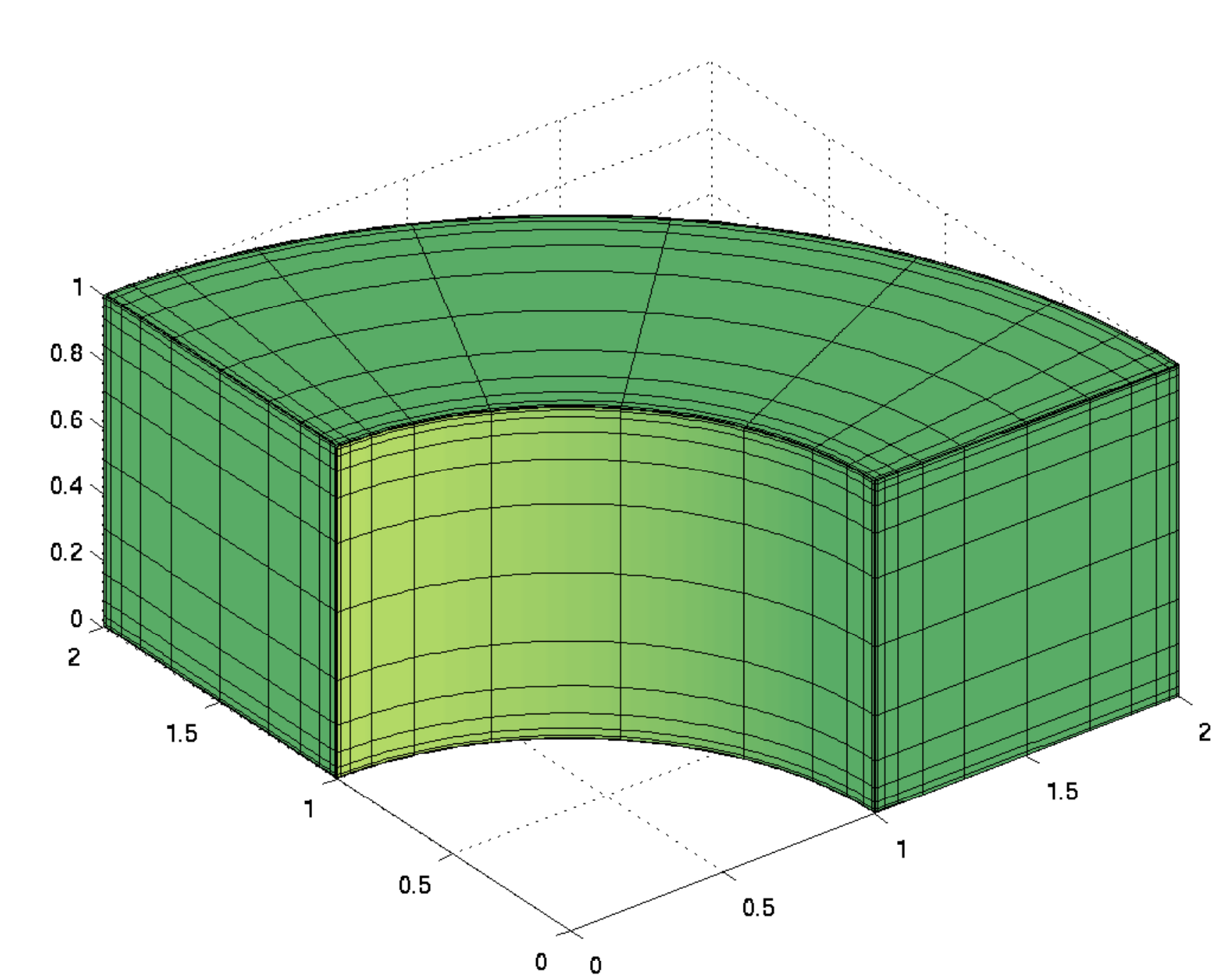}
  \caption{Knotlines of the quarter-of-annulus, $d=3$, for radical meshes with $\gamma=1,2,4$ (from left to right).
  Observe that the number of elements is identical in the three meshes.}
  \label{fig:radical-meshes}
\end{figure}

The solutions for the previous problem lack of (mixed) Sobolev regularity because of their singular
behavior on the corners ($d=2$ and $d=3$) and edges ($d=3$) of the domain $\Omega$
(see \cite{Grisvard}). This causes the slower convergence rate for both
IGA and SG-IGA. To improve convergence in this situation,
we adopt non-uniform radical meshes, see \cite{Babuska_book,Beirao_Cho_Sangalli};
see also \cite{garcke:graded,griebel.thurner:graded} 
for previous use of graded meshes in a sparse grids context. 
To this end, we compute the position of the $i$-th node on the $j$-th parametric dimension as $\xi_{i,j}=f(t_i)$, 
where $t_i$ are knots of an equispaced grid over $[0,1]$ and 
$f$ is defined as follows, depending on a parameter $\gamma$: 
\[
  f(t) = 
  \begin{cases}
   \displaystyle \frac{t^\gamma}{1/2^{\gamma-1}} & \text{ if } t\leq \frac{1}{2}, \\[8pt]
   \displaystyle  1 - \frac{(1-t)^{\gamma}}{1/2^{\gamma-1}} & \text{ if } t > \frac{1}{2}.
  \end{cases}
\]
Note that for $\gamma=1$ one recovers the equispaced mesh. We show some radical meshes
for the $d=3$ version of the quarter-of-annulus domain in Figure \ref{fig:radical-meshes}.

\begin{figure}[tp]
  \centering
  \includegraphics[width=\GammaSize\linewidth]{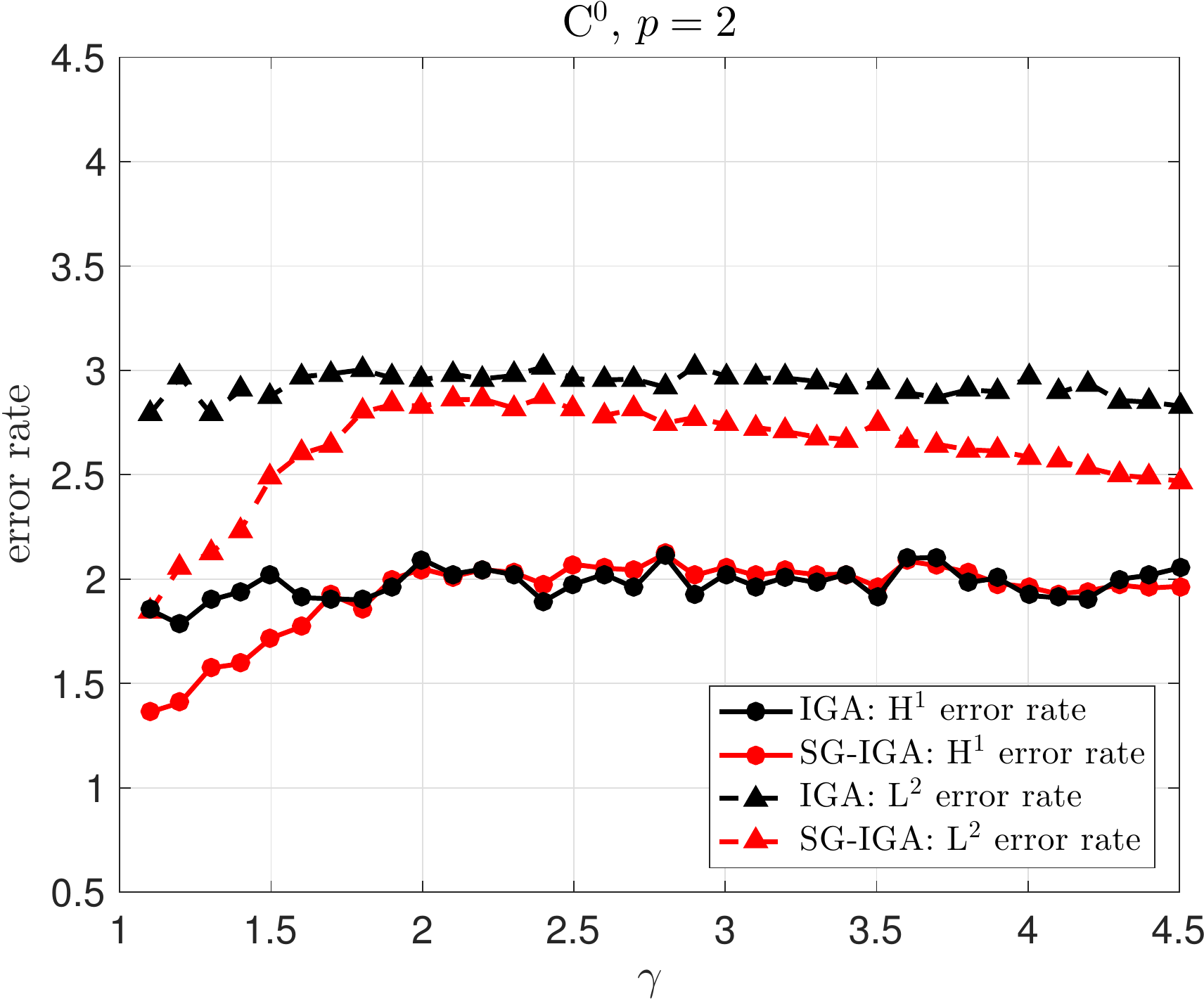} 
  \includegraphics[width=\GammaSize\linewidth]{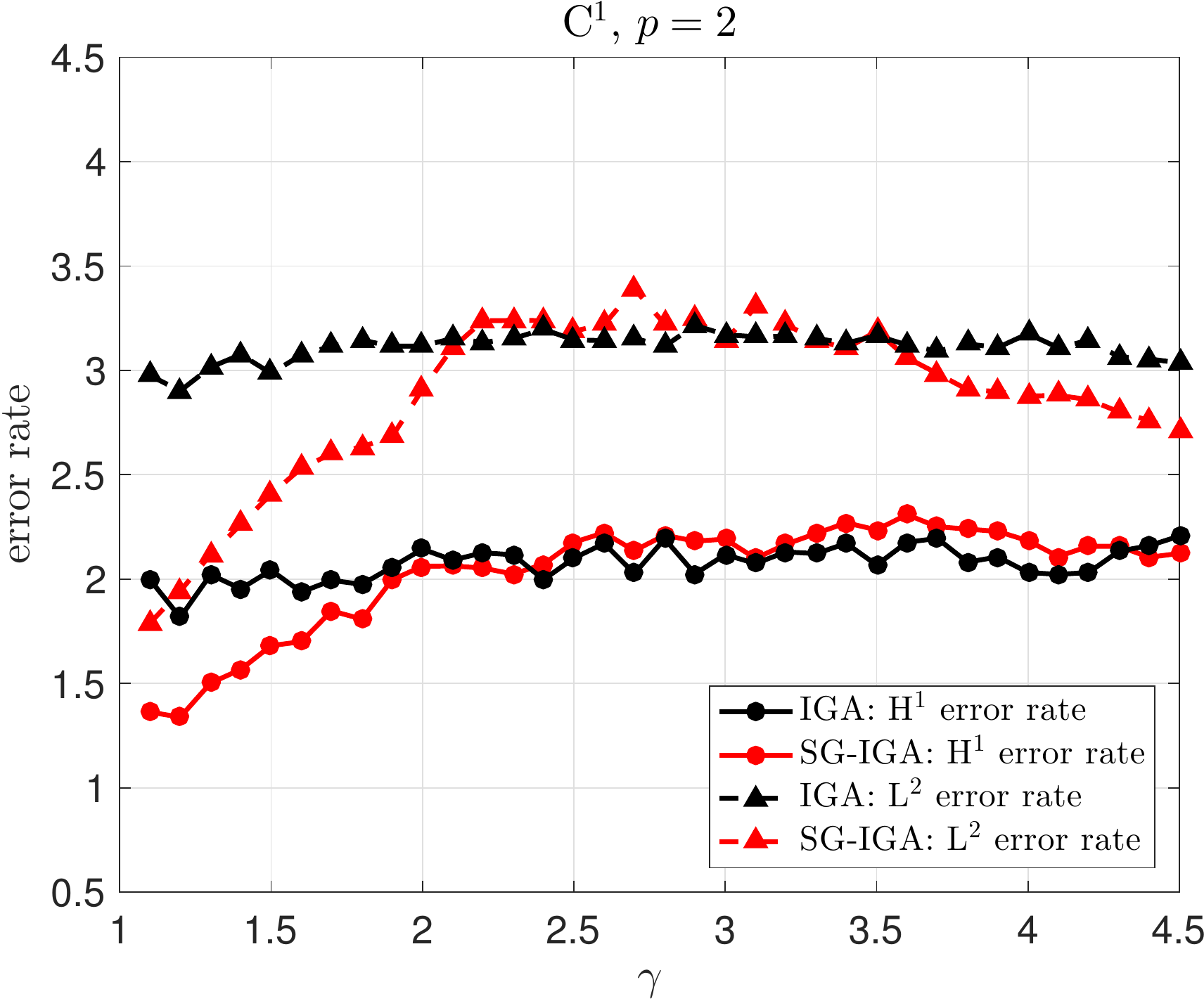}  \\[\SpaceSize]
  \includegraphics[width=\GammaSize\linewidth]{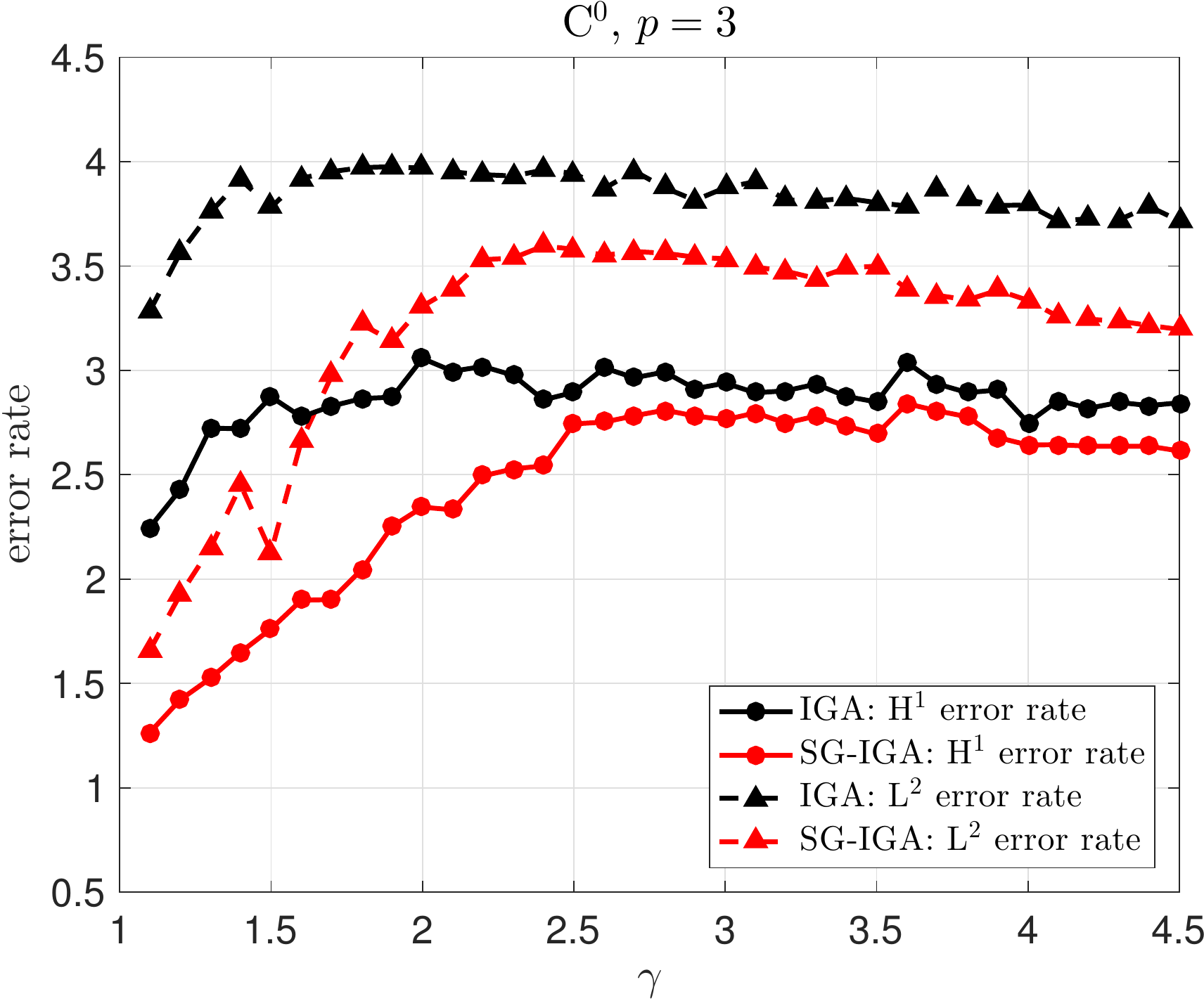} 
  \includegraphics[width=\GammaSize\linewidth]{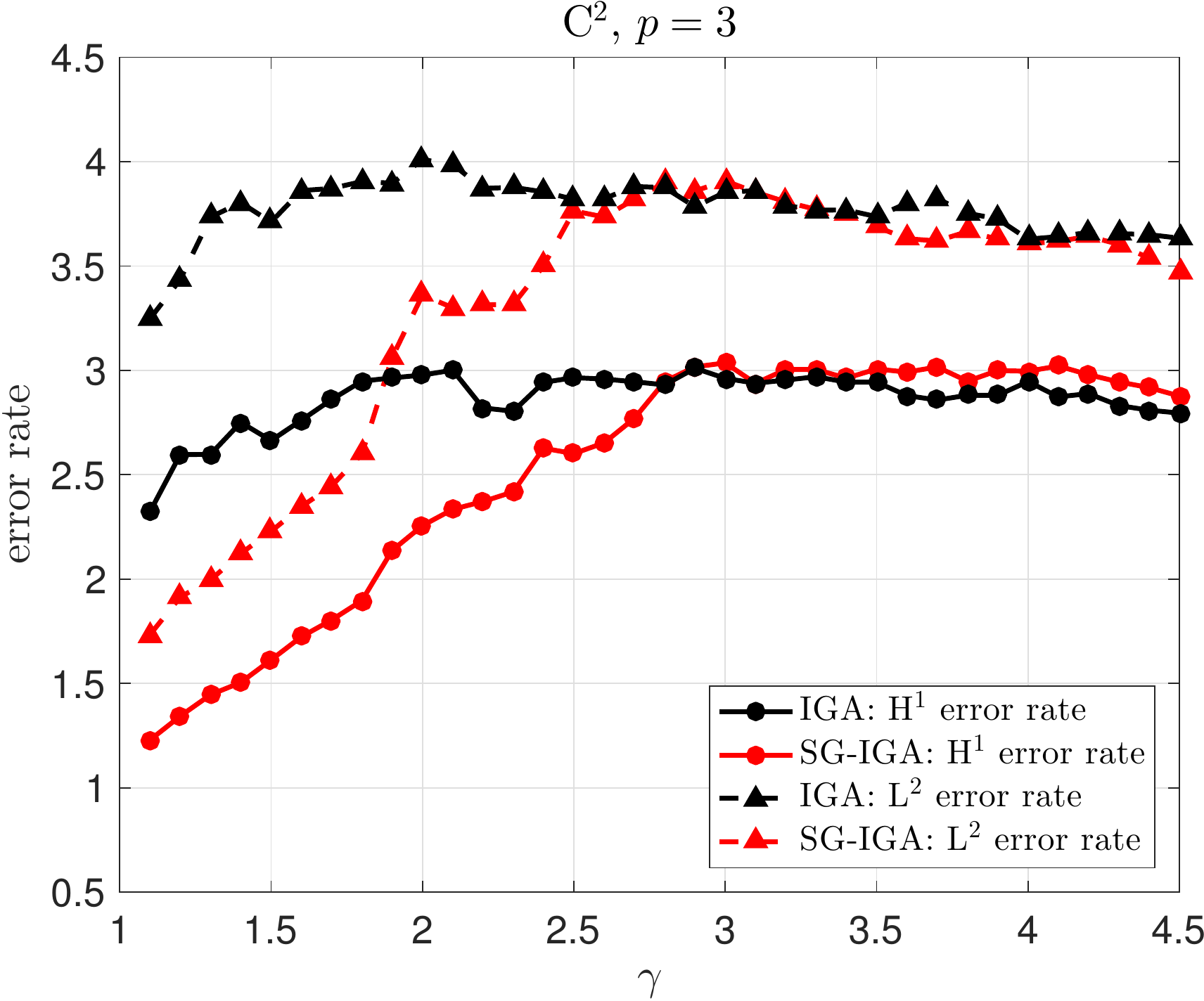} \\[\SpaceSize]
  \includegraphics[width=\GammaSize\linewidth]{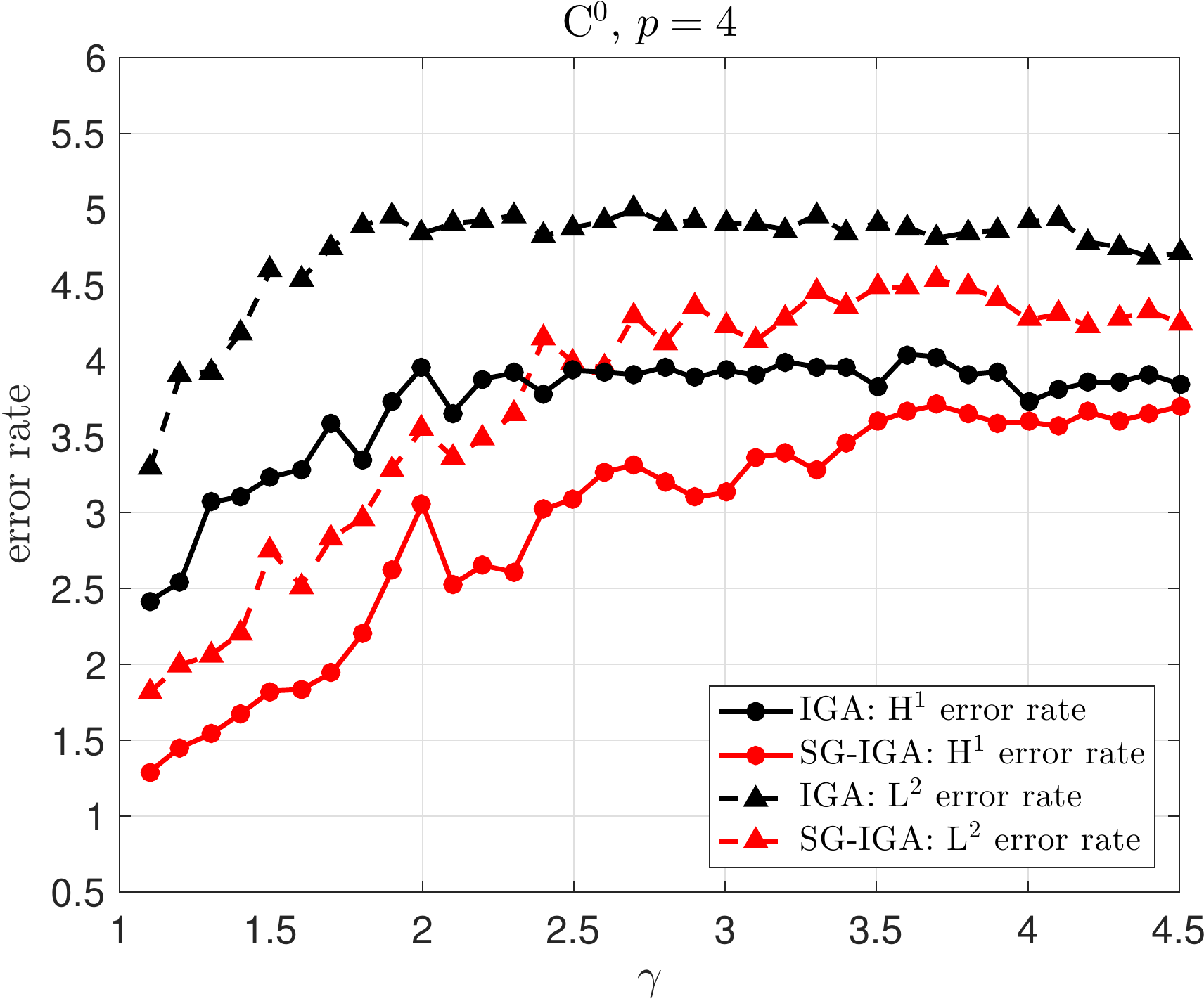} 
  \includegraphics[width=\GammaSize\linewidth]{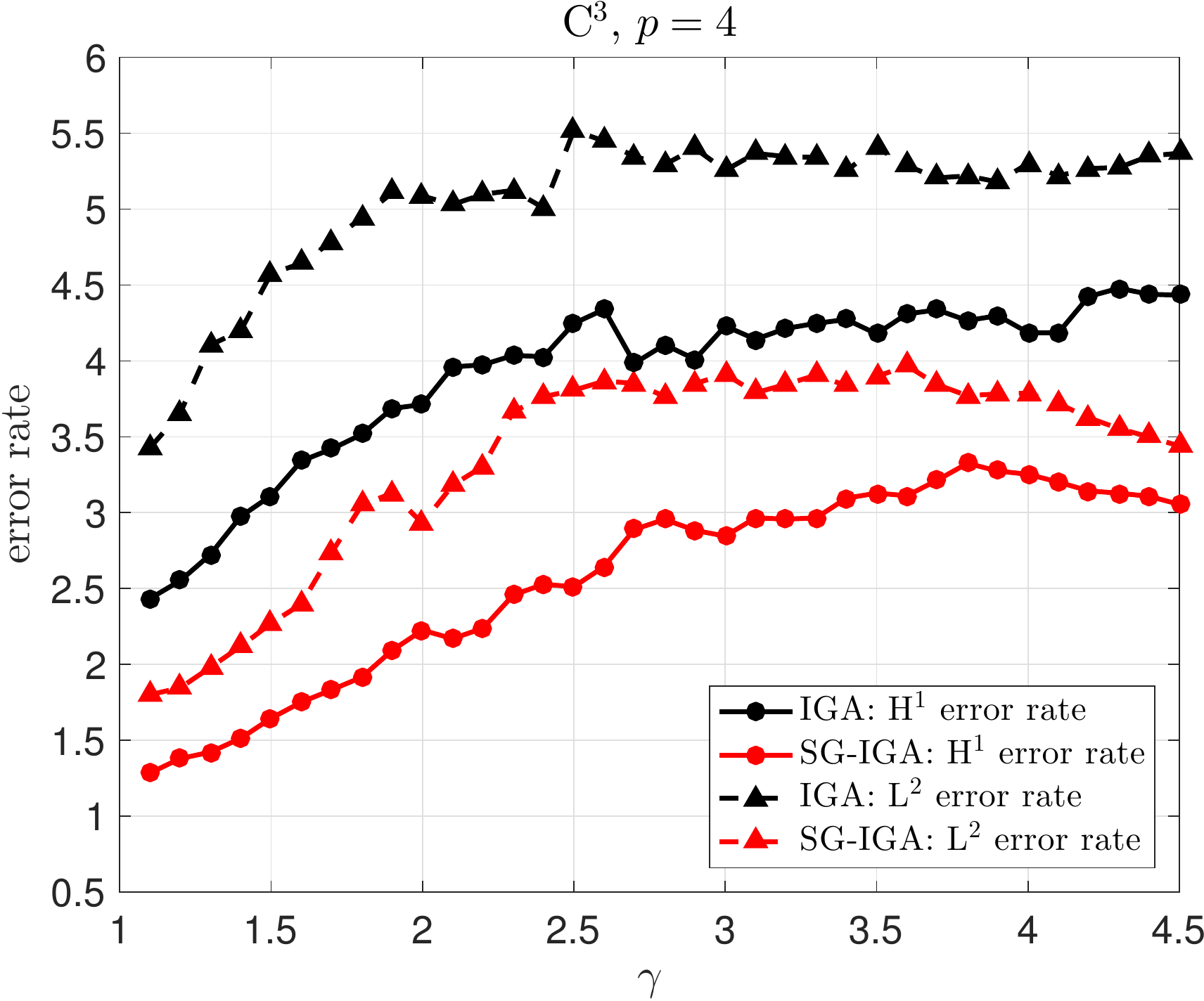}  \\
  \caption{Measured convergence rate for increasing values of $\gamma$ for the quarter 
    of annulus problem ($d=2$). From top to bottom: $p=2$ to $p=4$. 
    The left column shows $C^0$ B-splines, the right column $C^{p-1}$ B-splines
    (maximal continuity).}
  \label{fig:rates_vs_gamma-ring}
\end{figure}

\begin{figure}[tp]
  \centering
  \includegraphics[width=\OptSetSize\linewidth]{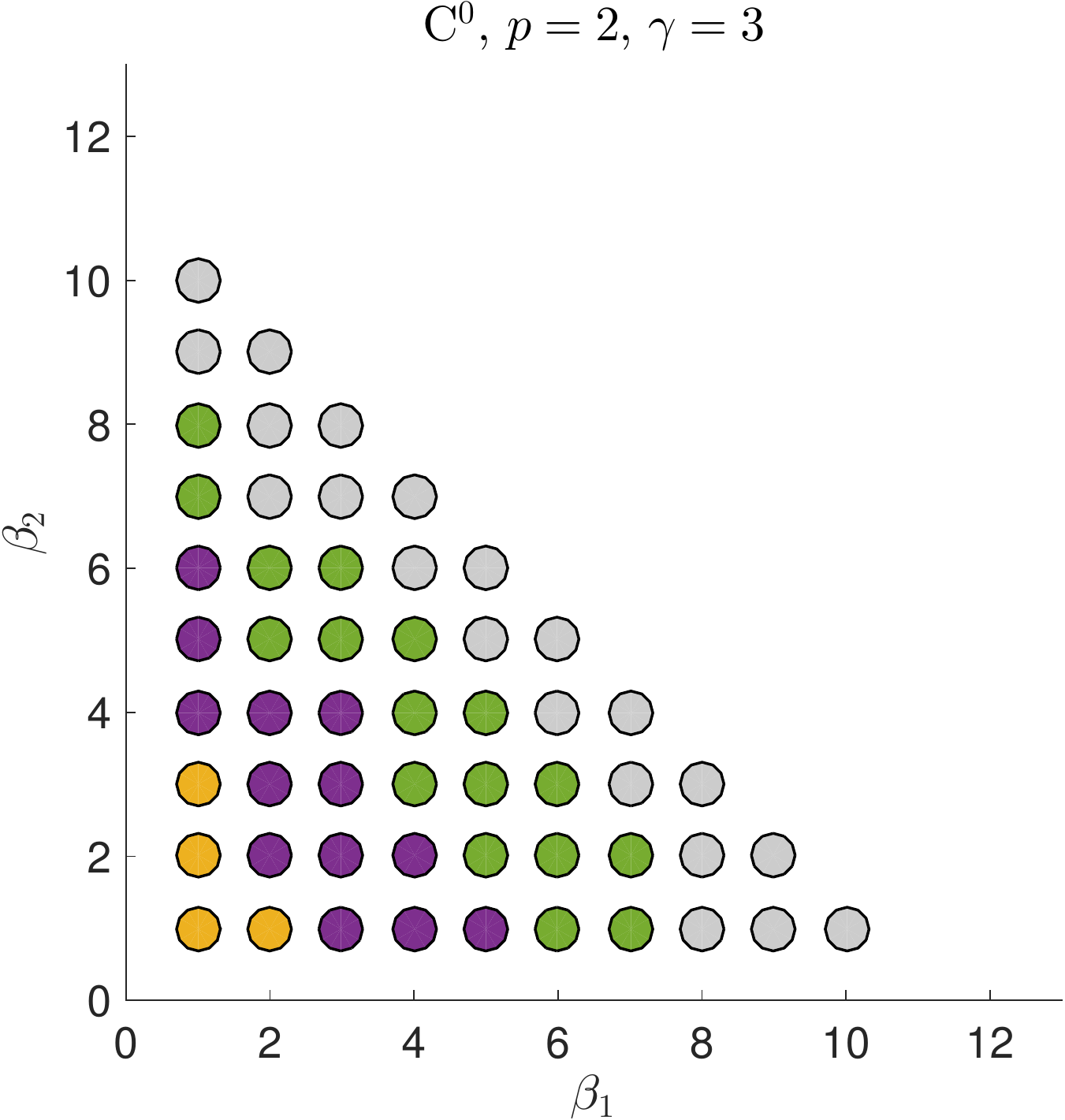}
  \includegraphics[width=\OptSetSize\linewidth]{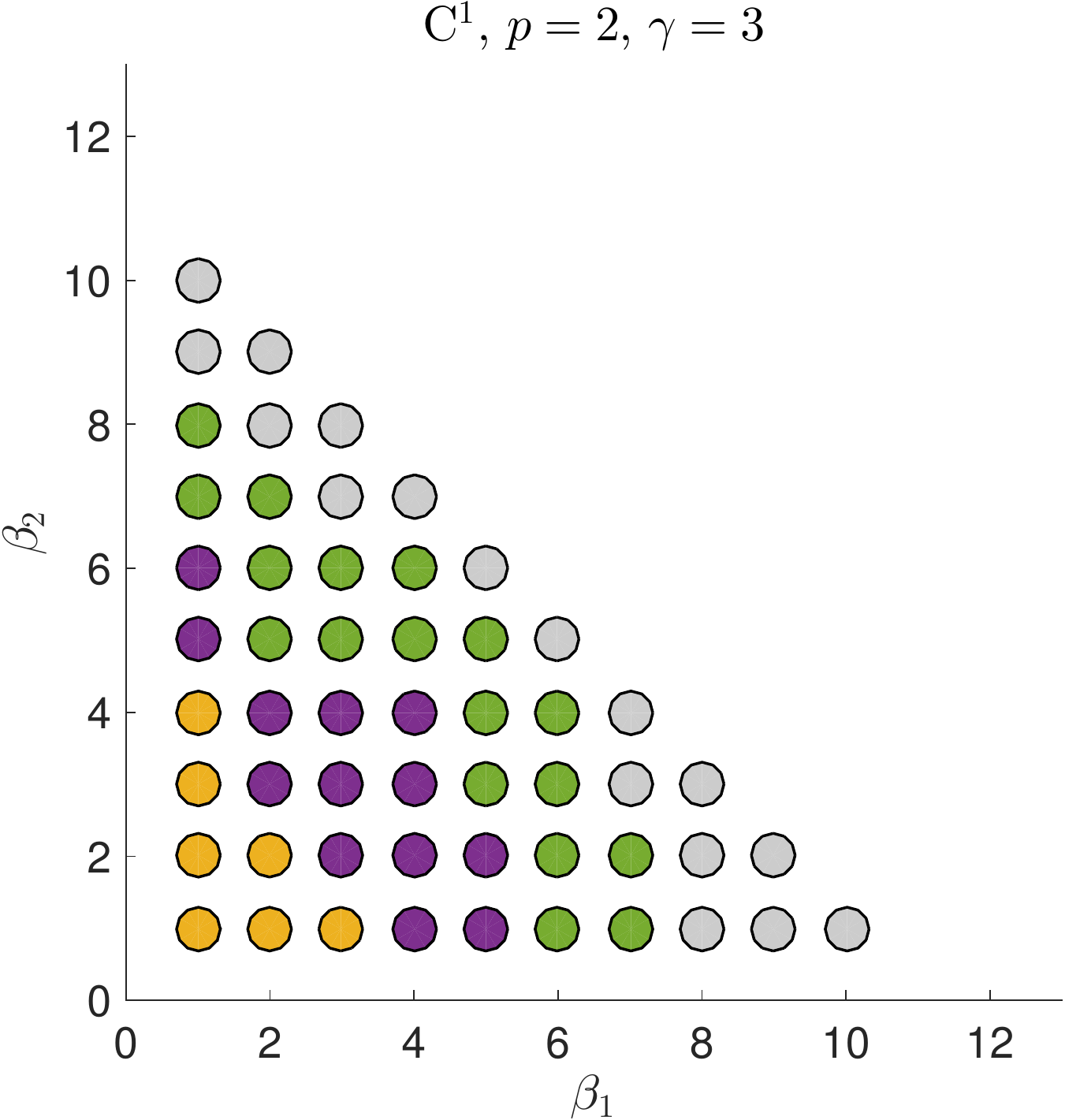}\\[\SpaceSize]
  \includegraphics[width=\OptSetSize\linewidth]{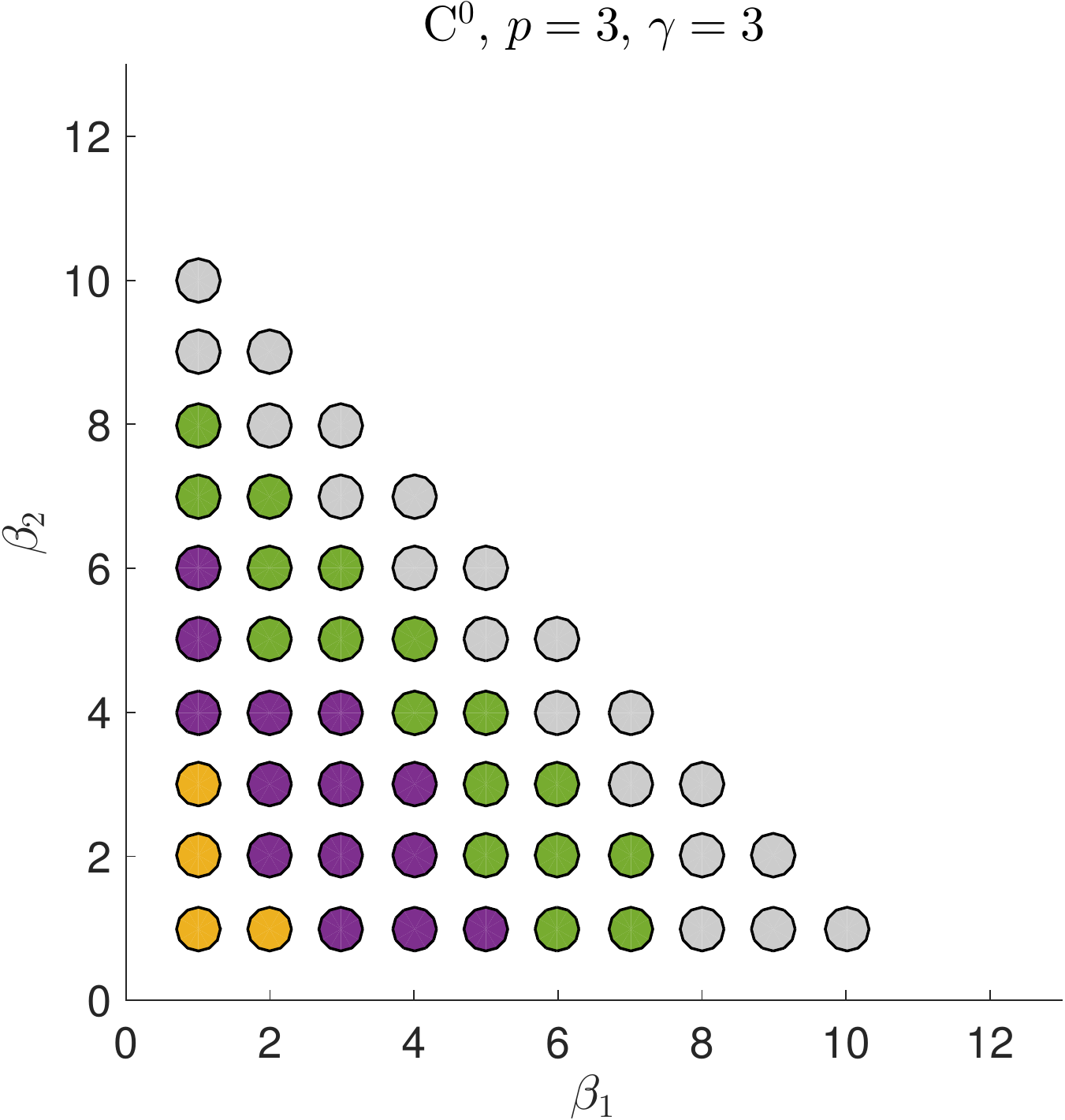}
  \includegraphics[width=\OptSetSize\linewidth]{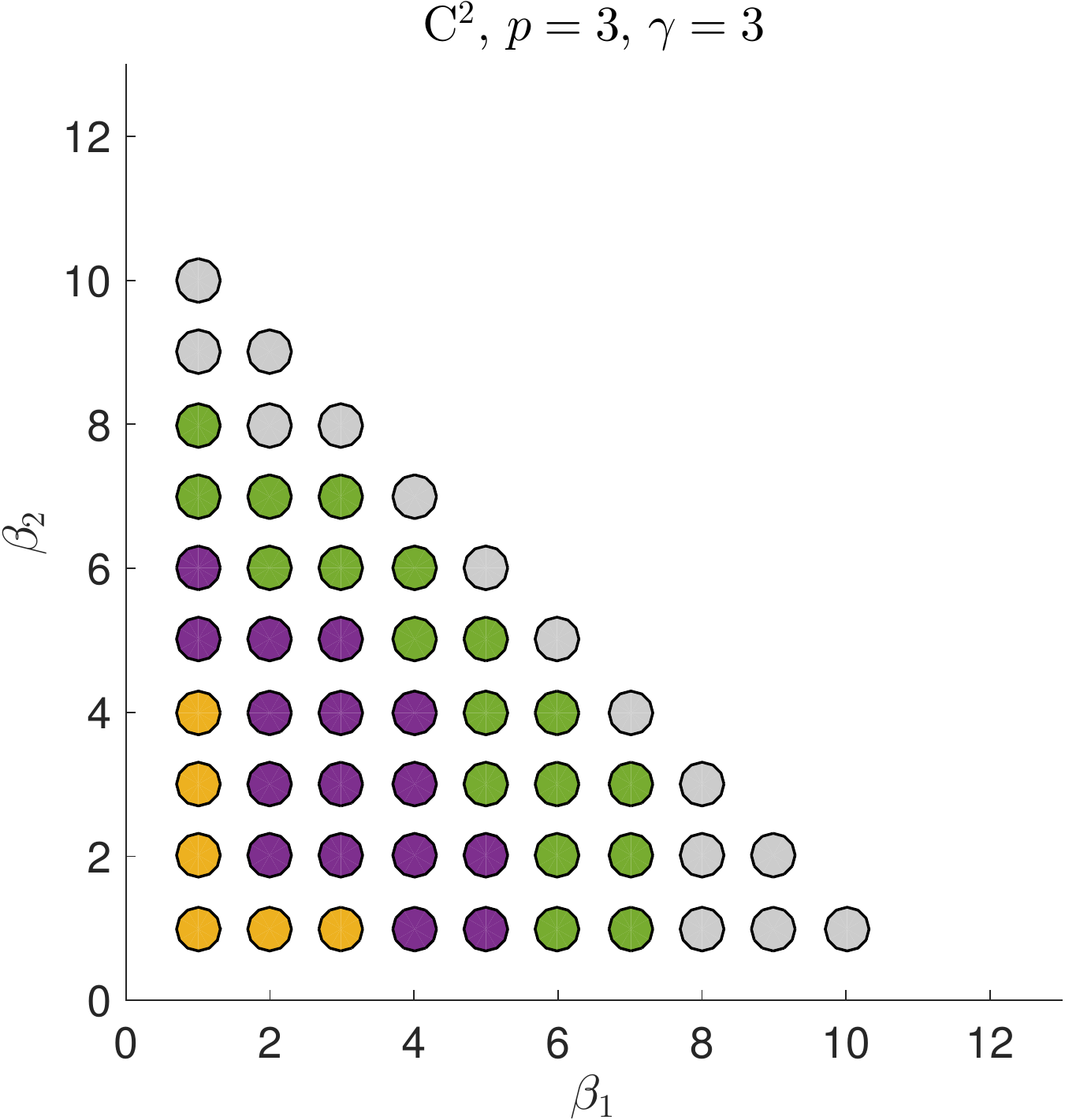}\\[\SpaceSize]
  \includegraphics[width=\OptSetSize\linewidth]{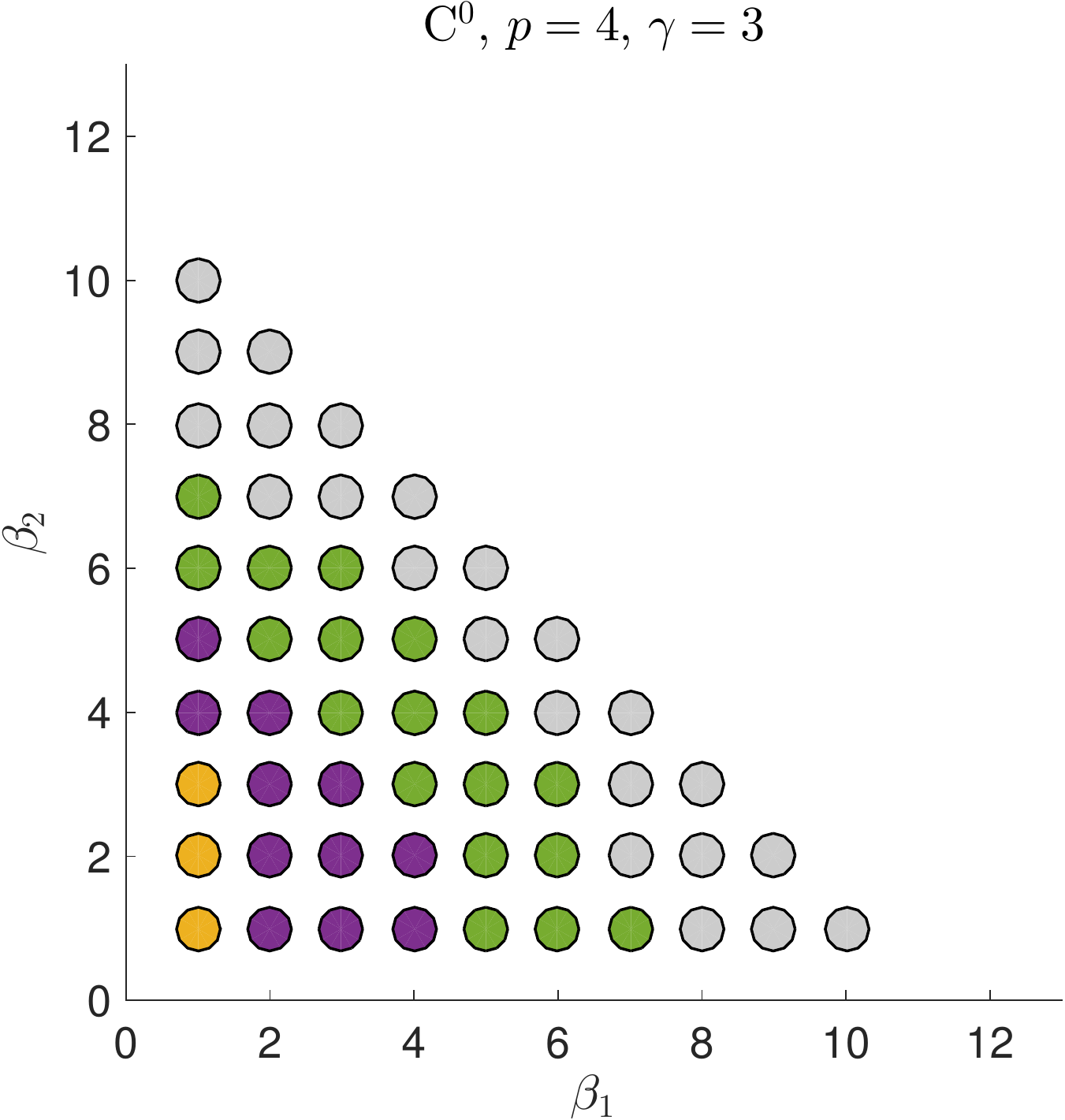}
  \includegraphics[width=\OptSetSize\linewidth]{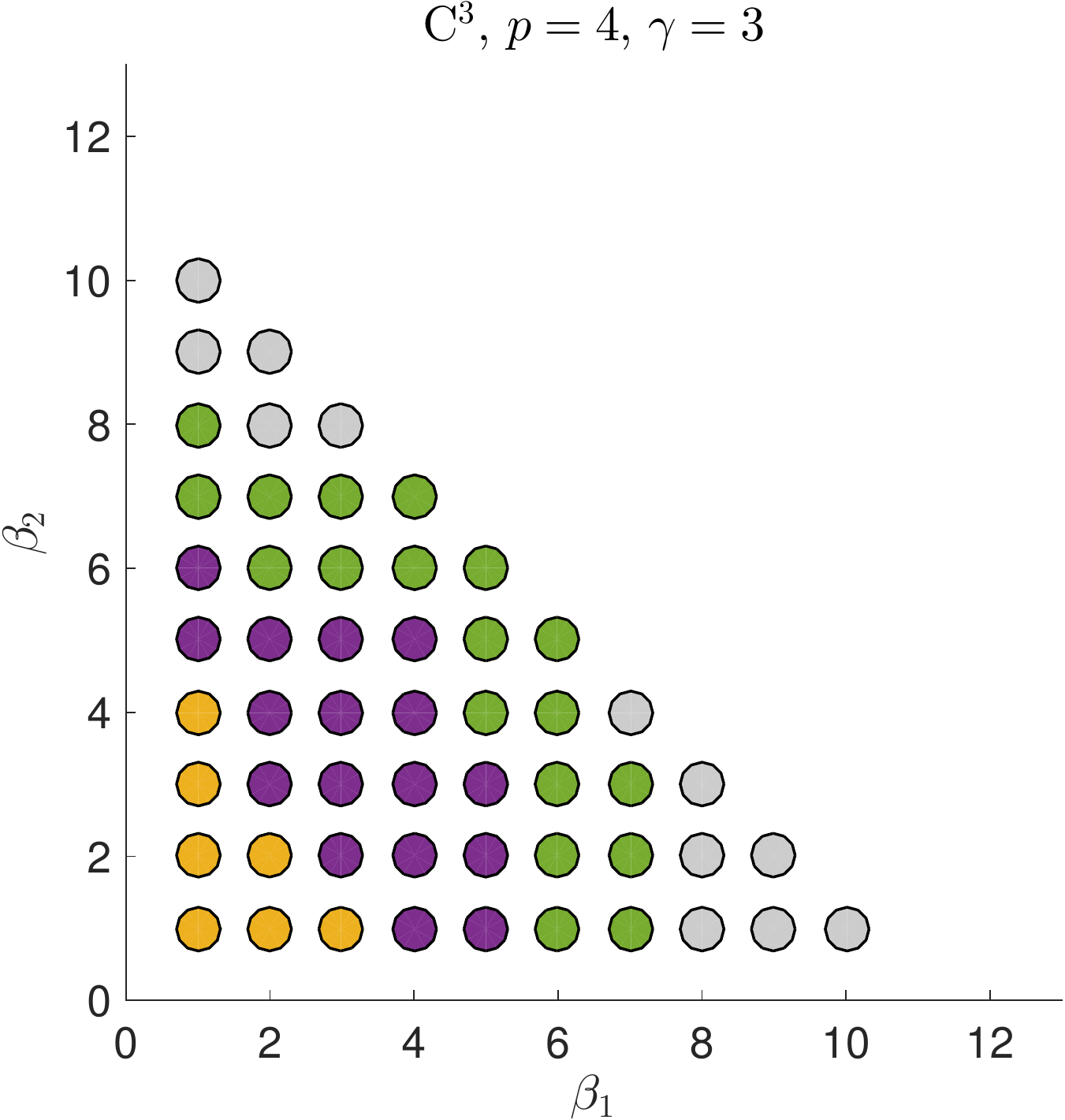}
  \caption{Optimal sets for the quarter-of-annulus problem with low-regular solution, for radical meshes with $\gamma=3$.
    From top to bottom row: $p=2,3,4$. Left column: $C^0$ B-spline basis; Right column: maximal continuity
    B-splines basis.}
  \label{fig:err-vs-time-non-reg-gamma3-ring-qoi} 
\end{figure}

In Figure \ref{fig:rates_vs_gamma-ring} we show the value of the rate for
the $H^1$ and $L^2$ error for IGA and SG-IGA applied to the
two-dimensional version of the problem, 
for B-splines of order $p=2,3,4$, with minimal and maximal continuity,
with radical meshes and $1 \leq \gamma \leq 4.5$. 
As expected, the rate increases as $\gamma$ increases, until it reaches
a maximum (then decreases in some cases).
However, it is not easy to draw a general conclusion on the connection
between the choice of $\gamma$, the degree $p$ and the regularity of the 
B-splines basis. Remarkably, for $p=2,3$ it seems that
only B-splines with maximal regularity reach the optimal convergence rate.
This is not the case for $p=4$. 
A possible explanation is that even for larger values of $\gamma$,
the optimal set is not $\beta_1+\beta_2 = constant$ but it becomes a more 
pronounced ``convex shape'' as $p$ increases, regardless of the continuity
of the B-spline basis, see Figure \ref{fig:err-vs-time-non-reg-gamma3-ring-qoi}. 
After inspection of Figure \ref{fig:rates_vs_gamma-ring},
we now repeat the analysis of error vs time and degrees-of-freedom setting $\gamma=3$.

\begin{figure}[tp]
  \centering
  \includegraphics[width=\ConvSize\linewidth]{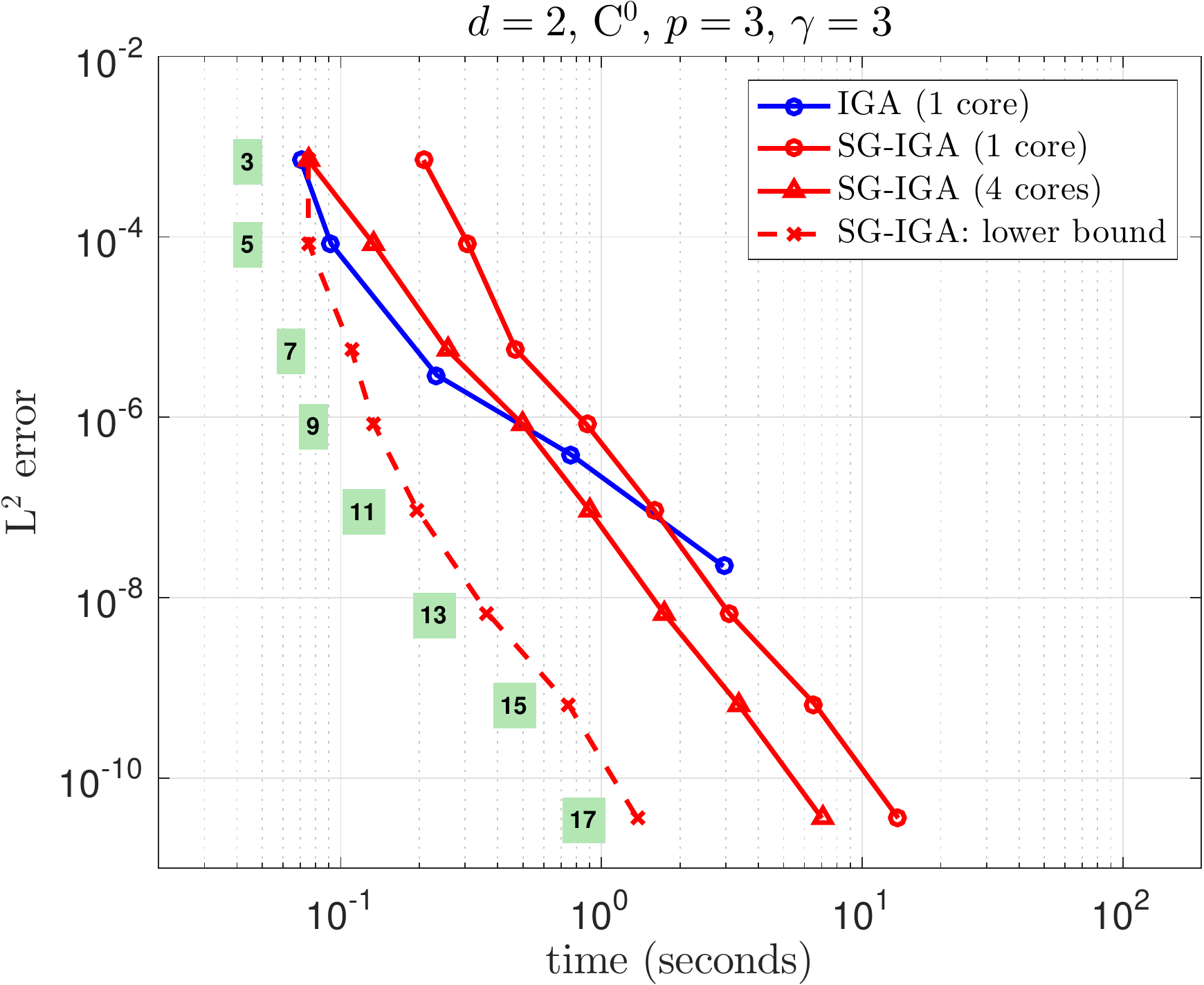}
  \includegraphics[width=\ConvSize\linewidth]{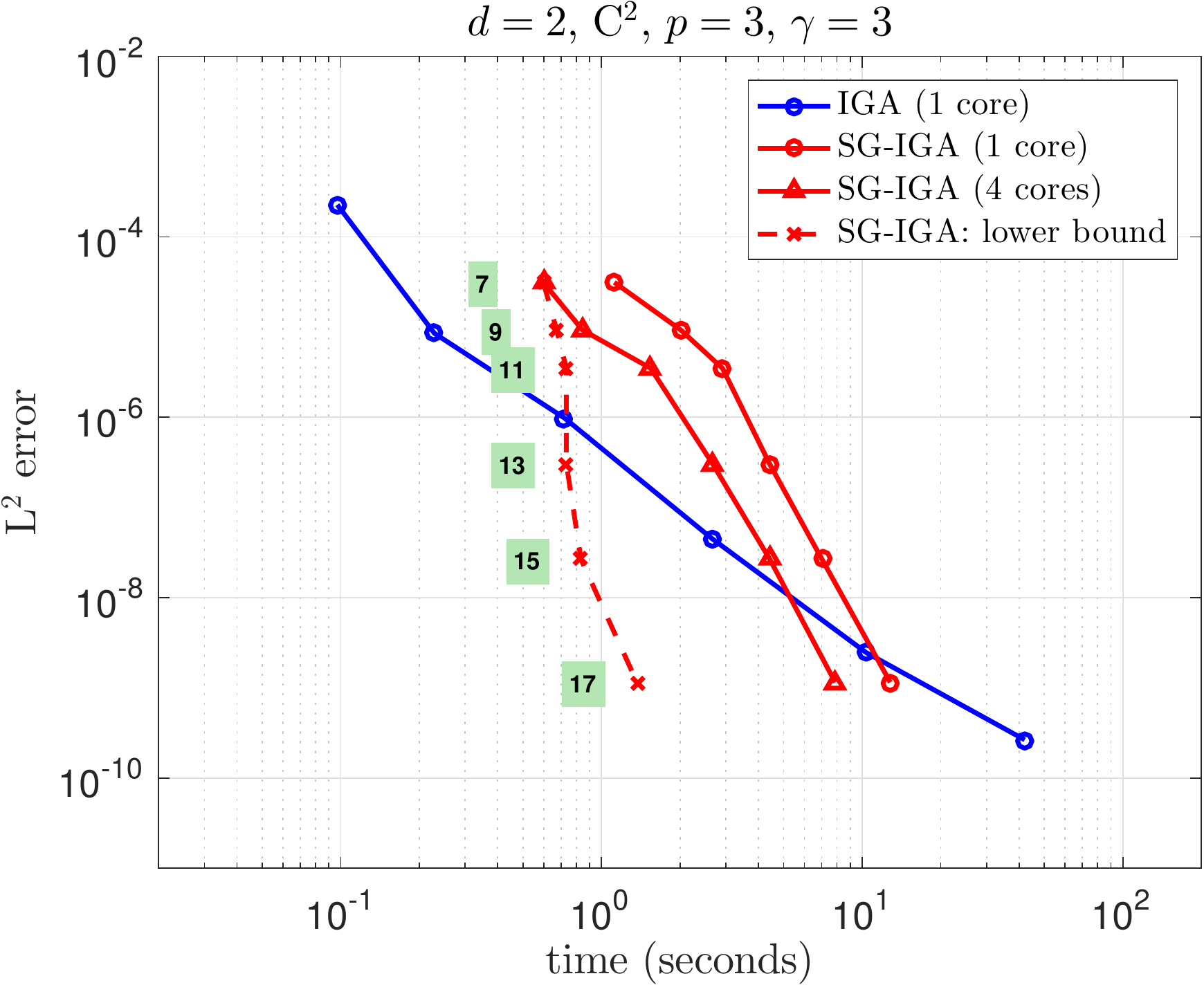} \\[\SpaceSize]
  \includegraphics[width=\ConvSize\linewidth]{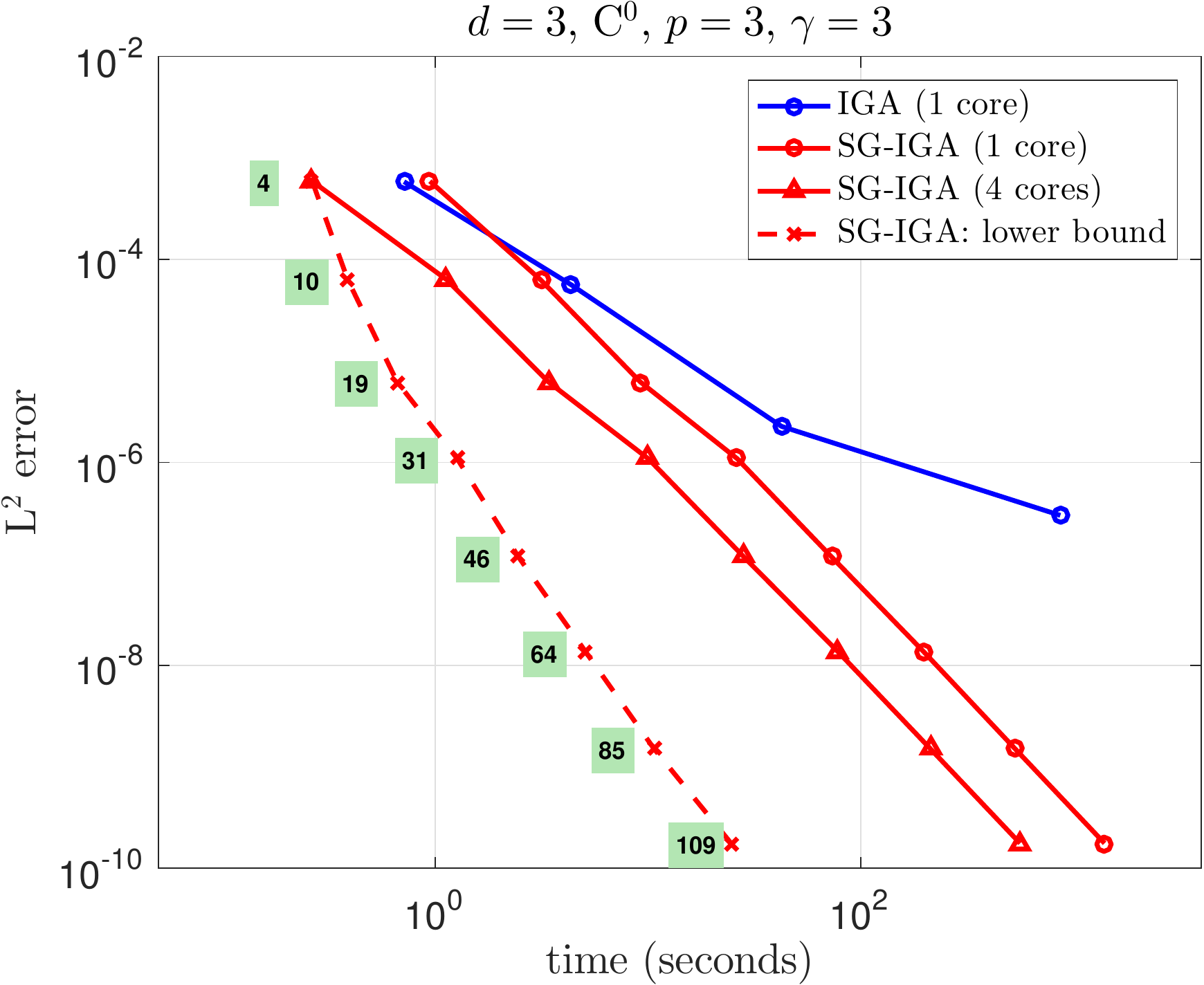}
  \includegraphics[width=\ConvSize\linewidth]{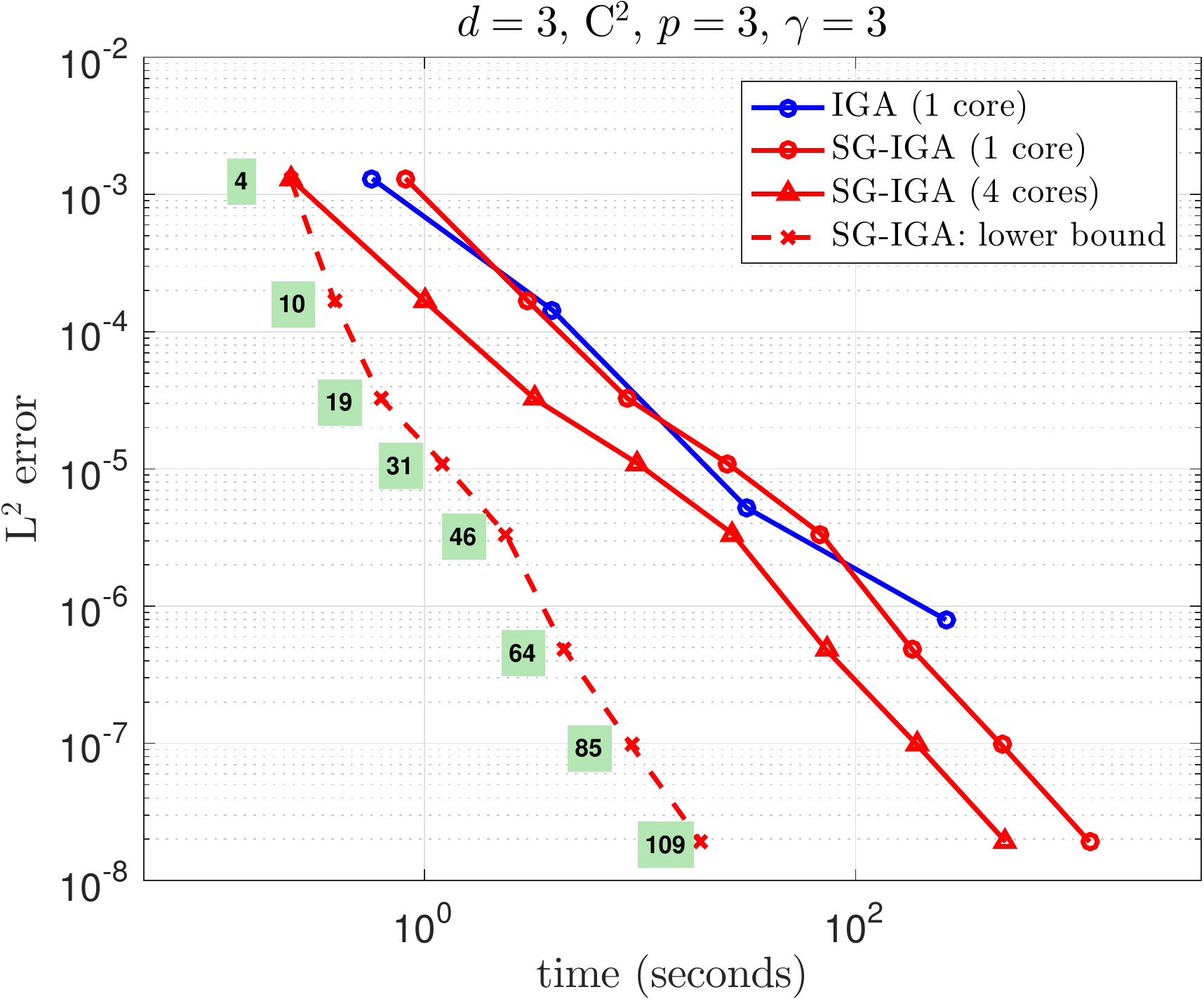}
  \caption{$L^2$ error vs. time for the quarter-of-annulus problem with low-regular solution,
    with radical meshes and $\gamma=3$.
    Here we fix $p=3$ and change $d$ (top row: $d=2$, bottom row: $d=3$)
    and the regularity of the B-splines basis
    \La{(left column: $C^0$, right column: $C^2$)}.
    The dashed line is the lower bound that 
    can be achieved if the number of available cores is at least equal to the number of components 
    of the combination technique for each level, given by the numbers in green boxes.}
  \label{fig:err-vs-time-non-reg-gamma3-ring}

  \bigskip

  \includegraphics[width=\ConvSize\linewidth]{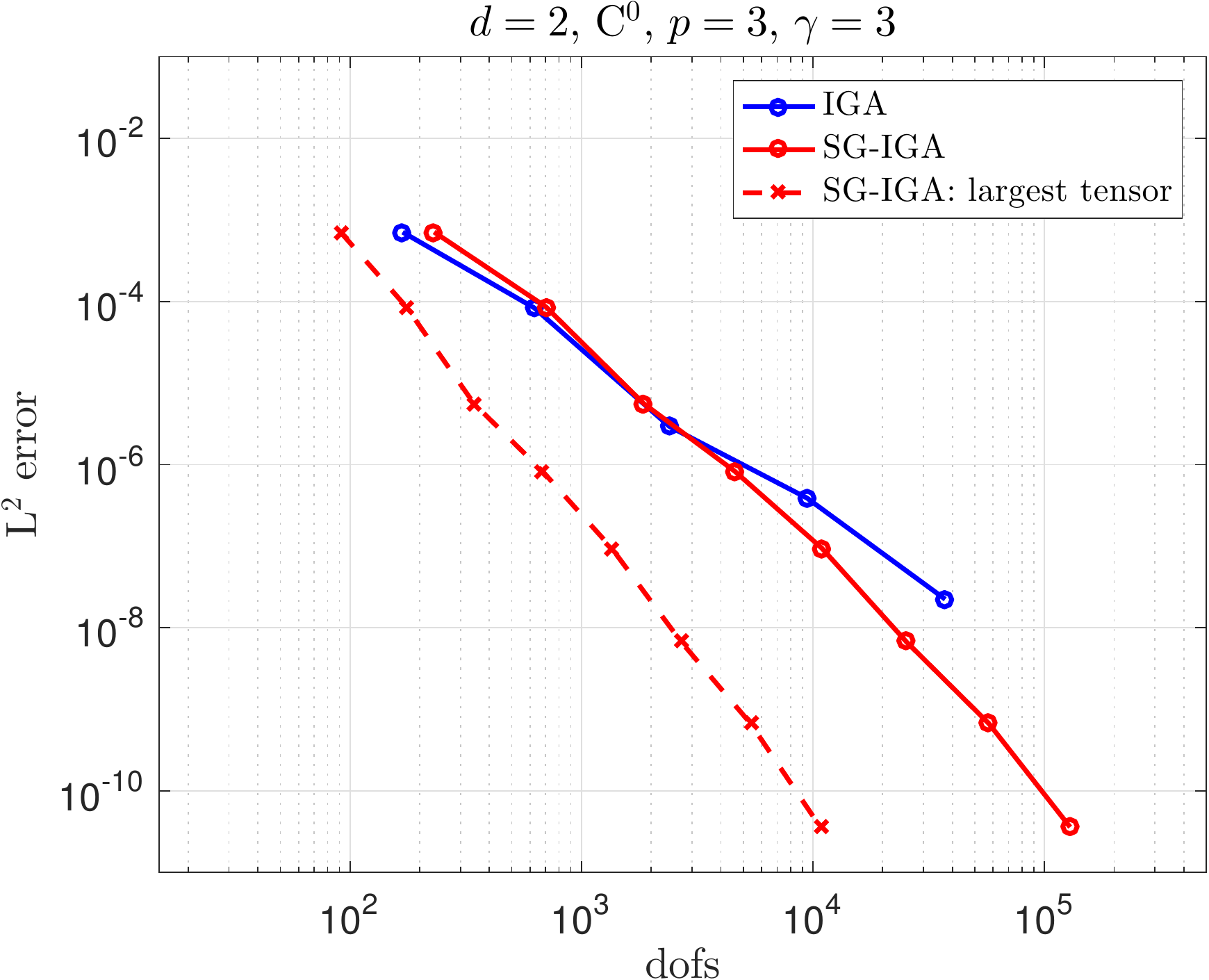}
  \includegraphics[width=\ConvSize\linewidth]{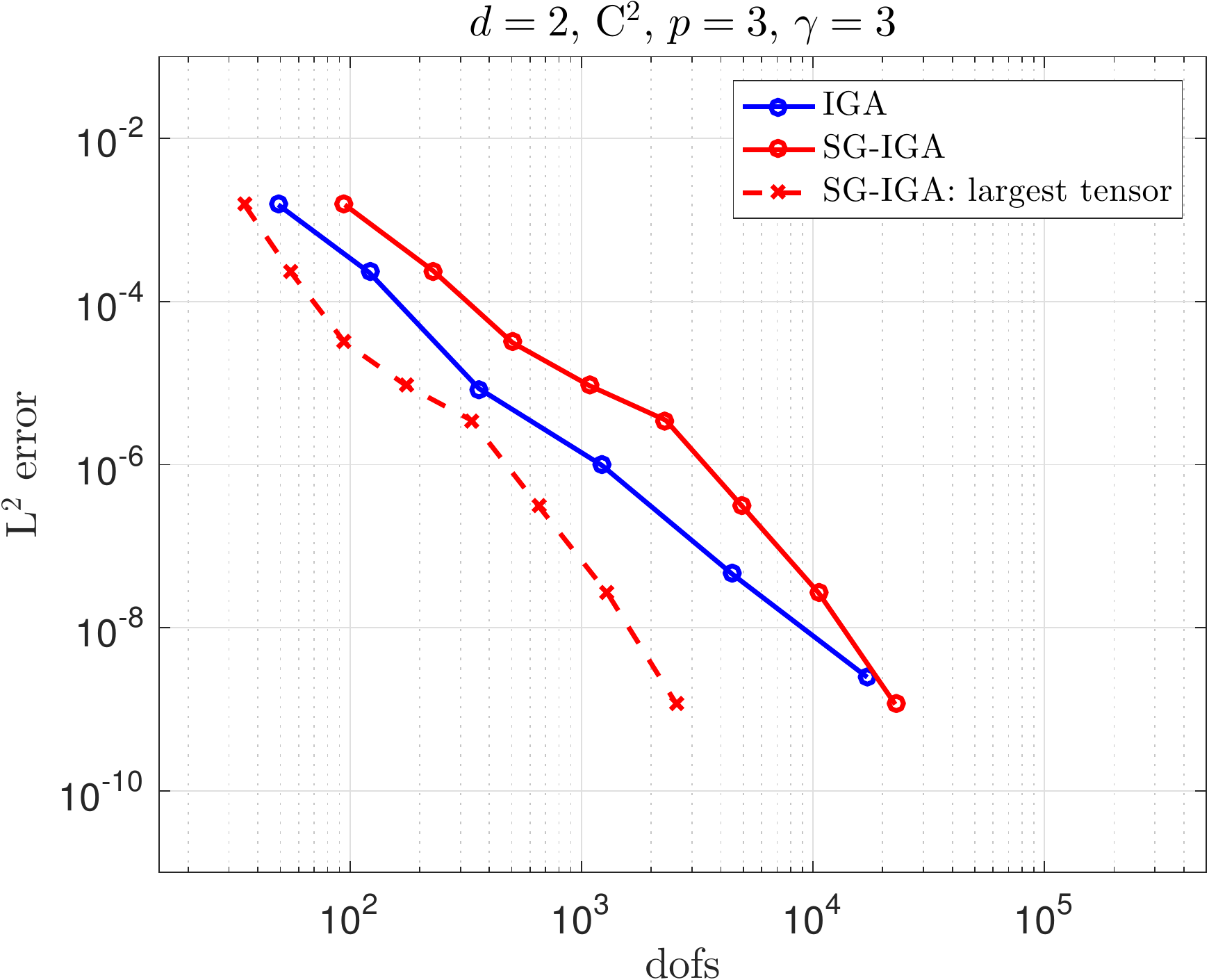} \\[\SpaceSize]
  \includegraphics[width=\ConvSize\linewidth]{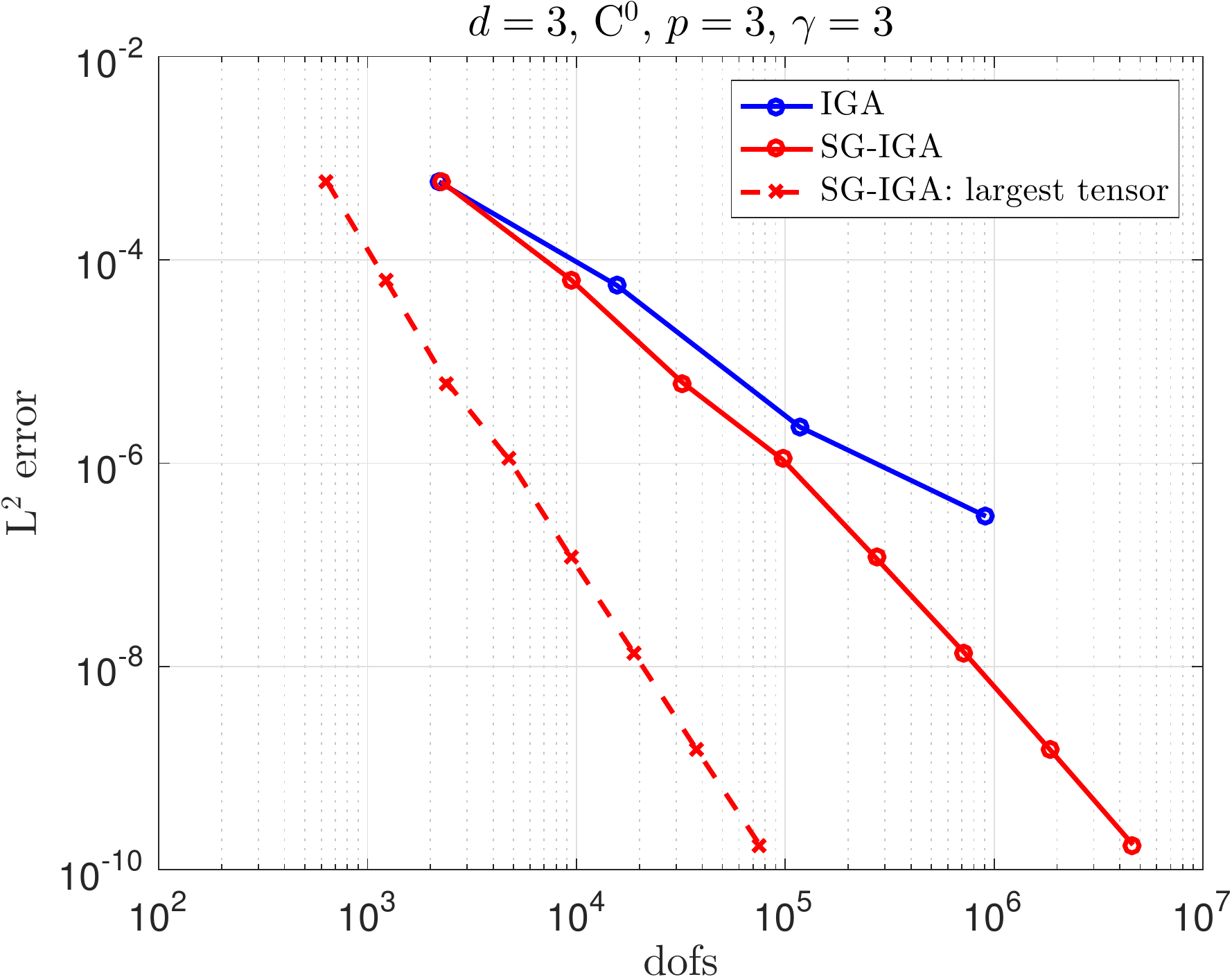}
  \includegraphics[width=\ConvSize\linewidth]{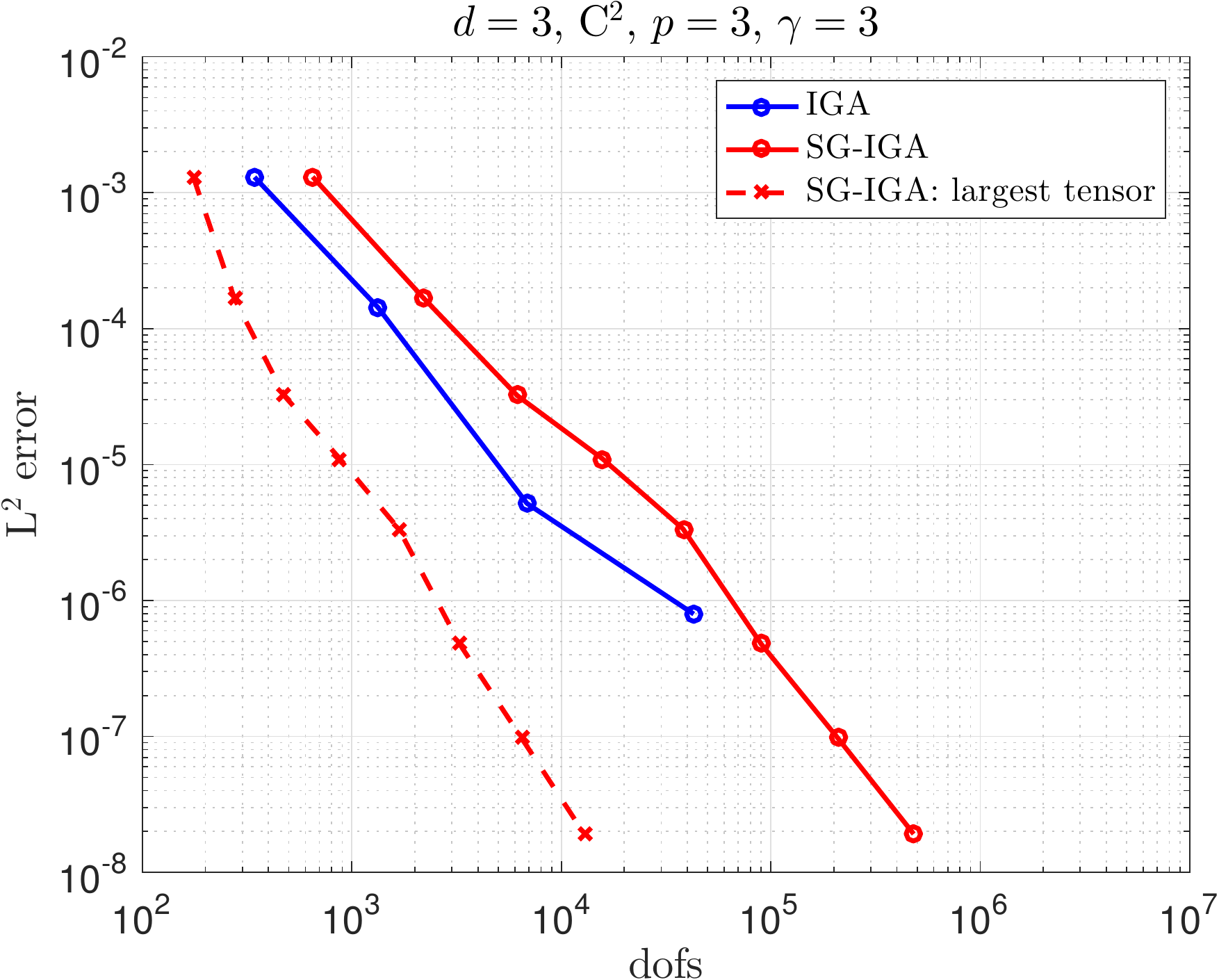}
  \caption{$L^2$ error vs. dofs for the quarter-of-annulus problem with low-regular solution,
    with radical meshes and $\gamma=3$.
    Here we fix $p=3$ and change $d$ (top row: $d=2$, bottom row: $d=3$)
    and the regularity of the B-splines basis \La{(left column: $C^0$, right column: $C^2$)}.
    In each figure, the dashed line represents the error of SG-IGA vs. the number
    of the degrees-of-freedom for the largest tensor grid in the combination technique.}\label{fig:err-vs-dofs-non-reg-gamma3-ring}
\end{figure}

The results obtained by introducing graded grids $\gamma=3$ are shown in
Figures \ref{fig:err-vs-time-non-reg-gamma3-ring} and
\ref{fig:err-vs-dofs-non-reg-gamma3-ring}, and should be 
compared to the results for non-radical meshes reported in Figures
\ref{fig:err-vs-time-non-reg-ring} and \ref{fig:err-vs-dofs-non-reg-ring}.
Inspecting the plots, it can be seen that SG-IGA is now 
competitive with IGA, especially for $d=3$, 
and indeed the results are now analogous to those reported in 
Figures \ref{fig:err-vs-time-reg-ring} and \ref{fig:err-vs-dofs-reg-ring}
for the problem with regular solution.


\subsection{Horseshoe domain}

\begin{figure}[tp]
  \centering
  \includegraphics[width=\ConvSize\textwidth]{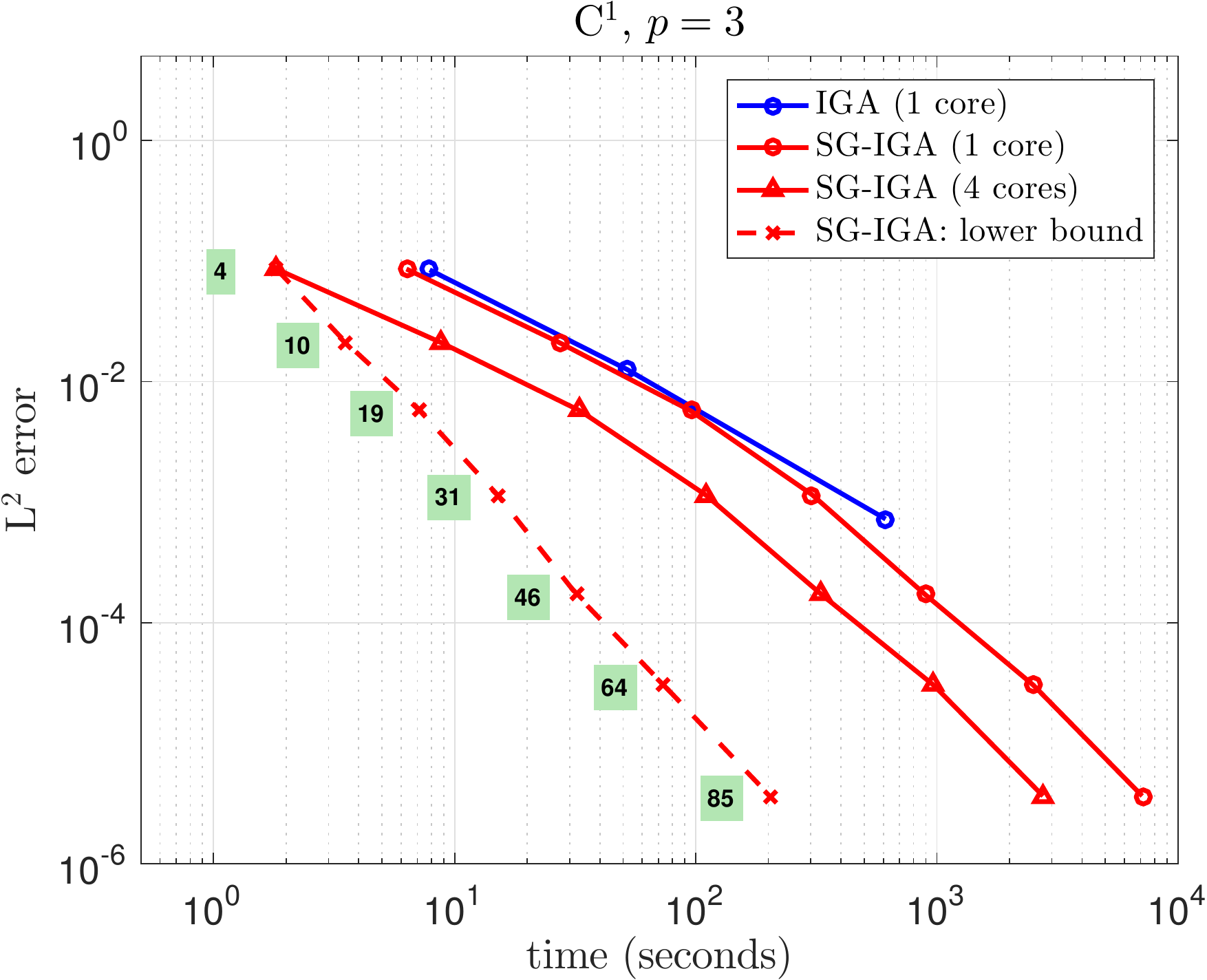}
  \includegraphics[width=\ConvSize\textwidth]{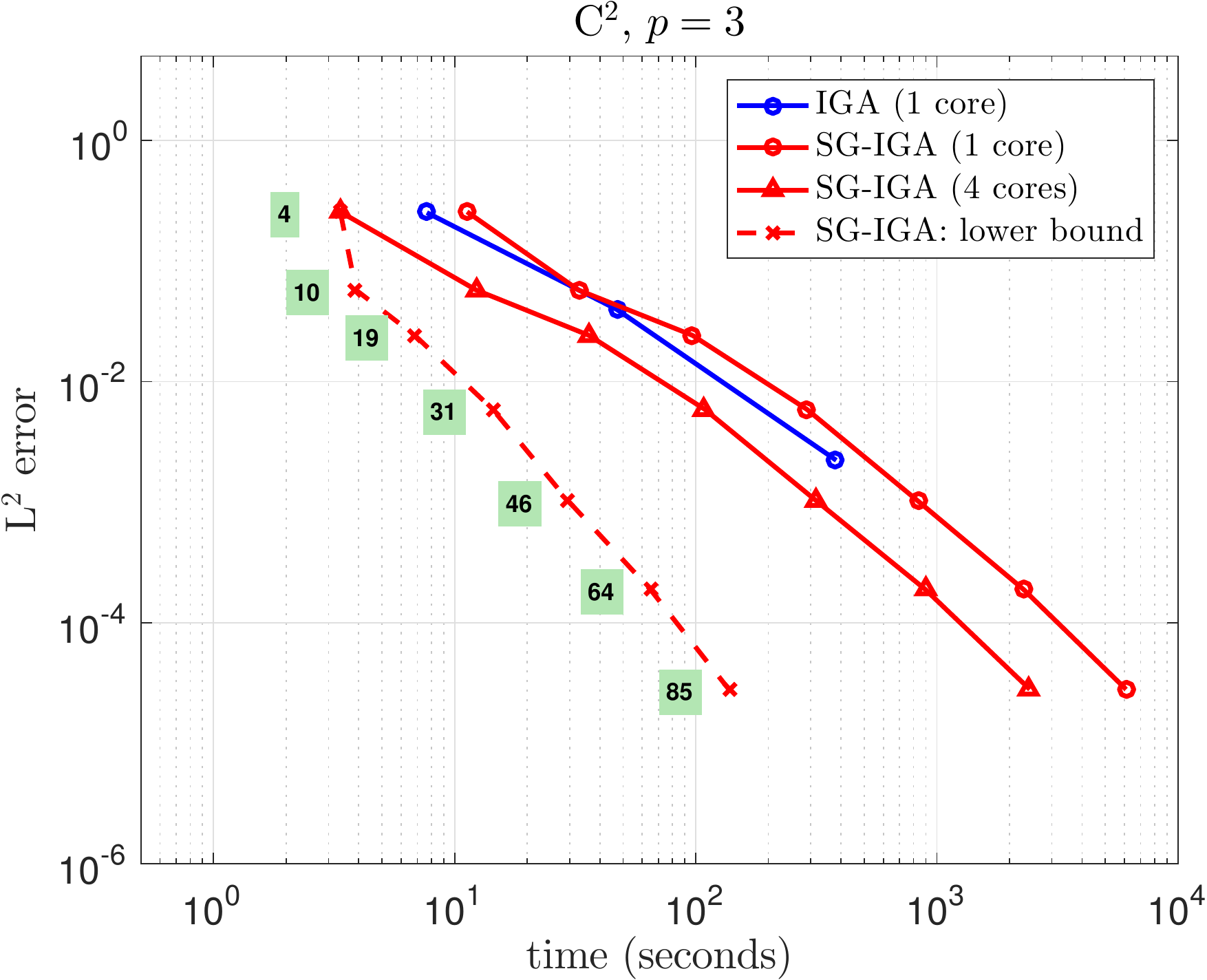} \\[\SpaceSize]
  \includegraphics[width=\ConvSize\textwidth]{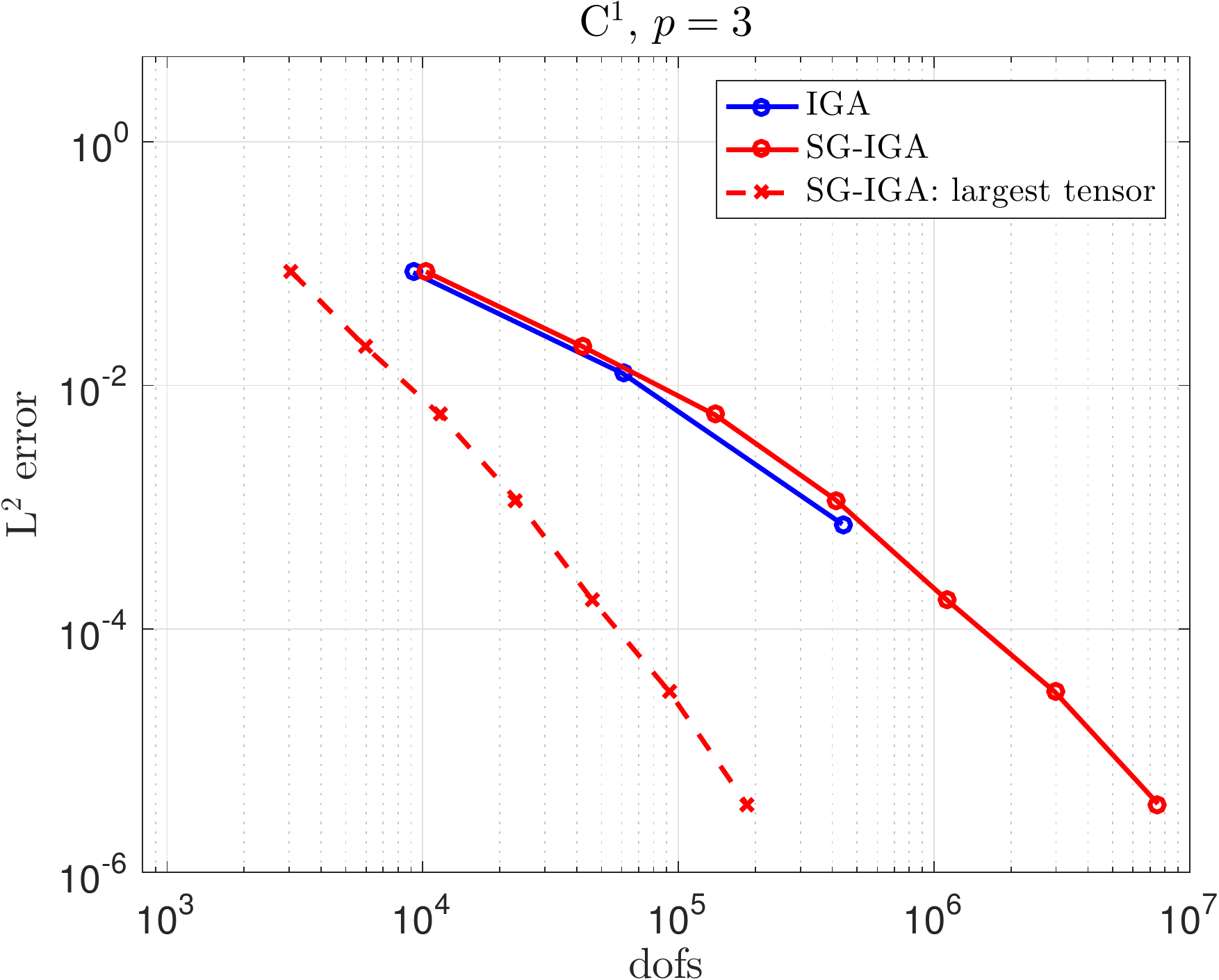}
  \includegraphics[width=\ConvSize\textwidth]{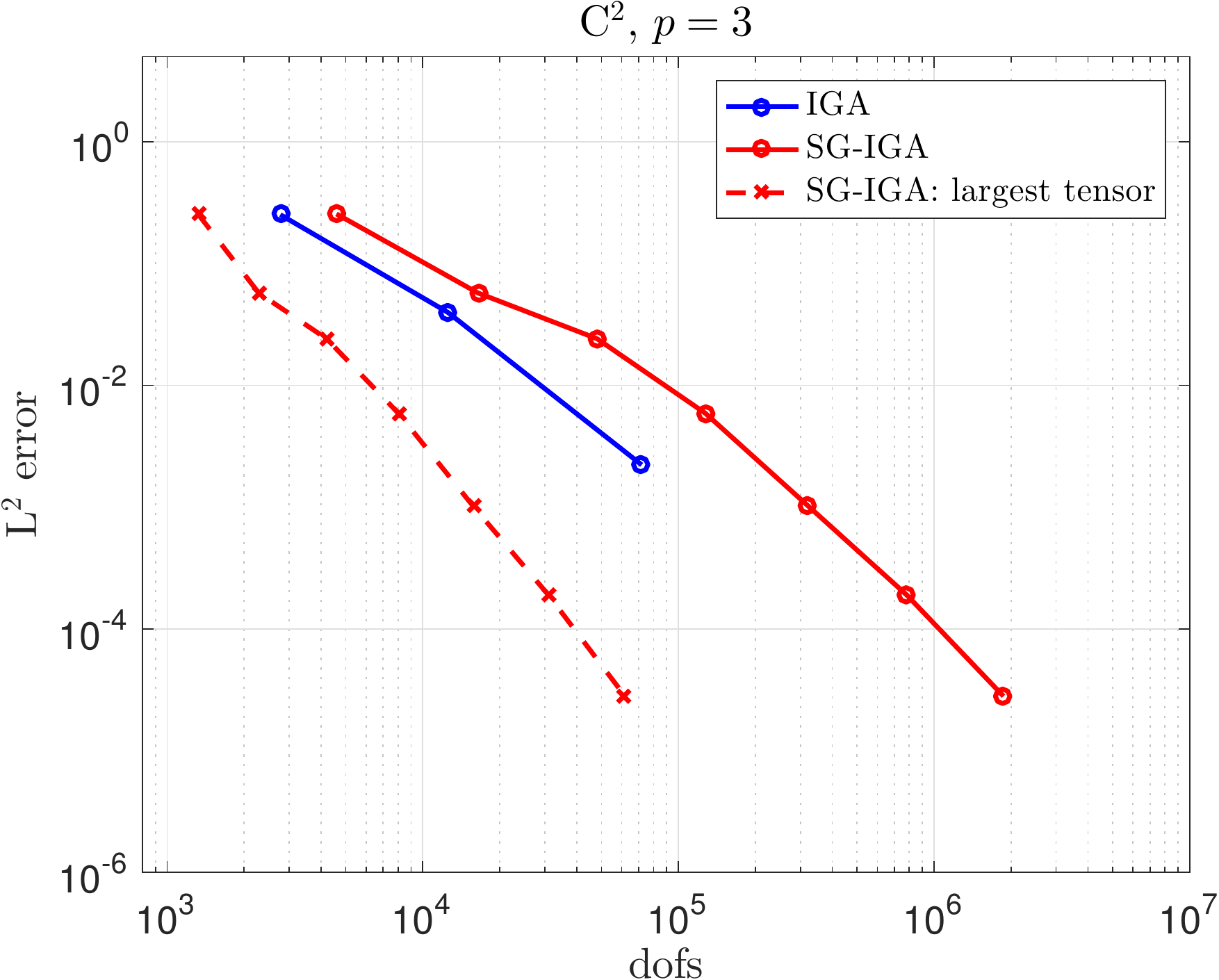} 
  \caption{$L^2$ error vs. time (top row) and degrees-of-freedom (bottom row) 
    for the horseshoe problem; \La{the left column shows results with $C^1$ B-splines, while
    the right column with $C^2$ B-splines.}
    For the bottom plots, the dashed line represents the error of SG-IGA vs. the number
    of the degrees-of-freedom for the largest tensor grid in the combination technique.
    Conversely, for the top plots the dashed line is the lower bound that 
    can be achieved if the number of available cores is at least equal to the number of components 
    of the combination technique for each level, given by the numbers in green boxes.
}
  \label{fig:err-vs-time-dofs-horseshoe}
\end{figure}

The last test we consider is on a more complex geometry, i.e., the 
horseshoe domain in Figure \ref{fig:domains}-right.
Also in this case, we set $f(\xx)=1$, which will result in a solution 
with reduced regularity. 
Therefore, based on the previous experience, we tackle this problem with a graded mesh with $\gamma=3$.
The results shown in Figure \ref{fig:err-vs-time-dofs-horseshoe} for
error versus time and degrees-of-freedom 
are similar to those for the three-dimensional version of quarter-of-annulus problem
of Figure \ref{fig:err-vs-time-non-reg-gamma3-ring}. 





\section{Conclusions}\label{sec:conclusions}

In this work we have shown how IGA solvers naturally fit in the sparse-grid construction framework,
obtaining a method that on the one hand can improve the performances of standard full-tensor IGA approximations
(if certain regularity requirements for the solution are met) and on the other hand extends
the spectrum of applications of the sparse-grid technology to arbitrary domains and to
basis functions with arbitrary order and regularity in a very convenient way.
Moreover, the combination-technique version of the sparse-grid approximation allows to 
maximize the reuse of pre-existing IGA solver in a black-box fashion, and provides the user
with a simple yet effective parallelization strategy for serial solvers.
For problems whose solution does not feature the needed regularity for sparse grids to 
converge in a satisfactory way, a remedy consists in introducing properly tuned
radical meshes in the parametric domain, even if more investigations are needed
to explore the interplay between $\gamma$, $p$, and regularity of the basis.

This preliminary work opens to possible extensions 
and research directions. On the ``methodological'' level, it would be interesting
to consider $p$ and $h-p$ sparsification for IGA solvers, as well as local adaptivity
``\`a la sparse-grid'', see e.g. \cite{b.griebel:acta,jakeman.archi.xiu:discont,Pfluger:2013}, 
which might be borrowed and interwined with the ``classical'' IGA refinement strategies such as
\cite{Giannelli2012485,LR-splines,GJS13}. 
Other interesting work directions are the sparsified version of IGA collocation solvers, and
the combination of fast IGA implementations \cite{sangalli2016isogeometric,calabro2016fast}
with the sparse-grid technology. Moreover, it might be worthwhile investigating sparse-adaptive solvers, 
in the spirit of what proposed while discussing Figures \ref{fig:opt-sets-ring-reg} and \ref{fig:opt-sets-ring-nonreg}.
Let us also remark that, even if the combination technique sometimes does not yield
impressive gains over the regular full-tensor solver, it is nonetheless a fast way to build
an approximation of the solution, and in this spirit it could be interesting to use it
as a preconditioner for complex problems.



\section*{Acknowledgement}
Giancarlo Sangalli and Lorenzo Tamellini were partially supported by the European Research Council
through the FP7 ERC consolidator grant n.616563 \emph{HIGEOM}
and by the GNCS 2017 project 
``Simulazione numerica di problemi di Interazione Fluido-Struttura (FSI) con metodi agli elementi finiti ed isogeometrici''.
Lorenzo Tamellini also received support from the scholarship 
``Isogeometric method'' granted by the Universit\`a di Pavia and by 
the European Union's Horizon 2020 research and innovation program 
through the grant no. 680448 CAxMan. 
Joakim Beck received support from the KAUST CRG3 Award Ref:2281 and the KAUST CRG4 Award Ref:2584.

\section*{References}
\bibliographystyle{elsarticle-num}
\bibliography{sparse_grids_biblio,IGA_biblio_gian,UQ_biblio}

\def\cprime{$'$}
\begin{thebibliography}{10}
\expandafter\ifx\csname url\endcsname\relax
  \def\url#1{\texttt{#1}}\fi
\expandafter\ifx\csname urlprefix\endcsname\relax\def\urlprefix{URL }\fi
\expandafter\ifx\csname href\endcsname\relax
  \def\href#1#2{#2} \def\path#1{#1}\fi

\bibitem{Hughes:2005}
T.~J.~R. Hughes, J.~A. Cottrell, Y.~Bazilevs, Isogeometric analysis: {CAD},
  finite elements, {NURBS}, exact geometry and mesh refinement, Computer
  Methods in Applied Mechanics and Engineering 194~(39) (2005) 4135--4195.

\bibitem{IGA-book}
J.~A. Cottrell, T.~J.~R. Hughes, Y.~Bazilevs, Isogeometric {A}nalysis: toward
  integration of {CAD} and {FEA}, John Wiley \& Sons, 2009.

\bibitem{acta-IGA}
L.~Beirao Da~Veiga, A.~Buffa, G.~Sangalli, R.~V\'{a}zquez, Mathematical
  analysis of variational isogeometric methods, Acta Numerica 23 (2014)
  157--287.

\bibitem{Gao2013}
L.~Gao, Kronecker products on preconditioning, Ph.D. thesis, King Abdullah
  University of Science and Technology (2013).

\bibitem{Gao2014}
L.~Gao, V.~M. Calo, Fast isogeometric solvers for explicit dynamics, Computer
  Methods in Applied Mechanics and Engineering 274 (2014) 19--41.

\bibitem{Gao2015}
L.~Gao, V.~M. Calo, Preconditioners based on the alternating-direction-implicit
  algorithm for the 2d steady-state diffusion equation with orthotropic
  heterogeneous coefficients, Journal of Computational and Applied Mathematics
  273 (2015) 274--295.

\bibitem{sangalli2016isogeometric}
G.~Sangalli, M.~Tani, {Isogeometric preconditioners based on fast solvers for
  the Sylvester equation}, SIAM Journal on Scientific Computing 38~(6) (2016)
  A3644--A3671.

\bibitem{calabro2016fast}
F.~Calabr\`o, G.~Sangalli, M.~Tani, Fast formation of isogeometric galerkin
  matrices by weighted quadrature, Computer Methods in Applied Mechanics and
  Engineering.

\bibitem{mantzaflaris2016low}
A.~Mantzaflaris, B.~J{\"u}ttler, B.~N. Khoromskij, U.~Langer, Low rank tensor
  methods in galerkin-based isogeometric analysis, Computer Methods in Applied
  Mechanics and Engineering 316 (2017) 1062--1085.

\bibitem{antolin2015efficient}
P.~Antolin, A.~Buffa, F.~Calabro, M.~Martinelli, G.~Sangalli, Efficient matrix
  computation for tensor-product isogeometric analysis: {T}he use of sum
  factorization, Computer Methods in Applied Mechanics and Engineering 285
  (2015) 817--828.

\bibitem{donatelli2017symbol}
M.~Donatelli, C.~Garoni, C.~Manni, S.~Serra-Capizzano, H.~Speleers,
  Symbol-based multigrid methods for {G}alerkin {B}-spline isogeometric
  analysis, SIAM Journal on Numerical Analysis 55~(1) (2017) 31--62.

\bibitem{hofreither2015robust}
C.~Hofreither, S.~Takacs, W.~Zulehner, A robust multigrid method for
  isogeometric analysis in two dimensions using boundary correction, Computer
  Methods in Applied Mechanics and Engineering 316 (2017) 22 -- 42, special
  Issue on Isogeometric Analysis: Progress and Challenges.

\bibitem{zenger91sparse}
C.~Zenger, Sparse grids, in: W.~Hackbusch (Ed.), Parallel Algorithms for
  Partial Differential Equations, Vol.~31 of Notes on Numerical Fluid
  Mechanics, Vieweg, 1991, pp. 241--251.

\bibitem{bungartz:PhD.thesis}
H.~Bungartz, D\"unne gitter und deren anwendung bei der adaptiven l\"osung der
  dreidimensionalen poisson-gleichung, Ph.D. thesis, Institut f\"ur Informatik,
  Technische Universit\"at Munchen (1992).

\bibitem{griebel:first}
M.~Griebel, A parallelizable and vectorizable multi-level-algorithm on sparse
  grids, in: W.~Hackbusch (Ed.), Parallel Algorithms for Partial Differential
  Equations. Notes on Numerical Fluid Mechanics, Vol.~31, Verlag Vieweg,
  Braunschweig, 1991, pp. 94--100.

\bibitem{Griebel.schneider.zenger:combination}
M.~Griebel, M.~Schneider, C.~Zenger, A combination technique for the solution
  of sparse grid problems, in: P.~de~Groen, R.~Beauwens (Eds.), {Iterative
  Methods in Linear Algebra}, IMACS, Elsevier, North Holland, 1992, pp.
  263--281.

\bibitem{b.griebel:acta}
H.~Bungartz, M.~Griebel, Sparse grids, Acta Numer. 13 (2004) 147--269.

\bibitem{Garcke2013}
J.~Garcke, Sparse grids in a nutshell, in: J.~Garcke, M.~Griebel (Eds.), Sparse
  Grids and Applications, Springer Berlin Heidelberg, Berlin, Heidelberg, 2013,
  pp. 57--80.

\bibitem{Dornseifer:curvilinear.sparse.grids}
T.~Dornseifer, C.~Pflaum, Discretization of elliptic differential equations on
  curvilinear bounded domains with sparse grids, Computing 56~(3) (1996)
  197--213.

\bibitem{bungartz:general-deg-and-dom}
H.~J. Bungartz, T.~Dornseifer, Sparse grids: {R}ecent developments for elliptic
  partial differential equations, in: W.~Hackbusch, G.~Wittum (Eds.), Multigrid
  Methods V, Vol.~3 of Lecture Notes in Computational Science and Engineering,
  Springer, Berlin/Heidelberg, 1998.

\bibitem{achatz:high.order.sparse.grids}
S.~Achatz, Higher order sparse grid methods for elliptic partial differential
  equations with variable coefficients, Computing 71~(1) (2003) 1--15.

\bibitem{Babuska_book}
I.~Babu{\v{s}}ka, T.~Strouboulis, The finite element method and its
  reliability, Numerical Mathematics and Scientific Computation, The Clarendon
  Press Oxford University Press, New York, 2001.

\bibitem{Beirao_Cho_Sangalli}
L.~Beir{\~a}o~da Veiga, D.~Cho, G.~Sangalli, Anisotropic {NURBS} approximation
  in {I}sogeometric {A}nalysis, Comput. Methods Appl. Mech. Engrg. 209-212
  (2012) 1--11.

\bibitem{garcke:graded}
J.~Garcke, M.~Griebel, On the computation of the eigenproblems of hydrogen and
  helium in strong magnetic and electric fields with the sparse grid
  combination technique, Journal of Computational Physics 165~(2) (2000) 694 --
  716.

\bibitem{griebel.thurner:graded}
M.~Griebel, V.~Thurner, {The efficient solution of fluid dynamics problems by
  the combination technique}, International Journal of Numerical Methods for
  Heat \& Fluid Flow 5~(3) (1995) 251--269.

\bibitem{Grisvard}
P.~Grisvard, Elliptic problems in nonsmooth domains, Vol.~24 of Monographs and
  Studies in Mathematics, Pitman (Advanced Publishing Program), Boston, MA,
  1985.

\bibitem{BeiraoDaVeiga:2014}
L.~Beira{\~o}~da Veiga, A.~Buffa, G.~Sangalli, R.~V{\'a}zquez, Mathematical
  analysis of variational isogeometric methods, Acta Numerica 23 (2014)
  157--287.

\bibitem{DeBoor}
C.~de~Boor, A practical guide to splines, revised Edition, Vol.~27 of Applied
  Mathematical Sciences, Springer-Verlag, New York, 2001.

\bibitem{wasi.wozniak:cost.bounds}
G.~Wasilkowski, H.~Wozniakowski, Explicit cost bounds of algorithms for
  multivariate tensor product problems, Journal of Complexity 11~(1) (1995) 1
  -- 56.

\bibitem{Bungartz.Griebel.Roschke.ea:pointwise.conv}
H.-J. Bungartz, M.~Griebel, D.~R\"oschke, C.~Zenger, Pointwise convergence of
  the combination technique for the {L}aplace equation, East-West J. Numer.
  Math. 2 (1994) 21--45.

\bibitem{Hegland:combination}
M.~Hegland, J.~Garcke, V.~Challis, The combination technique and some
  generalisations, Linear Algebra and its Applications 420~(2–3) (2007) 249
  -- 275.

\bibitem{griebel.harbrecht:combi-conv}
M.~Griebel, H.~Harbrecht, On the convergence of the combination technique, in:
  J.~Garcke, D.~Pfl\"uger (Eds.), Sparse Grids and Applications - Munich 2012,
  Vol.~97 of Lecture Notes in Computational Science and Engineering, Springer
  International Publishing, 2014, pp. 55--74.

\bibitem{reisinger:analysis}
C.~Reisinger, Analysis of linear difference schemes in the sparse grid
  combination technique, IMA Journal of Numerical Analysis 33~(2) (2013)
  544--581.

\bibitem{reisinger:efficient}
C.~Reisinger, G.~Wittum, Efficient hierarchical approximation of
  high-dimensional option pricing problems, SIAM Journal on Scientific
  Computing 29~(1) (2007) 440--458.

\bibitem{griebel.harbrecht:tensor-spaces}
M.~Griebel, H.~Harbrecht, On the construction of sparse tensor product spaces,
  Mathematics of Computation 82~(282) (2013) 975--994.

\bibitem{griebel.harbrecht:L-fold}
M.~Griebel, H.~Harbrecht, A note on the construction of {L}-fold sparse tensor
  product spaces, Constructive Approximation 38~(2) (2013) 235--251.

\bibitem{pflaum-zhou:combi_conv}
C.~Pflaum, A.~Zhou, Error analysis of the combination technique, Numerische
  Mathematik 84~(2) (1999) 327--350.

\bibitem{Bungartz:two.proofs}
H.-J. Bungartz, M.~Griebel, D.~Röschke, C.~Zenger, Two proofs of convergence
  for the combination technique for the efficient solution of sparse grid
  problems, in: In Domain decomposition methods in scientific and engineering
  computing (University Park, PA), 1993.

\bibitem{geopdesv3}
R.~V\'azquez, A new design for the implementation of isogeometric analysis in
  {O}ctave and {M}atlab: {G}eo{PDE}s 3.0, Computer \& Mathematics with
  Applications 72~(3) (2016) 523--554.

\bibitem{sangalli2017matrix}
G.~Sangalli, M.~Tani, Matrix-free isogeometric analysis: the computationally
  efficient $ k $-method, arXiv preprint arXiv:1712.08565.

\bibitem{griebel:parallerl92}
M.~Griebel, The combination technique for the sparse grid solution of {PDE}s on
  multiprocessor machine, Parallel Processing Letters 2 (1992) 61--70.

\bibitem{heene:efficient}
M.~Heene, D.~Pfl\"uger, Efficient and scalable distributed-memory
  hierarchization algorithms for the sparse grid combination technique, in:
  Parallel Computing: On the Road to Exascale, Vol.~27 of Advances in Parallel
  Computing, IOS Press, 2016, pp. 339--348.

\bibitem{Heene:scalable}
M.~Heene, D.~Pfl{\"u}ger, Scalable algorithms for the solution of
  higher-dimensional {PDE}s, in: H.-J. Bungartz, P.~Neumann, W.~E. Nagel
  (Eds.), Software for Exascale Computing - SPPEXA 2013-2015, Springer
  International Publishing, Cham, 2016, pp. 165--186.

\bibitem{garcke_hegland_nielsen_2006}
J.~Garcke, M.~Hegland, O.~Nielsen, Parallelisation of sparse grids for large
  scale data analysis, The ANZIAM Journal 48~(1) (2006) 11–22.

\bibitem{martello:knapsack}
S.~Martello, P.~Toth, Knapsack problems: algorithms and computer
  implementations, Wiley-Interscience series in discrete mathematics and
  optimization, J. Wiley \& Sons, 1990.

\bibitem{griebel.knapek:optimized}
M.~Griebel, S.~Knapek, Optimized general sparse grid approximation spaces for
  operator equations, Math. Comp. 78~(268) (2009) 2223--2257.

\bibitem{nobile.eal:optimal-sparse-grids}
F.~Nobile, L.~Tamellini, R.~Tempone, Convergence of quasi-optimal sparse-grid
  approximation of {H}ilbert-space-valued functions: application to random
  elliptic {PDE}s, Numerische {M}athematik 134~(2) (2016) 343--388.

\bibitem{dauge2006elliptic}
M.~Dauge, {Elliptic Boundary Value Problems on Corner Domains: Smoothness and
  Asymptotics of Solutions}, Lecture Notes in Mathematics, Springer Berlin
  Heidelberg, 2006.

\bibitem{jakeman.archi.xiu:discont}
J.~Jakeman, R.~Archibald, D.~Xiu, Characterization of discontinuities in
  high-dimensional stochastic problems on adaptive sparse grids, J. Comput.
  Phys. 230~(10).

\bibitem{Pfluger:2013}
D.~Pfl\"uger, Spatially adaptive refinement, in: J.~Garcke, M.~Griebel (Eds.),
  Sparse Grids and Applications, Springer Berlin Heidelberg, Berlin,
  Heidelberg, 2013, pp. 243--262.

\bibitem{Giannelli2012485}
C.~Giannelli, B.~J{\"u}ttler, H.~Speleers, {THB}-splines: The truncated basis
  for hierarchical splines, Comput. Aided Geom. Design. 29~(7) (2012) 485 --
  498.

\bibitem{LR-splines}
T.~Dokken, T.~Lyche, K.~F. Pettersen, Polynomial splines over locally refined
  box-partitions, Comput. Aided Geom. Design 30~(3) (2013) 331--356.

\bibitem{GJS13}
C.~Giannelli, B.~J{\"u}ttler, H.~Speleers, Strongly stable bases for adaptively
  refined multilevel spline spaces, Adv. Comput. Math. (2013) 1--32.

\end{thebibliography}



\end{document}